\documentclass[12pt,letterpaper]{article}
\usepackage{fullpage}

\usepackage[dvips]{graphics, epsfig, color}

\usepackage{graphics}
\usepackage{graphicx}
\usepackage{amsmath, amsthm}
\usepackage{amsfonts}

\newtheorem{thm}{Theorem}[section]
\newtheorem{cor}[thm]{Corollary}
\newtheorem{lem}[thm]{Lemma}
\newtheorem{prop}[thm]{Proposition}
\theoremstyle{fact}
\newtheorem{fct}[thm]{Fact}
\theoremstyle{definition}
\newtheorem{defn}[thm]{Definition}
\theoremstyle{remark}
\newtheorem{rem}[thm]{Remark}
\theoremstyle{example}

\numberwithin{equation}{subsection}

\newcommand{\reals}{\mathbb{R}}
\renewcommand{\Re}{\reals}
\newcommand{\rsym}{\mathbb{S}}
\newcommand{\trace}{\textup{\textrm{trace}}}
\newcommand{\closure}{\textup{\textrm{cl}}\:}

\newcommand{\proj}{\textup{\textrm{proj}}}
\newcommand{\Diag}{\textup{\textrm{Diag}}}
\newcommand{\nul}{\textup{\textrm{null}}}

\newcommand{\rank}{\textup{\textrm{rank}}}
\newcommand{\conv}{\textup{\textrm{conv}}}

\title{
Shrink-Wrapping trajectories for Linear Programming}
\author{Yuriy Zinchenko
\footnote{Department of Mathematics and Statistics, University of
Calgary, MS 446, 2500 University Drive NW, Calgary AB T2N 1N4,
Canada, e-mail: yzinchen@ucalgary.ca}}

\begin{document}
\maketitle \abstract{Hyperbolic Programming (HP) --minimizing a
linear functional over an affine subspace of a finite-dimensional
real vector space intersected with the so-called hyperbolicity
cone-- is a class of convex optimization problems that contains
well-known Linear Programming (LP). In particular, for any LP one
can readily provide a sequence of HP relaxations. Based on these
hyperbolic relaxations, a new Shrink-Wrapping approach to solve LP
has been proposed by Renegar. The resulting Shrink-Wrapping
trajectories, in a sense, generalize the notion of central path in
interior-point methods.

We study the geometry of Shrink-Wrapping trajectories for Linear
Programming. In particular, we analyze the geometry of these
trajectories in the proximity of the so-called central line, and
contrast the behavior of these trajectories with that of the central
path for some pathological LP instances.

In addition, we provide an elementary real proof of convexity of hyperbolicity cones.}

\section{Introduction}
We consider LP in its standard form
\begin{equation*}
\min_x\{c^T x:\: Ax=b, x\in\Re^n_+\}
\end{equation*}
where $c\in\Re^n, b\in\Re^m, A\in\Re^{m\times n}, m<n$, and
$\Re^n_+$ denotes $n$-dimensional nonnegative orthant. LP is
paramount in many applications of mathematical programming today.

Amongst numerical methods employed to solve LP instances in
practice, the two most notable classes are the so-called pivot-type
methods and the interior-point methods. For a given method, of
particular theoretical and practical interest is dependence of
number of iterations (and elementary arithmetic operations) required
to solve an LP. Assuming that the LP data, namely, the triple
$\{A,b,c\}$ are rational, one may measure its bit input complexity
as the number of $\{0,1\}$-bits required to store the data. If any
LP instance may be solved by the method in at most polynomial number
of arithmetic operations in $m,n$ and $L$, such a method is called
polynomial time algorithm for LP.

Pivot-type methods, such as simplex method, follow the faces of the
polytope that define the problem's feasible region $\{x: Ax=b,
x\in\Re^n_+\}$. Although, some variants of pivot-type methods
require polynomial number of iterations on average, and perform well
in practice, it is not known wether there exists a polynomial time
algorithm within this class.

In contrast, the interior-point methods follow some continuous
trajectory typically inside the feasible region; many variants of
these methods are known to be polynomial time algorithms and perform
well in practice especially for very large (in terms of $m$ and $n$)
LP instances.

A variant of the question of wether there exists a
strongly-polynomial algorithm for solving LP --in simple terms, the
algorithm whose running time would depend only on $m$ and $n$-- was
cited by S. Smale as on of the 18 greatest unsolved problems of the
21st century, and ``is the main unsolved problem of linear
programming theory.''

Linear Programming, when viewed from the point of view of the
so-called hyperbolic polynomials, exhibits rich and beautiful
algebraic structure, which seems not to be exploited yet by any of
the existing methods. It is our motivation to take advantage of this
structure in hope to develop potentially more efficient methods to
solve this optimization problem.

The extensive study of hyperbolic polynomials begins with the work
of Lars G\aa rding \cite{Garding:Poly}, which dates back to 1950's,
in the context of partial-differential equations; here the author
established a number of important results about the hyperbolic
polynomials including the convexity of the associated hyperbolicity
cones. The notion of hyperbolic programming was first introduced
in~\cite{Guler:Poly}; here the author demonstrated, in particular,
that the hyperbolic programming problems can be efficiently solved
using the interior point methods, and gave a first characterization
of the hyperbolicity cones as a set of polynomial inequalities
(although, quite different and more complicated then the one
in~\cite{Renegar:Poly} that we partially rely on). Further study of
hyperbolic polynomials in the context of convex optimization was
done by the group of authors of~\cite{Lewis:Poly}; a number of
important observations were made regarding the connections of
hyperbolic polynomials with the symmetric functions, and in
particular, the elementary symmetric functions. This latter
reference is an excellent introduction to hyperbolic polynomials in
the context of mathematical programming. This line of research was
continued in~\cite{Renegar:Poly}, where many important properties of
the boundary of the hyperbolicity cones are revealed together with
the relevance of the so-called hyperbolic derivative cones. In the
new paper~\cite{Renegar:SW} generalized trajectories for solving
hyperbolic programming problems based on hyperbolic relaxations are
introduced.

We study the so-called Shrink-Wrapping algorithm for LP by analyzing
the local behavior of its trajectories. In Section~\ref{SecMatOfDef}
we review the notions of hyperbolic polynomial and hyperbolicity
cones giving the first proof of G\aa rding's key result on cones'
convexity that does not rely on complex variables; in
Section~\ref{SecSWIntro} we introduce Shrink-Wrapping for LP; in
Section~\ref{SecSWLoc} we analyze the behavior of Shrink-Wrapping
trajectories in the proximity of a certain invariant set that
contains optimal LP solution and, as a consequence, describe a
simple idealized locally super-quadratically convergent discrete
bi-section scheme; in Section~\ref{SecSWvsIPM} we contrast
Shrink-Wrapping trajectories with the so-called central path for
some pathological LP instances.

\section{Basics}\label{SecMatOfDef}
\subsection{Hyperbolic polynomials and hyperbolicity cones}~\label{Hyper}

Mainly, we follow the exposition of elementary properties of hyperbolic
polynomials found in~\cite{Renegar:Poly}. However, unlike the original
proof of convexity of hyperbolicity cones in~\cite{Garding:Poly} or its
later version, e.g., in~\cite{Renegar:Poly}, our approach is rather
geometric and does not rely on complex numbers, thus, bringing it closer
to the spirit of continuous optimization.

Let $X$ be a finite-dimensional real vector space. Recall a
polynomial $p: X \rightarrow \Re$ is \emph{homogeneous of degree
$m$} if $p(tx)=t^mp(x)$ for all $t\in\Re$ and every $x\in X$.

\begin{defn}
Let $p: X \rightarrow \Re$ be a homogeneous polynomial of degree $m$
and $d\in X$ is such that $p(d)> 0$. $p$ is \emph{hyperbolic with
respect to $d$} if the univariate polynomial $t\mapsto
p(x+t d)$ has $m$ real roots for every $x\in X$.
\end{defn}

Examples:
\begin{itemize}
\item $X=\Re^m$, $d=\mathbf{1}\in\Re^m$ -- the vector of all ones.
The $m^{th}$ elementary symmetric polynomial $E_m(x)=\prod_{i=1}^m
x_i$ is a hyperbolic polynomial with respect to $\mathbf{1}$, since
$t\mapsto E_m(x+t\mathbf{1})=\prod_{i=1}^m (x_i+t)$ has roots
$-x_i,\:i=1,\ldots,m$,
\item $X=\rsym^m$ -- the space of real symmetric $m\times m$
matrices, $d=I$ -- the identity matrix. The determinant $\det(x)$ is
a hyperbolic polynomial with respect to $I$, since the eigenvalues
of $x\in X$ are minus the roots of $t\mapsto\det(x+t I)$ and are
real.
\end{itemize}
By analogy with the last example, given a hyperbolic polynomial $p$
and its hyperbolicity direction $d$, the roots of $\lambda\mapsto p(x-\lambda d)$
are called the \emph{eigenvalues} of $x$ in direction $d$, i.e., the eigenvalues
are precisely the roots of $t\mapsto p(x+td)$ with signs reversed. Ordering
eigenvalues in non-decreasing order, we denote them by
\[
\lambda_1(x)\leq\lambda_2(x)\leq\cdots\lambda_m(x).
\]
\begin{rem}
If a homogeneous polynomial is such that $t\mapsto
p(x+t d)$ has $m$ real roots for every $x\in X$ but $p(d)<0$, then
$-p(x)$ is hyperbolic with respect to $d$. Therefore, our definition may have been
augmented to require only $p(d)\neq 0$. In turn, $p(d)\neq 0$ is
an essential requirement to preserve the resulting cone's convexity and thus may not be
further relaxed; for an illustrative example see~\cite{Renegar:Poly}.
\end{rem}

Recall that a set is a cone if it is closed under multiplication by
nonnegative reals.
\begin{defn}
The \emph{hyperbolicity cone of $p$} with respect to $d$, written
$\mathcal{C}(d)$, is the set $\{x\in X:p(x+t d)\neq 0,\forall
t\geq 0\}$.
\end{defn}
We omit $p$ from the notation above as it will be clear which polynomial we refer to.
$\mathcal{C}(d)$ is a cone by homogeneity of $p$.

Examples:
\begin{itemize}
\item $X=\Re^m$, $d=\mathbf{1}$, $p(x)=E_m(x)$, then $\mathcal{C}(d)=\Re^m_{++}$ is strictly positive orthant,
\item $X=\rsym^m$, $d=I$, $p(x)=\det(x)$,
then $\mathcal{C}(d)$ 
is the cone of positive definite matrices.
\end{itemize}

\begin{prop}~\label{Hyper:ElemProp} Given a hyperbolic polynomial $p$
and its hyperbolicity direction $d$,
\begin{itemize}
\item[\textup{(A)}] for fixed real $\alpha\geq 0$ we have
$\lambda_i(\alpha x)=\alpha\lambda_i(x),\:i=1,\ldots,m$,
\item[\textup{(B)}] for fixed real $\beta$ we have $\lambda_i(x+\beta d)=\lambda_i(x)+\beta$,
\item[\textup{(C)}] if $e\in\mathcal{C}(d)$, then the linear segment $[e,d]\subset\mathcal{C}(d)$,
and, more generally, $[\gamma e,\delta d]\subset\mathcal{C}(d)$ for any $\gamma,\delta>0$.
\end{itemize}
\end{prop}
\begin{proof}
Part A follows immediately from the definition of eigenvalues.
Just as the case of symmetric matrices, part B is readily established
by a simple regrouping of variables $p((x+\beta d)-\lambda d)=p(x-(\lambda-\beta) d)$.
Part C follows from A and B: observe that for any $\xi\in(0,1)$ we have
\begin{eqnarray*}
\lambda_i(\xi e+(1-\xi)d)=\xi \lambda_i\left(e+\frac{1-\xi}{\xi}d\right)=
\xi \left(\lambda_i(e)+\frac{1-\xi}{\xi}\right)> 0,
\end{eqnarray*}
for all $i=1,\ldots,m$, since $\lambda_i(e)>0$, and so $\xi e+(1-\xi)d \in\mathcal{C}(d)$;
similarly, a more general statement follows.
\end{proof}
As a straightforward consequence, we can make two important observations.
\begin{prop}
$\mathcal{C}(d)=\{x\in X:\lambda_1(x)>0\}$.
\end{prop}
\begin{proof}
Follows from Proposition~\ref{Hyper:ElemProp} part B and $p(x)=p(d)\prod_{i=1}^m \lambda_i(x)$,
where the coefficient $p(d)>0$ in the identity is a consequence of considering
$\lim_{t\uparrow\infty}p(x+td)$ and homogeneity of $p$.
\end{proof}
\begin{prop}~\label{Hyper:ConeP}
$\mathcal{C}(d)$ is a (linearly) connected component of $\{x\in X: p(x)>0\}$ containing $d$.
\end{prop}
\begin{proof}
From the previous proposition it follows that $\mathcal{C}(d)\subset\{x\in X: p(x)>0\}$.
Clearly $\lambda_i(d)=1$ for all $i=1,\ldots,m$,
so $d\in\mathcal{C}(d)$. To establish connectivity of $\mathcal{C}(d)$,
for $e,f\in\mathcal{C}(d)$ observe $[e,d]\cup[d,f]\subset\mathcal{C}(d)$
by Proposition~\ref{Hyper:ElemProp} part C.
\end{proof}

Another important conclusion to be made from Proposition~\ref{Hyper:ElemProp}
is that the cone $\mathcal{C}(d)$ is defined \emph{locally} around its
hyperbolicity direction $d$. That is, in order to describe $\mathcal{C}(d)$ it
suffices only to know the behavior of the eigenvalues in the small ball around $d$, while
the properties A, B and C tell us explicitly how to compute the boundary of the closure of the
cone $\mathcal{C}(d)$ having this information. For ease of reference, we
distill the above mentioned computational procedure into the following statement.
\begin{prop}~\label{Hyper:LineRoot}
For fixed real $\omega$, the roots $t_i$ of $t\mapsto p(d+\omega(e-d)+td)$ satisfy
\begin{equation*}
t_i = \omega(1-\lambda_i(e))-1.
\end{equation*}
\end{prop}
\begin{proof}
Observe
$p(d+\omega(e-d)+td)=p(\omega e+(1-\omega+t)d)=\omega^m\: p\left(e+\frac{1-\omega+t}{\omega}d\right)$.
\end{proof}
That the cone is defined \emph{locally} is also a straightforward
consequence of analyticity of $p$. This \emph{localization} around $d$ is
a key principle that allows us to establish the convexity of $\mathcal{C}(d)$ as a
corollary to the following important property.
\begin{thm}~\label{Hyper:KeyProp}
If $e\in\mathcal{C}(d)$ then $p$ is hyperbolic with respect to $e$ and $\mathcal{C}(e)=\mathcal{C}(d)$.
\end{thm}

Example:
\begin{itemize}
\item $X=\rsym^4$, $d=I$, $p(x)=\det(x)$; note that not every element $e\in X, \:p(e)>0$,
gives rise to a hyperbolicity direction, e.g., consider eigenvalues of a linear matrix pencil
corresponding to $\det(A-\lambda B)=0$ where
\begin{eqnarray*}
A=\left(
\begin{array}{cccc}
0 & 0 & 0 & 1\\
0 & 0 & 1 & 0\\
0 & 1 & 0 & 0\\
1 & 0 & 0 & 0
\end{array}
\right),
&
B=\left(
\begin{array}{cccc}
0 & 1 & 1 & 0\\
1 & 0 & 0 & -1\\
1 & 0 & 0 & 0\\
0 & -1 & 0 & 0
\end{array}
\right).
\end{eqnarray*}
The eigenvalues $\lambda$ are $\pm i$ each with multiplicity two, although $\det(B)=1$.
(We would like to thank Prof. Peter Lancaster for providing this example.)
Interestingly, in case of $X=\rsym^2$ the condition $\det(B)>0$ suffices to ensure
that all the eigenvalues of a matrix pencil $A-\lambda B$ are real, as it amounts
to $B$ being either positive or negative definite -- the latter is a standard
sufficient condition for the linear matrix pencil spectra to be real, which may be easily
established by say pre and post-multiplying $A-\lambda B$ by inverse Cholesky factors
of $B$ or $-B$.
\end{itemize}
Given a homogeneous polynomial $p$, an intriguing question is to
characterize all $e\in X$ giving hyperbolicity directions to $p$, if
such exist.

Before we prove Theorem~\ref{Hyper:KeyProp}, we start with more elementary
but illustrative exercise of showing that one may perturb $d$ ever so slightly
maintaining hyperbolicity of $p$. This proof, with some minor modifications, essentially
carries over to the proof of our theorem in question.
\begin{prop}
Given a hyperbolic polynomial $p$ and its hyperbolicity direction $d$, there exists $\varepsilon>0$
such that for any $\Delta, \:\|\Delta\|\leq\varepsilon$, polynomial $p$ is hyperbolic with respect to
$d+\Delta$.
\end{prop}
\begin{proof}
By Proposition~\ref{Hyper:ElemProp} parts A and B, it suffices to
show that there exists an open neighborhood around
$\widetilde{d}=d+\Delta$ such that for any point $x$ in this
neighborhood $t\mapsto p(x+t\:\widetilde{d})$ has all real roots;
any point $z\in X$ may be shown to have associated roots of
$t\mapsto p(z+t\:\widetilde{d})$ real by translating $z$ along
vector $\widetilde{d}$ to a properly scaled version of this
neighborhood, where every point in the neighborhood including
$\widetilde{d}$ is carried into its multiple by some fixed positive
constant, see Figure~\ref{Hyper:Fig1}.
\begin{figure}[h!]
\begin{center}
\epsfig{file=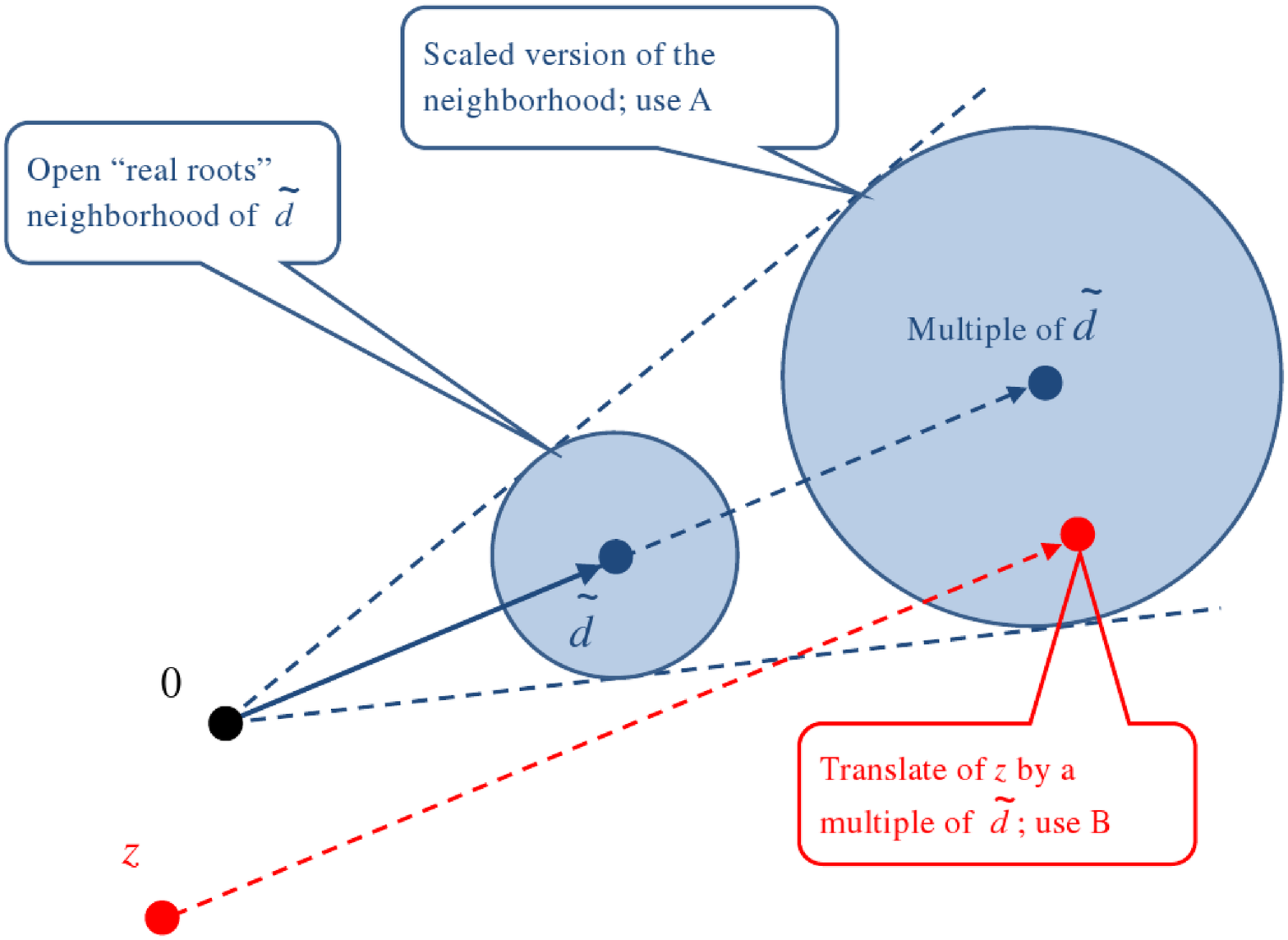,height=10.0cm} \caption{Localization of real roots around $\widetilde{d}$}\label{Hyper:Fig1}
\end{center}
\end{figure}

Consider $p$ of degree $m$, hyperbolic with respect to $d$.
Clearly $\lambda_i(d)=1,\: i=1,\ldots,m$.
Let $B_{\varepsilon}(d)$ be an open ball of radius $\varepsilon>0$ around $d$ such that
$\forall y\in B_{\varepsilon}(d)$ we have $|\lambda_i(y)-1|<\frac{1}{2}, \: i=1,\ldots,m$.
Such a ball exists by continuity of $\lambda_i$.

Fix $\Delta\in X,\:\|\Delta\|\leq\varepsilon$, and consider a mapping
\begin{equation*}
\tau\mapsto p(x+\tau\Delta+td)
\end{equation*}
for each fixed $x$ and real $\tau$ producing a polynomial in $t$. Observe that for any
$x\in B_{\frac{\varepsilon}{2}}(d+\Delta)=\left\{x\in X: x=(d+\Delta)+\delta, \|\delta\|<\frac{\varepsilon}{2}\right\}$
and $\tau\in[-3/2,-1/2]$ we have
\[
\|x+\tau\Delta-d\|=\|(d+\Delta)+\delta+\tau\Delta-d\|\leq
\|\delta\|+|1+\tau| \|\Delta\|< \frac{\varepsilon}{2}+\frac{\varepsilon}{2}=\varepsilon.
\]
So, for any fixed $x\in B_{\frac{\varepsilon}{2}}(d+\Delta)$ and $\tau\in[-3/2,-1/2]$,
a polynomial $t\mapsto p((x+\tau\Delta)+td)$ has all real roots in the interval $[-3/2,-1/2]$ by the choice of $\varepsilon$.

Now, in order to show that $t\mapsto p(x+t(d+\Delta))$ has all real
roots $\forall x\in B_{\frac{\varepsilon}{2}}(d+\Delta)$, start increasing
$\tau$ in $t\mapsto p(x+\tau\Delta+td)$ from -3/2 to -1/2, see
Figure~\ref{Hyper:Fig2}; whenever $t=\tau$ intersects $t=\lambda_i(x+\tau\Delta)$, we capture one of the
desired real roots, increasing $\tau$ until we extract $m$ roots.
\begin{figure}[h!]
\begin{center}
\epsfig{file=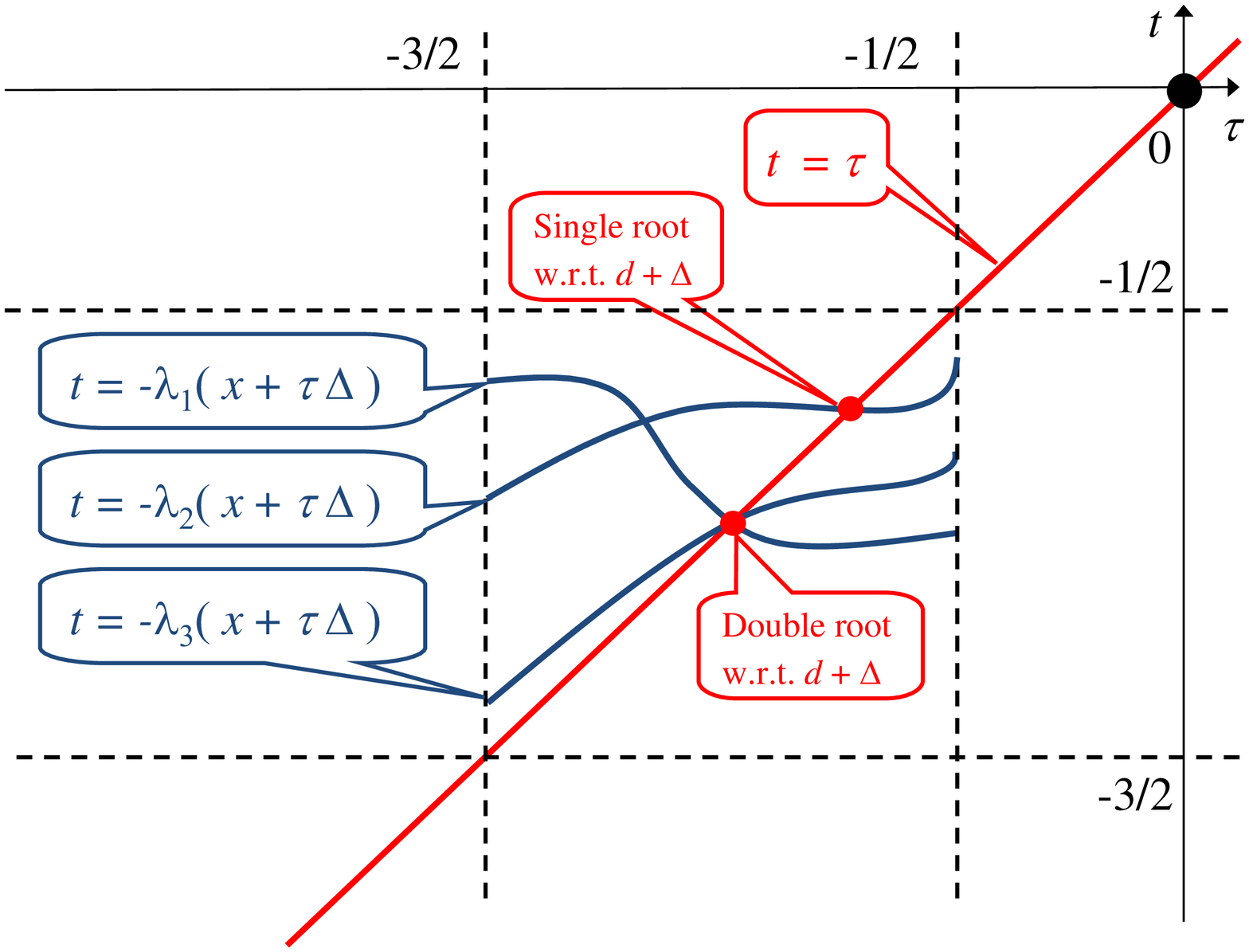,height=10.0cm} \caption{Identifying $m$
real roots of $t\mapsto p(x+t\:\widetilde{d})$}\label{Hyper:Fig2}
\end{center}
\end{figure}
\end{proof}

Note that the proof of the proposition does not rely on the initial
neighborhood of $d$ being a ball. Thus, we extend the proof to
establish hyperbolicity of $p$ with respect to $e\in\mathcal{C}(d)$.

\begin{proof}[Proof of Theorem~\ref{Hyper:KeyProp}]
To establish hyperbolicity of $p$ with respect to $e\in\mathcal{C}(d)$,
just as before, it suffices to show that there is an open neighborhood of $e$ such that
$t\mapsto p(x+te)$ has $m$ real roots for all $x$ in the neighborhood.

By homogeneity of $p$, without loss of generality we may assume
$0<2\gamma<\lambda_i(e)<1-2\gamma<1,\: i=1,\ldots,m$. By Proposition~\ref{Hyper:LineRoot} for $e'=d-(e-d)$
we have $0<2\gamma<\lambda_i(e)<1-2\gamma<1,\: i=1,\ldots,m$, and so the linear segment $[e,e']$ belongs
to $\mathcal{C}(d)$ by Proposition~\ref{Hyper:ElemProp} part C. Moreover, by
Proposition~\ref{Hyper:LineRoot} and continuity of $\lambda_i$
there exists $\varepsilon>0$ such that the open ``tubular'' neighborhood $\mathcal{N}$ around
$[e,e']$ consisting of convex combination of two open balls of radius $\varepsilon$ around $e$ and $e'$,
$\mathcal{N}=\textup{conv}(B_{\varepsilon}(e),B_{\varepsilon}(e'))$, satisfies
$\forall y\in \mathcal{N}$ we have $\lambda_i(y)\in(\gamma,2-\gamma), \: i=1,\ldots,m$.

Consider $\Delta=e-d$ and $\tau\mapsto p(x+\tau\Delta+td)$. Note
$x+\tau\Delta\in\mathcal{N}$ for all $x\in B_{\varepsilon}(e)$ and $\tau\in[-2,0]$. So, for any
fixed $x\in B_{\varepsilon}(e)$ and $\tau\in[-2,0]$ the polynomial $t\mapsto p(x+\tau\Delta+td)= p(x+\tau(e-d)+td)$
has $m$ real roots in the interval $(-2+\gamma,-\gamma)$.
Therefore, by increasing $\tau$ from -2 to 0 we may identify $m$ real roots of
$t\mapsto p(x+t(d+\Delta)=p(x+te)$ as intersections of $t=\tau$ and $t=-\lambda_i(x+\tau\Delta)$.

Finally, that $\mathcal{C}(e)=\mathcal{C}(d)$ easily follows from Proposition~\ref{Hyper:ConeP}.
\end{proof}

We are in position to prove the convexity of $C(d)$; as observed in~\cite{Renegar:Poly} this is a
consequence of Theorem~\ref{Hyper:KeyProp}, we restate the proof for completeness.
\begin{thm}
$\mathcal{C}(d)$ is an open convex cone.
\end{thm}
\begin{proof}
$\mathcal{C}(d)$ is an open set by continuity of $p$ and Proposition~\ref{Hyper:ConeP}.
Consider $x,y\in\mathcal{C}(d)$. Note $\mathcal{C}(y)=\mathcal{C}(d)$ and so by
Proposition~\ref{Hyper:ElemProp} part C we have $[x,y]\in\mathcal{C}(y)=\mathcal{C}(d)$.
\end{proof}

In~\cite{Garding:Poly} the convexity of $\mathcal{C}(d)$ was established as a corollary to the following result.
\begin{fct}\label{Gaarding:concave}
$\lambda_1(x)$ is a concave function of $x$.
\end{fct}
\noindent Indeed, if $x,y\in\mathcal{C}(d)$ and $\lambda_1(x)$ is concave, then for any $\xi\in(0,1)$ we have
$\lambda_1(\xi x+(1-\xi)y)\geq\xi\lambda_1(x)+(1-\xi)\lambda_1(y)>0$. Later in~\cite{Renegar:IPM}
it was shown that conversely, the concavity of $\lambda_1(x)$ follows from
convexity of $\mathcal{C}(d)$, using much simpler proofs yet still relying on complex numbers.
%
%
If we introduce sums of the smallest $k$ eigenvalues
\[
s_k=\sum_{i=1}^k \lambda_i,
\]
a more general statement regarding the eigenvalues may be established~\cite{Lewis:Poly}.
\begin{fct}\label{Lewis:concave}
$s_k(x)$ is a concave function for any $k=1,\ldots,m$.
\end{fct}

Next, we turn our attention to the so-called hyperbolic derivatives.

\subsection{Hyperbolic derivatives and cone characterization}

Given a hyperbolic polynomial $p$ of degree
$m$ and its hyperbolicity direction $d$, the directional derivative
of $p(x)$ along $d$ is called
\emph{the hyperbolic derivative polynomial of $p$ with respect to $d$}, denoted
\[
p'(x)=p'_t(x+td)|_{t=0}.
\]
We refer to $p'$ simply as the \emph{derivative polynomial} of $p$,
omitting $d$ for brevity of notation. By the root interlacing
property for the polynomials with all real roots --by continuity,
for fixed $x$ between any two roots of $t\mapsto p(x+td)$ there is a root of
$t\mapsto p'_t(x+td)$-- it follows that $p'(x)$ is also hyperbolic with respect to $d$.

Similarly, for a fixed hyperbolicity direction $d$, we can define
higher derivatives $p'',p''',\ldots,p^{(m)}$ as $p^{(k)}(x)=p^{(k)}_t(x+td)|_{t=0}$.
Since $p$ is of degree $m$, $p^{(m-1)}$ is linear and $p^{(m)}(x)$ is constant.

Examples:
\begin{itemize}
\item $X=\Re^m$, $d=\mathbf{1}$, $p(x)=E_m(x)$, then
\[E_m^{(k)}(x)=k! E_{m-k}(x)\] where $E_j(x)$ is the $j^{th}$ elementary
symmetric polynomial
\[
\begin{array}{cccc}
E_1(x) = \sum_{1 \leq i \leq m} x_i, & E_2(x) = \sum_{1 \leq i < j
\leq m} x_i x_j, & \ldots, & E_m(x)=\prod_{1\leq i\leq m} x_i,
\end{array}
\]
\item $X=\Re^m$, $d\in\Re^n_{++}$, $p(x)=E_m(x)$; by a similar inductive argument
as above one can show that
\[
E_m^{(k)}(x)=k! E_m(d) E_{m-k}\left(\left[
\frac{x_1}{d_1},\frac{x_2}{d_2},\ldots,\frac{x_m}{d_m}\right]\right).
\]
\end{itemize}

\begin{rem}
The elementary symmetric polynomials in the example above play
an important role in representing the derivative polynomials via the
eigenvalues at $x\in X$. Namely, since
$p(x+td) = p(d) \prod_{1\leq i\leq m} (t+\lambda_i(x))$ we have
\begin{eqnarray*}
p'(x) = \left.\frac{\partial}{\partial t} \left(p(d) \prod_{1\leq i\leq m}
(t+\lambda_i(x))\right)\right|_{t=0} = p(d) \sum_{1\leq i\leq m}
\prod_{j\neq i} \lambda_j(x) = p(d) E_{m-1}(\lambda(x))
\end{eqnarray*}
where $\lambda(x)$ is the vector of $m$ eigenvalues of $x$, and more generally
\[
p^{(k)}(x) = k! p(d) E_{m-k}(\lambda(x)).
\]
\end{rem}

For the $k^{th}$ hyperbolic derivative of $p$, we use $\mathcal{C}^{(k)}(d)$ to denote
the associated hyperbolicity cone; note $\mathcal{C}^{(m-1)}(d)$ is
an open half-space and $\mathcal{C}^{(m-1)}(d)=X$. Although, $\mathcal{C}(e)=\mathcal{C}(d)$
for any $e\in \mathcal{C}(d)$, the hyperbolicity cones corresponding to derivative polynomials
$p',p'',\ldots$ with respect to $e\neq d$ might not coincide with one another,
as in the last example where, for instance, $k=m-1$.

It turns out that these derivative polynomials come in handy in
characterization of the hyperbolicity cone $\mathcal{C}(d)$ itself as observed in~\cite{Renegar:Poly}.
We note that if all $\lambda_i(x)>0$ then clearly $p^{(k)}(x)>0$ for all $k=1,\ldots,m$. Conversely,
by Taylor series of $p(x+td)$,
\begin{eqnarray*}
p(x+td) = p(x)+p'(x) t+p''(x) \frac{t^2}{2!}+\cdots+p^{(m-1)}(x)\frac{t^{m-1}}{(m-1)!}+p^{(m)}(x)\frac{t^m}{m!},
\end{eqnarray*}
observe that if all $p^{(k)}(x)>0$, then $p(x+td)>0$ for all $t>0$, and thus $x\in\mathcal{C}(d)$.
\begin{fct}~\label{Hyper:ConeChar}
The hyperbolicity cone satisfies
\[
\mathcal{C}(d)=\{x\in X:\: p(x) > 0, p'(x) > 0, p''(x) > 0, \ldots, p^{(m-1)}(x) > 0 \}.
\]
\end{fct}
As an important consequence of this fact we have the following cone inclusion.
\begin{cor} ~\label{Hyper:ConeInclusion}
\[
\mathcal{C}(d) \subseteq \mathcal{C}'(d) \subseteq \mathcal{C}''(d) \subseteq \cdots \subseteq
\mathcal{C}^{(m-1)}(d).
\]
\end{cor}

Throughout the rest of the manuscript we will be concerned with the closure of a hyperbolicity cone,
$\closure\mathcal{C}(d)$. Due to continuity of $p(x)$ all the results in this and previous
subsections naturally extend to $\closure\mathcal{C}(d)$ by replacing strict inequalities
with corresponding inequalities when necessary. To this end, we note that, for example, the closed cone
$\closure\mathcal{C}(d)=\{x\in X:\lambda_1(x)\geq0\}=\{x\in X:\: p^{(k)}(x)\geq0, k=1,\ldots,m-1\}\subseteq \closure\mathcal{C}'(d)\subseteq \cdots\subseteq \closure\mathcal{C}^{(m-1)}(d)$ is convex, etc.;
likewise, one may easily characterize the boundary of $\closure\mathcal{C}(d)$ as follows, see~\cite{Renegar:Poly}.
\begin{cor}~\label{Hyper:ConeBdry}
$\partial\left( \closure\mathcal{C}(d)\right)=\{x\in X:\: p(x)=0, p'(x)\geq0,\ldots,p^{(m-1)}(x)\geq0\}.$
\end{cor}

\begin{prop}~\label{Hyper:MultRoot}
If $x\in \closure\mathcal{C}^{(r)}(d) \bigcap \closure\mathcal{C}^{(r+1)}(d)$ for some $r>0$, then $x\in\closure\mathcal{C}(d)$.
\end{prop}
\begin{proof}
By the inclusion property for derivative cones, $x$ must belong to the boundary of both cones,
$x\in \partial\left(\closure\mathcal{C}^{(r)}(d)\right) \bigcap \partial\left(\closure\mathcal{C}^{(r+1)}(d)\right)$.
Consequently, by the root interlacing property for polynomials with all real roots,
it follows that 0 is a root of multiplicity $\mu\geq2$ corresponding to $t\mapsto p^{(r)}(x+td)$:
by contradiction, if 0 has multiplicity 1, then the derivative polynomial $t\mapsto p^{(r+1)}(x+td)$
cannot have 0 as its root. Analogously, if $r>1$ then $t\mapsto p^{(r-1)(x+td)}$ must have
0 as its root of multiplicity $\mu+1$, etc. So, indeed, $x\in\closure\mathcal{C}(d)$, and, in particular,
$x$ lies on the boundary of the cone.
\end{proof}

Example:
\begin{itemize}
\item $X=\Re^n$, $d\in\Re^n_{++}$, $p(x)=E_n(x)$; let $\mathcal{K}_{r,d}$
denotes the closure of the hyperbolicity cone associated with $r^{th}$ derivative
polynomial of $p$ with respect to $d$, $\mathcal{K}_{r,d}=\closure\mathcal{C}^{(r)}(d)$.
The following cone inclusion
\[
\Re^n_+=\mathcal{K}_{0,d} \subseteq \mathcal{K}_{1,d} \subseteq \cdots
\subseteq \mathcal{K}_{n-1,d} \subseteq \mathcal{K}_{n,d}=\Re^n
\]
gives a natural sequence of relaxations of the nonnegative orthant $\Re^n_+$,
a pivotal observation for building Shrink-Wrapping framework for linear programming.
Note that $K_{r,d}$ coincides with the closure of hyperbolicity cone associated
with $E_{n-r}(x./d)$, where $x./d$ is a componentwise ratio of vectors $x$ and $d$;
observe $K_{r,d}=\{x\in\Re^n: x=d.\cdot z, z\in\mathcal{K}_{r,\mathbf{1}}\}$, where $d.\cdot z$
is a componentwise product of two vectors; in particular,
$K_{n-1,d}$ is a half-space passing through
the origin with normal vector $\mathbf{1}./d$.
\end{itemize}

\begin{rem}
Interestingly, for $x\in\partial\left(\closure\mathcal{C}(d)\right)$
we have $\nabla\lambda_1(x)$ parallel to $\nabla p(x)$. Let
$X=\Re^n$; considering
$x\in\partial\left(\closure\mathcal{C}(d)\right)$ so that
$0=\lambda_1(x)<\lambda_2(x)$, recall
$p(x)=p(d)\prod_{j=1,m}\lambda_j(x)$ and so
\[
\frac{\partial}{\partial x_i}\:p(d)\prod_{j=1,m}\lambda_j(x)=p(d)\sum_{j=1,m}\frac{\partial}{\partial x_i} \lambda_j(x)\prod_{k\neq j}\lambda_k(x)=p(d)\prod_{k\neq 1} \lambda_k(x) \cdot \frac{\partial}{\partial x_i} \lambda_1(x)
\]
giving
\begin{eqnarray*}
\nabla p(x) = p(d)\prod_{k\neq1}\lambda_k(x)\cdot \nabla\lambda_1(x).
\end{eqnarray*}
Similarly, if $t\mapsto p(x+td)$ has 0 as its root of multiplicity $\mu>1$, one may consider the boundary of
the corresponding derivative cone $\closure\mathcal{C}^{(\mu-1)}(d)$ instead.
\end{rem}
Although, there is a simple algebraic characterization of the
hyperbolicity cones, their dual cones are poorly understood, with
some exceptions, e.g.,~\cite{Chua:Poly, Zinchenko:ESP}.

\subsection{Hyperbolic programs and relaxations}

The significance of hyperbolicity cones in convex optimization
becomes evident once we introduce the three most prominent instances
of the so-called conic programming problems. Letting $X$ be equipped
with an inner product $\langle \cdot, \cdot \rangle$, a {\em conic
programming problem} is an optimization problem of the form
\[
\inf_x \{\langle c,x \rangle: Ax=b, x\in K\}
\]
where $K\subset X$ is a closed convex cone, $c\in X$, $b\in\Re^m$
and $A:X\rightarrow \Re^m$ -- a linear operator. It is well known
that any convex optimization problem can be recast as conic
programming problem.

The three most prominent instances of conic programming are:
\begin{itemize}
\item LP, $X=\Re^n, \langle x,y \rangle=x^T
y, K=\Re^n_+$,
\item Second-Order Conic
Programming (SOCP), $X=\Re^n, \langle x,y \rangle=x^T y, K=K_1
\times K_2 \times \cdots \times K_{\ell}$ with second-order cones
$K_i=\{(x,t)\in\Re^{n_i-1}\times\Re:\|x\|\leq t\}$,
$\sum_{i=1}^{\ell} n_i=n$, and
\item positive Semi-Definite Programming (SDP), $X=\rsym^m$,
$\langle x, y \rangle=\trace(xy)$ and $K$ -- the cone of positive semi-definite
matrices.
\end{itemize}
In applications, these three types of problems provide an extremely
powerful modeling framework, ranging from production planning,
relaxations to hard combinatorial problems, mathematical finance and
Markov chains, to control theory and polynomial programming
\cite{Chvatal:LP},\cite{GW:Maxcut},\cite{BP:Finance},\cite{Boyd:MC},\cite{BenTalN:ConvOpt},\cite{Lasserre:Poly}.
Also, they naturally arise as robust
counterparts~\cite{BenTalN:ConvOpt} to one another in the presence
of uncertainty in the initial data, e.g.,~\cite{SMY:IMRT}.

A Hyperbolic Programming (HP) problem is a conic programming problem
where $K$ is a closure of hyperbolicity cone. Note LP, SOCP and SDP
are instances of HP.

\begin{rem}
When implementing an interior-point method for SDP it is frequently
required to determine how far one may advance along a given vector $h\in\rsym^m$
from some point $e\in\rsym^m$ in the cone of positive definite matrices, before hitting the boundary
of the closure of this cone. Typically, the procedure is considered to be computationally expensive
due to its implementation as ``trial and error'' testing on whether a given vector $e+\omega h$ is still in the cone, $\omega\in\Re$. Theorem~\ref{Hyper:KeyProp} combined with Proposition~\ref{Hyper:LineRoot}
gives an elegant basis for an alternative relatively inexpensive procedure.
Note that with respect to $\det(\cdot)$, $\mathcal{C}(e)$
coincides with the cone of positive definite matrices. Using Cholesky factors of $e=LL^T$, one
may compute the largest eigenvalue of $L^{-1}(e-h)L^{-T}$ or its approximation, say, using Lanczos-type
algorithm provided $e-h$ is also positive-definite, and subsequently use this value to determine
the maximum allowed step-length along $h$ using Proposition~\ref{Hyper:LineRoot}.
\end{rem}

In what follows, within Shrink-Wrapping framework, together with a linear programming instance
\begin{eqnarray*}\label{Problem:LP}
\min_x\{c^T x:\: Ax=b, x\in\Re^n_+\} & (LP)
\end{eqnarray*}
where $c\in\Re^n, b\in\Re^m, A\in\Re^{m\times n}, m<n$,
we consider its \emph{$r^{th}$ hyperbolic relaxation} with respect to some fixed $d\in\Re^n_{++}$, $0<r<n-1$,
\begin{eqnarray*}\label{Problem:HP}
\min_x\{c^T x:\: Ax=b, x\in\mathcal{K}_{r,d}\} & (HP_{r,d})
\end{eqnarray*}
recalling that $\mathcal{K}_{r,d}$ is the closure of hyperbolicity cone corresponding to
$r^{th}$ hyperbolic derivative $r! E_n(d) E_{n-r}(x./d)$ of $E_n(x)$
with respect to $d$. Let $x^*$ and $x(d)$ denote optimal solutions for $LP$ and $HP_{r,d}$ respectively;
for convenience, we are assuming $x^*$ is a unique minimizer for $LP$.

Example:
\begin{itemize}
\item consider linear programming problem $\min_x\{c^T x:\: \mathbf{1}^Tx=3, x\in\Re^3_+\}$
together with its first-order hyperbolic relaxation $\min_x\{c^T
x:\: \mathbf{1}^Tx=3, x\in\mathcal{K}_{1,d}\}$ where $d=\mathbf{1}$,
see Figure~\ref{Hyper:Fig3}; note that the feasible region of
$HP_{1,\mathbf{1}}$ is inscribed by a circle in
$\{x\in\Re^3:\:\mathbf{1}^Tx=3\}$ centered around $d=\mathbf{1}$.
\begin{figure}[h!]
\begin{center}
\epsfig{file=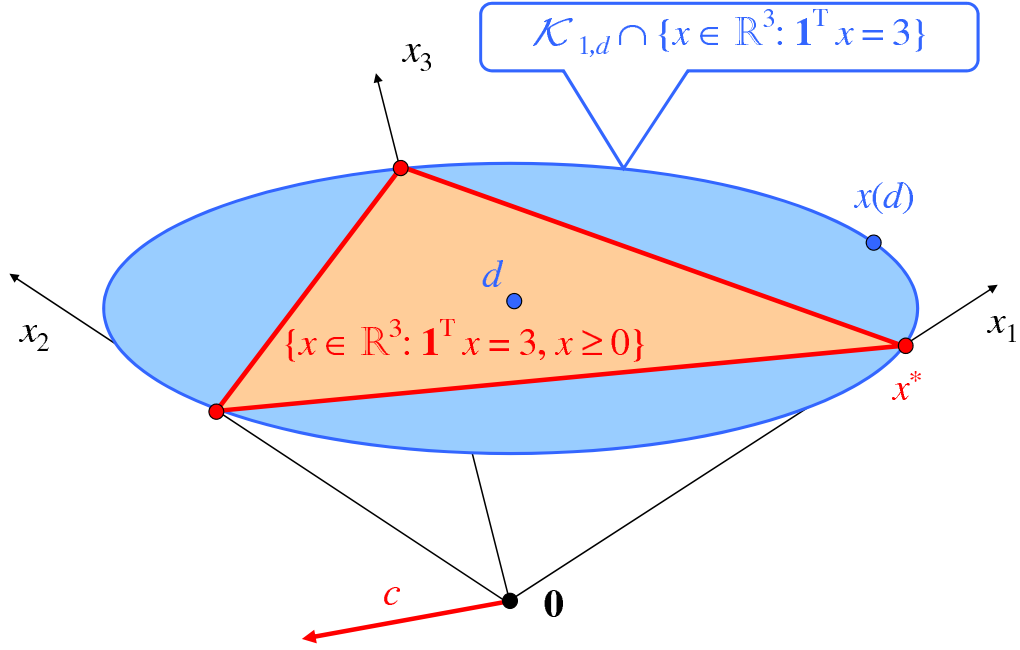,height=10.0cm} \caption{$LP$ and its hyperbolic relaxation}\label{Hyper:Fig3}
\end{center}
\end{figure}
\end{itemize}
Although, the above example is fairly simple, it illustrates
a few key geometric concepts of Shrink-Wrapping throughout
the manuscript .

Hyperbolic relaxations $HP_{r,d}$ will be used to define a
family of continuous trajectories terminating at the $LP$ optimum.
In a similar fashion, one may define hyperbolic relaxations
for any other HP instance besides LP, including SOCP and SDP.

\section{Shrink-Wrapping approach for LP}\label{SecSWIntro}
\subsection{Main ingredients}

\begin{prop}
If bounded, $HP_{r,d}$ has a unique solution $x(d)$ unless $x(d)$ solves $LP$.
\end{prop}
\begin{proof}
The boundary of $\mathcal{K}_{r,d}$ at $x$ has strict curvature except along $x$ itself, unless $x\in\Re^n_+$,
see~\cite{Renegar:Poly}, Theorem 14; note that the only flat faces of $\mathcal{K}_{r,d}$ are precisely
$n-r-1$ and lower dimensional faces of the nonnegative orthant, since in the latter case $E_{n-r}(x./d)=0$.
\end{proof}

We are interested in recovering $LP$ solution using hyperbolic
relaxations. Although, our present investigation is mostly of
theoretical nature, we would like to comment on practicality of the
underlying assumptions to indicate potential usability of this new
setting. To this extent, we assume that
\begin{itemize}
\item[(A)] $LP$ is bounded,
\item[(B)] we know an initial strictly $LP$-feasible point $d\in\{x\in\Re^n: Ax=b, x\in\Re^n_{++}\}$,
\item[(C)] the corresponding hyperbolic relaxation $HP_{r,d}$ is bounded as well,
\item[(D)] $x(d)$ is not $LP$-optimal,
\item[(E)] we can easily solve $HP_{r,d}$ to find $x(d)$.
\end{itemize}
Hypotheses (A) and (B) are fairly standard assumptions for linear
programming, in particular in the context of interior-point methods.
In fact, instead of (A) and (B) one frequently relies on even more
restrictive hypothesis (A1) that the feasible region of $LP$ is
bounded, and (B1) that there is an \emph{affine feasible}
$x\in\Re^n_{++}:\: Ax=b$ and the $LP$ is strictly feasible, i.e.,
remains feasible under all infinitesimal perturbations of $b$ -- the
latter implied, for example, by having $\rank (A)=m$ and a feasible
point $x\in\Re^n_{++}$. Note that even stronger (A1) and (B1) are
quite reasonable from a practical point of view: if $LP$ is used to
model a certain physical phenomenon, it is natural to assume the
compactness of its feasible region; in addition, a well thought
through model is typically feasible and avoids unnecessary state
variables and constraints, leading to strict feasibility. Also, with
regards to recovering the optimum, (A) and (A1) are not that much
different from one another: (A1) clearly implies (A), conversely, if
$\underline{val}$ is an a priori bound on the optimal value of $LP$,
then we might as well augment the feasible region of $LP$ by adding
a constraint $\mathbf{1}^Tx\leq n\frac{|\underline{val}|}{\|c\|}$,
thus making it compact. Shortly we will indicate that for all
reasonable LP instances, accommodating (C) should not pose
significant practical difficulties either; here, by a reasonable LP
instance we understand the feasible problem satisfying (A1). We use
(D) since otherwise we solved $LP$ already; henceforth, we refer to
$x(d)$ as \emph{the solution} of $HP_{r,d}$. Since our focus is on
analyzing continuous Shrink-Wrapping trajectories, we employ (E); in
more practical terms, one may think of setting up a Newton's method
based path-following scheme, e.g., similar to the so-called
short-step interior-point method, to recover a sufficiently good
approximation to $x(d)$.

Example:
\begin{itemize}
\item consider $\min_x\{(1, 1, 0)^T x:\: \mathbf{1}^Tx=3, x\in\Re^3_+\}$ and its relaxations
$\min_x\{(1, 1, 0)^T x:\: \mathbf{1}^Tx=3, x\in\mathcal{K}_{1,d}\}$ for three choices of $d\in\Re^3_{++}$:
(a) $d=\mathbf{1}$, (b) $d=(1, .1, 1.9)$, (c) $d=(.1, .1, 2.8)$. All $d$ are chosen affine
feasible, $\mathbf{1}^T d=3$; $x^*=(0,0,3)$ is $LP$ optimum.

Observe $2E_2(x./d)=(x./d)^T (\mathbf{1}\mathbf{1}^T-I) (x./d)$. The
boundary of $HP_{1,d}$ feasible region corresponds to $E_2(x./d)=0$
where $\mathbf{1}^Tx=3$. So, for affine feasible $x$ in the basis of
$x_1, x_2$ the boundary satisfies
\[
(x./d)^T (\mathbf{1}\mathbf{1}^T-I) (x./d)=\frac{1}{2}\xi^T Q
\xi+r^T\xi+s=0
\]
where $\xi=(x_1, x_2)$, and denoting $\Diag(z)$ the diagonal matrix with $\Diag(z)_{i,i}=z_i$,
\begin{eqnarray*}
Q=2
\left(
\begin{array}{rr}
1 & 0\\
0 & 1\\
-1 & -1
\end{array}
\right)^T
\Diag(\mathbf{1}./d)(\mathbf{1}\mathbf{1}^T-I)\Diag(\mathbf{1}./d)
\left(
\begin{array}{rr}
1 & 0\\
0 & 1\\
-1 & -1
\end{array}
\right),\\
r=2
\left(
\begin{array}{rr}
1 & 0\\
0 & 1\\
-1 & -1
\end{array}
\right)^T
\Diag(\mathbf{1}./d)(\mathbf{1}\mathbf{1}^T-I)\Diag(\mathbf{1}./d)
\left(
\begin{array}{c}
0\\
0\\
3
\end{array}
\right)
,\\
s=
\left(
\begin{array}{c}
0\\
0\\
3
\end{array}
\right)^T
\Diag(\mathbf{1}./d)(\mathbf{1}\mathbf{1}^T-I)\Diag(\mathbf{1}./d)
\left(
\begin{array}{c}
0\\
0\\
3
\end{array}
\right)
\end{eqnarray*}
In turn, for affine feasible $d>0$ writing $d_3=3-d_1-d_2$ we have
\[
Q=
\frac{3}{d_1d_2(3-d_1-d_2)}
\left(
\begin{array}{cc}
-2d_2 & 3-2(d_1+d_2)\\
3-2(d_1+d_2) & -2d_1
\end{array}
\right)=\frac{3\widetilde{Q}}{d_1d_2(3-d_1-d_2)}.
\]
Let us analyze the sign pattern for the eigenvalues of $\widetilde{Q}$; observe $\widetilde{Q}$ has at least one
negative eigenvalue as its diagonal is negative as well, also
\[
\det(\widetilde{Q})=-9+12(d_1+d_2)-4(d_1^2+d_2^2)-4d_1d_2.
\]
Note that for (a) where $d_1=d_2=1$, the matrix $Q$ is negative
definite, and thus the boundary of $HP_{1,d}$ feasible region indeed
corresponds to an ellipse. In both cases (b) and (c), where $d_1=1,
d_2=.1$ or $d_1=d_2=.1$, we have $\det(Q)<0$ and so the boundary
corresponds to a branch of hyperbola. In fact, for any sufficiently
small $d_1, d_2$ the boundary assumes hyperbolic shape, in
particular, when $d$ approaches $LP$ optimum, $d\rightarrow x^{*}$;
see Figure~\ref{Hyper:Fig4}. Note that in both cases (a) and (c)
$x(d)=x^*$ with the corresponding hyperbolicity directions belonging
to an open line segment $\mathcal{L}$ which extends to $x^*$.
\begin{figure}[h!]
\begin{center}
\epsfig{file=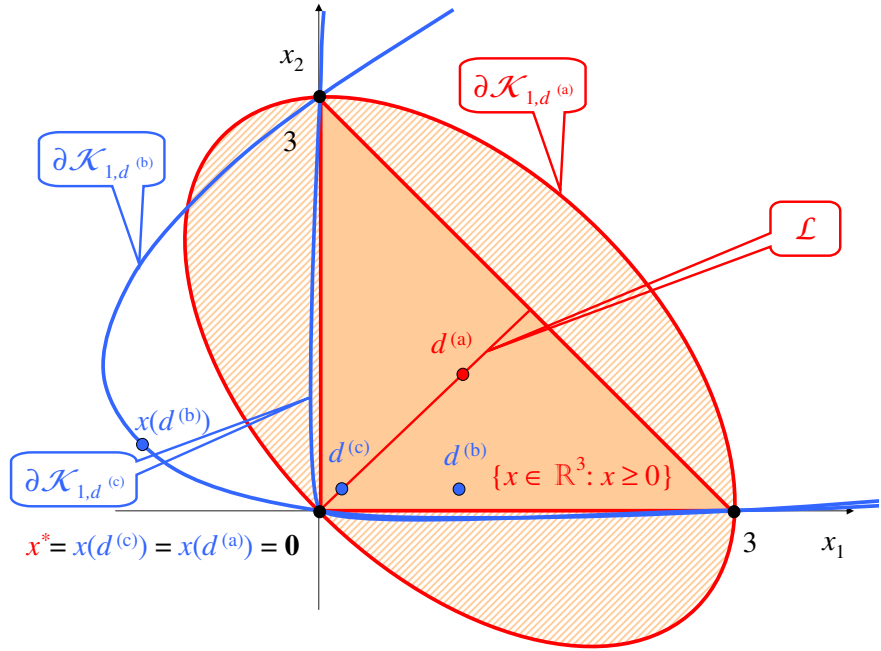,height=10.0cm} \caption{$HP_{1,d}$ in the basis of $x_1, x_2$ with varying $d$}\label{Hyper:Fig4}
\end{center}
\end{figure}
\end{itemize}
From this we make an important observation: although $LP$ has a bounded feasible region,
the feasible region corresponding to $HP_{r,d}$ may become unbounded.

Next, we are going to discuss (C). For a cone $K\subseteq\Re^n$, the
\emph{dual cone} is defined as $K^*=\{y\in\Re^n: x^T y \geq 0,
\forall x\in K\}$. More generally, the dual cone may be defined with
respect to an arbitrary inner product on $\Re^n$. A closed, convex
cone is \emph{regular} if both it has non-empty interior and its
lineality space is $\{\mathbf{0}\}$.
\begin{prop}
$\{x\in\Re^n:\mathbf{1}^Tx=n,\:x\in\mathcal{K}_{r,d}\}$ is bounded if and only if $d$ is in the interior of $\mathcal{K}_{r,\mathbf{1}}^*$.
\end{prop}
\begin{proof}
Recall that we consider $0<r<n-1$ and that $K_{r,d}=\{x\in\Re^m: x=d.\cdot z, z\in\mathcal{K}_{r,\mathbf{1}}\}$.
Note $\frac{n}{\mathbf{1}^Td}d\in \{x\in\Re^n:\mathbf{1}^Tx=n,\:x\in\mathcal{K}_{r,d}\}$; $\mathbf{1}^Td>0$ since $d\in\Re^n_{++}$.

Suppose $d\notin\mathcal{K}_{r,\mathbf{1}}^*$. There are two alternatives:
\begin{itemize}
\item[(i)] there exists $z\in\mathcal{K}_{r,\mathbf{1}}$ such that $0>z^T d=\mathbf{1}^T (d.\cdot z)$;
note we can take $\alpha\in(0,1)$ such that $u=\alpha (d.\cdot z) +(1-\alpha)d\in\mathcal{K}_{r,d}$ satisfies $\mathbf{1}^Tu=0$, $u\neq\mathbf{0}$ since $\mathcal{K}_{r,d}$ is regular; and, therefore, $\frac{n}{\mathbf{1}^Td}d+\tau u\in\{x\in\Re^n:\mathbf{1}^Tx=n,\:x\in\mathcal{K}_{r,d}\}$ for any $\tau\geq0$,
\item[(ii)] there exists $z\in\mathcal{K}_{r,\mathbf{1}},\:z\neq\mathbf{0}$ such that $z^T d=0$, simply
take $u=d.\cdot z$.
\end{itemize}
Conversely, if $d\in\mathcal{K}_{r,\mathbf{1}}^*$ the direction of unboundedness $u$ does not exist.
\end{proof}
\begin{cor}
Let $LP$ be such that $\mathbf{1}$
belongs to the range space of $A^T$.
If $d$ is in the interior of $\mathcal{K}_{r,\mathbf{1}}^*$, then
$HP_{r,d}$ has bounded feasible region.
\end{cor}
\noindent The condition is only sufficient, not necessary.

By strong LP duality the compactness of $LP$ feasible region is
equivalent to the existence of $e\in\Re^n_{++}$ such that $e$
belongs to the range space of $A^T$, i.e., $e=A^Ty$ for some
$y\in\Re^m$. So, assuming $LP$ has bounded feasible region with $e$
as above, if there exists $\bar{x}\neq\mathbf{0}$ feasible for $LP$,
we may add the constraint $e^T (x-\bar{x})=0$ to the $LP$ without
changing its feasible region. Therefore, by further scaling the $LP$
variables $x\mapsto e.\cdot x$ we may assume $e=\mathbf{1}$. The
last observation potentially allows us to identify some candidates
for the initial value of $d$ according to (C): $Ad=b$ with $d./e$ in
the interior of $\mathcal{K}_{r,\mathbf{1}}^*$, if such exist.

Alternatively, if $LP$ has a bounded feasible region, one may consider 
\[
\min_{x,\xi}\{c^T x+M\xi: Ax+\widetilde{b}\xi-b\gamma=0,\:\mathbf{1}^Tx+\xi+\gamma=n+2,\:x,\xi,\gamma\geq0\},
\]
where $\widetilde{b}=b-A\mathbf{1}$ and $M>0$ is a large number. The
vector $d=\mathbf{1}$ is feasible for this problem; by the
corollary, the hyperbolic relaxation of the problem above is
bounded. Our new optimization problem corresponds to first taking a
standard big-$M$ formulation of $LP$ followed by homogenizing the
variables using $\gamma$ and normalizing all variables to a standard
simplex. Due to normalization, not all $x, \xi, \gamma$ may be
zeroed simultaneously. For large enough $M$ at the optimum $\xi=0$;
also $\gamma>0$, for otherwise $LP$ must have unbounded feasible
region. To complete our justification of (C), observe that a
solution to $LP$ may be easily recovered from the solution to the
problem above.

\begin{rem}
Here we want to draw the first parallel between the proposed
Shrink-Wrapping setting and path-following interior-point methods.
Later in Section~\ref{SecSWvsIPM} we discuss this relationship in
more details. Observe that assumption (B) combined with an
additional requirement that $d$ lies on the central path
corresponding to standard log-barrier $f(x)=-\ln\prod_{i=1}^n x_i$
implies that indeed we may choose $d$ with $HP_{r,d}$ bounded. In
turn, note that (B) combined with existence of strictly dual
$LP$-feasible $s: A^Ty+s=c, s\in\Re^n_{++}$, implies the existence
of the central path. To see how to choose such $d$, consider $LP$
\[
\min\{c^T x: Ax=b, x\in\Re^n_+\}
\]
together with its dual
\[
\max\{b^Ty: A^Ty+s=c, s\in\Re^n_+\}.
\]
Recall that the central path may be characterized as $x.\cdot
s=\mu\mathbf{1}, \:\mu>0$. So, if $d\in\Re^n_{++}$ is on the central
path, then $\mu\mathbf{1}./d$ is dual $LP$-feasible for some
$\mu>0$. Moreover, note that $\mu\mathbf{1}./d$ is an element of the
dual cone $\mathcal{K}_{r,d}$; follows from the cone inclusion
\[
\mathcal{K}_{n-1,d}\subseteq\mathcal{K}_{n-2,d}\subseteq\cdots\subseteq\mathcal{K}_{1,d}\subseteq\Re^n_+\subseteq\mathcal{K}_{1,d}^*
\subseteq\cdots\subseteq\mathcal{K}_{n-2,d}^*\subseteq\mathcal{K}_{n-1,d}^*
\]
where $\mathcal{K}_{n-1,d}^*$ consists of all nonnegative multiples
of $\mathbf{1}./d$. Consequently $\mu\mathbf{1}./d$ is feasible for
the dual conic problem to $HP_{r,d}$, and by conic duality
$HP_{r,d}$ is bounded. The described argument easily generalizes to
SOCP and SDP. If no such points $d$ is readily available, think
``self-dual embedding'' for symmetric cones; this gives yet another,
this time more theoretical justification for (C) -- observe that the
self-dual embedding will nearly double the sizes of the matrices we
have to work with if we were to consider Newton's like scheme based
on linearization of, say,~\ref{HP:KKT_p} below for tracing $x(d)$,
and thus will increase the amount of computations roughly
$2^3=8$-fold.
\end{rem}

By convexity of $K_{r,d}$ and assumption (D) it follows that KKT
conditions are both necessary and sufficient for optimality in
$HP_{r,d}$, and so the solution $x=x(d)$ is characterized by a
system of polynomial equations
\begin{eqnarray}\label{HP:KKT_p}
\left\{
\begin{array}{l}
\nabla E_{n-r}(x./d)+A^T y=\tau c,\: \tau>0,\\
E_{n-r}(x./d) = 0,\\
Ax=b
\end{array}
\right.
\end{eqnarray}
with $x\in\mathcal{K}_{r,d}$ and $y\in\Re^m$. To see why (D) implies
necessity of KKT conditions, note that $\nabla E_{n-r}(x./d)$ does
not vanish at $x(d)$, for otherwise we must have that
$E_{n-r-1}(x(d)./d)=0$ and so by Proposition~\ref{Hyper:MultRoot}
$x(d)\in\Re^n_+$ and is optimal for $LP$. Strictly speaking, to
characterize $x(d)$ in the above we need to add another set of
constraints ensuring $x(d)\in\mathcal{K}_{r,d}$, e.g.,
Corollary~\ref{Hyper:ConeBdry}; $p(x)=0$ alone does not suffice.

In addition to assumptions (A)-(E), we will be assuming that
\begin{itemize}
\item{(F)} $\rank(A)=m$,
\item{(G)} $LP$ solution $x^*$ is unique and has precisely $m$ non-zeros.
\end{itemize}
(F) is a standard convenient assumption commonly underlying the
interior-point methods and practically may be ensured by, say,
performing a QR factorization of $A$. (G) is a convenient assumption
that greatly simplifies the subsequent analysis; note that (G) is
\emph{generic} in a sense that it holds true almost surely for all
infinitesimal perturbations of the constraint vector $b$ by strict
complementarity for LP, implying that even if (G) fails at first, it
may be easily restored by slightly perturbing the original problem.
We hypothesize that in fact (G) may be lifted altogether, but for
the sake of compactness and readability of the manuscript we do not
attempt to verify the latter claim now.

Observe that if we fix $r = n-m-1$, $HP_{r,d}$ produces a
\emph{tight fit} relaxation to $LP$: any $LP$ vertex, including the
optimum, as a nonnegative solution to $Ax=b$ having at most $m$
non-zero entries, belongs to $\partial\mathcal{K}_{r,d}$, thus
hypothetically even allowing $x(d)=x^*$ for some well chosen $d$. In
Section~\ref{SecSWLoc} we will see that such a choice is indeed
possible.

Note that most of the observations we made so far may be extended
beyond LP to other hyperbolic optimization problems, such as SOCP
and SDP.

\begin{rem}
In characterizing the solution of $HP_{r,d}$, in particular the
boundary of $\mathcal{K}_{r,d}$, rather than relying on
$p(x)=r!E_n(d)E_{n-r}(x./d)$ hyperbolic with respect to $d$, one may
rely the concave on $\mathcal{C}'(d)$ \emph{ratio functional}
\[
q(x)=\frac{p(x)}{p'(x)}.
\]
Assuming (F) and using $q(x)$, Renegar has observed the existence of
the so-called \emph{central line} -- a strictly feasible line
segment $\mathcal{L}$ whose closure contains $x^*$ with an
additional property that if $d\in\mathcal{L}$ then $x(d)=x^*$;
moreover, turns out that the Jacobian of $x(d)$ for
$d\in\mathcal{L}$ has a very special structure that allows for a
nice geometric interpretation.

From computational point of view, fixing $r = n-m-1$, the usage of
$q(x)$ instead of $p(x)$ might help one to better address potential
numerical ill-conditioning when considering the gradient and Hessian
in linearized KKT for $x(d)$ such as~\ref{HP:KKT_p}, as $q(x)$ is
proportional to
\[
\frac{\prod_{i=1}^m \lambda_i(x)}{\sum_{i=1}^m\prod_{j\neq i}\lambda_j(x)}=\left(\sum_{i=1}^m\frac{1}{\lambda_i(x)}\right)^{-1}
\]
while $p(x)$ is proportional to $\prod_{i=1}^m\lambda_i(x)$, and
thus $q(x)$ suffers from the additive effect of simultaneously
zeroing more than one eigenvalue of $x$, while for $p(x)$ this
effect is multiplicative, when (G) is lost. Recall that at least one
eigenvalue of $x$ approaches 0 as $x$ nears the boundary of
$\mathcal{K}_{r,d}$.
\end{rem}

\subsection{Choice of dynamics}

We start by recalling that $d\in\Re^n_{++}$, but
$x(d)\notin\Re^n_{++}$. From $\mathcal{K}_{r,d}=\{x\in\Re^m:
x=d.\cdot z, z\in\mathcal{K}_{r,\mathbf{1}}\}$ and
$\Re^n_{++}\subset \mathcal{K}_{r,\mathbf{1}}\subset \Re^n$ one may
conclude that in general, the closer $d$ is to a vertex of $LP$, the
tighter the feasible region of $HP_{r,d}$ fits around that vertex.
This last informal observation suggests that given some initial
value $d_{(0)}$ of $d$, it might be beneficial to update
$d_{(0)}\mapsto d_{(1)}$ so that $d_{(1)}$ is closer to the solution
to $LP$, hoping that $x(d_{(1)})$ gets closer to $x^*$. In
particular, one could consider obtaining $d_{(1)}$ by moving from
$d_{(0)}$ towards $x^*$. Since $x^*$ is not known a priori, we may
choose the next best possible candidate, namely $x(d_{(0)})$ as a
surrogate for $x^*$ in the above.

The suggested dynamics for $d$ and $x(d)$ my be formalized through the ODE
\begin{eqnarray}\label{Dynamics:ODE}
\begin{array}{c}
\dot{d}=x(d)-d,\\
d|_{t=0}=d_{(0)}.
\end{array}
\end{eqnarray}
Although, we chose this particular dynamics to govern the behavior of $d$ and $x(d)$, many
other choices are possible. We are interested in studying the continuous trajectories of $d(t),\:t\in[0,\infty)$,
where $d(t)$ solves~\ref{Dynamics:ODE}.

The following statement was conjectured by Renegar: ``under $LP$
strict dual feasibility $d(t)$ converges to $x^*$''; for more
details see the very recent~\cite{Renegar:SW}. We refine it by
observing that $x(d)$ might not even be defined if $d$ is chosen
poorly, i.e., $HP_{r,d}$ is unbounded; for convenience define
$x(\mathbf{0})=\mathbf{0}$.
\begin{thm}
If for all $t\geq0$
we have bounded $HP_{r,d(t)}$,
then $d(t)\rightarrow x^*$ as $t\rightarrow\infty$.
\end{thm}
\begin{proof}
Follows from $c^T\dot{d}=c^T(x(d)-d)<0$ since $HP_{r,d}$ is a relaxation of $LP$.
\end{proof}
\noindent We hypothesize that indeed for $HP_{r,d(t)}$ to stay
bounded for all $t\geq0$ it suffices to choose initial value
$d_{(0)}$ corresponding to bounded $HP_{r,d_{(0)}}$, in which case
both $d(t)$ and $x(d(t))$ converge to $x^*$ as $t\rightarrow\infty$.

Note that with choice of affine feasible $d_{(0)}$, the trajectory
$d(t)$ remains affine feasible for all $t\geq0$. Let columns of $B$
form a basis for $\nul(A)$. Any affine feasible point $d$ can be
written as $d=d_{(0)}+B\delta$ for some $\delta\in\Re^{n-m}$.
ODE~\ref{Dynamics:ODE} may be re-written as
\begin{eqnarray*}
\begin{array}{c}
B\dot{\delta}=d_{(0)}+B\xi(\delta)-d_{(0)}-B\delta,\\
B\delta|_{t=0}=d-d_{(0)},
\end{array}
\end{eqnarray*}
where $x(d)=d_{(0)}+B\xi(\delta)$. Since $B$ is injective, the above
is equivalent to
\begin{eqnarray*}
\begin{array}{c}
\dot{\delta}=\xi(\delta)-\delta,\\
\delta|_{t=0}=\delta_{(0)}.
\end{array}
\end{eqnarray*}
Consequently, we can pick an arbitrary affine coordinate system of
$\{x\in\Re^n: Ax=b\}$ to analyze~\ref{Dynamics:ODE}. Thus, we call
ODE~\ref{Dynamics:ODE} \emph{affine invariant}.

A corresponding discrete algorithm may be based on approximating the
trajectories $d(t)$ and $x(d(t))$ iteratively, generating a sequence
of pairs, $ (d_i,x_i) $, $ i = 1, \ldots, \infty $:
\begin{itemize}
\item given $d_i$, compute $x_i\approx x(d_i)$,
\item set $d_{i+1}=d_i+\alpha(x_i-d_i)$ for some properly chosen $\alpha\in(0,1)$, iterate.
\end{itemize}
Note that characterization~\ref{HP:KKT_p} suggests a way to trace
$x(d)$ when $d$ changes ever so slightly using, for example,
Newton's method.

\begin{rem}
Although, as we will see in Section~\ref{SecSWvsIPM}, trajectories
$d(t)$ generalize the notion of the central path, there appears to
be no known analogue for $x(d(t))$. In our limited numerical
experiments $x(d(t))$ typically converged to $x^*$ much sooner than
$d(t)$, which suggests that algorithmically one might focus on
tracing $x(d(t))$.
\end{rem}

In the subsequent section we formally introduce the notion of the
central line $\mathcal{L}$ which acts as an \emph{invariant} set
w.r.t. dynamics of the Shrink-Wrapping iterates $d$ (and $x(d)$) --
invariant in a sense that if $d(t_0)\in\mathcal{L}$ for some $t_0$,
then $d(t)\in\mathcal{L}$ for all $t\geq t_0$. Next, we devote our
attention to studying Shrink-Wrapping trajectories near the central
line; in turn, choosing the neighborhood of this invariant set
properly will enable us to lift the extra assumption on
$HP_{r,d(t)}$ to stay bounded for all $t\geq0$. We observe the
special structure of the Hessian of $x(d)$ where $d$ belongs to
$\mathcal{L}$ relying only on polynomials in characterization of
solution to $HP_{r,d}$. The latter allows us to significantly
simplify the analysis of Shrink-Wrapping trajectories in the
neighborhood of $\mathcal{L}$, as compared
to~\cite{Zinchenko:thesis}, and show that the central line acts as
\emph{attractor} set for $d$, provided the initial iterate $d_{(0)}$
was chosen significantly close to $\mathcal{L}$.

\section{On Shrink-Wrapping trajectories}\label{SecSWLoc}
\subsection{Invariant central line}

For a set of indices $\mathcal{B}$ let $x_{\mathcal{B}}$ denote a
vector with coordinates $x_i,\:i\in\mathcal{B}$, e.g.,
\[
\left(\begin{array}{c} 1\\2\\5\\6\end{array}\right)_{\{3,4\}}=
\left(\begin{array}{c} 5\\6\end{array}\right).
\]
We also simply write $x_{-i}$ when we want to obtain a vector from
$x$ of dimension one less by dropping $i^{th}$ coordinate. Let
$x.^2=x.*x$; as usual, power operation takes precedence over
multiplication or division. Component-wise vector operations take
precedence over standard operations on vectors, and otherwise occur
in order of appearance. For two vector-valued functions
$x,y:\mathcal{P}\rightarrow\Re^n$, we write $x=O(y)$ if there is a
constant $K>0$ such that $|x_i|\leq K|y_i|,\:\forall i=1,n,\:\forall
p\in\mathcal{P}$. $[A; B]$ denotes vertical block-matrix consisting
of $A,B$:
\[
[A; B]=\left(\begin{array}{c}A\\B\end{array}\right).
\]
For a matrix $A$ we use $A_{:,i}$ to denote its $i^{th}$ column.

From now on, fix $r=n-m-1$ in $HP_{r,d}$; note $x^*$ belongs to the
boundary of $\mathcal{K}_{r,d}$. Given (E), without loss of
generality we may assume the last $m$ coordinates ox $x^*$ to be
non-zero. We choose the following parametrization of
$\{x\in\Re^n:Ax=b\}$: let $\xi$ denote first $n-m$ components of
affine feasible $x$. Note that $\xi=\mathbf{0}$ at the $LP$ optimum
$x^*$. Fixing $\mathcal{B}=\{n-m+1,n-m+2,\ldots,n\}$,
$x_{\mathcal{B}}\in\Re^m$ corresponds to last $m$ components of $x$;
note $x^*_{\mathcal{B}}$ is a vector of (non-zero) basic components
of $x^*$, compare this notation with the example in
Section~\ref{SecSWIntro}. Similarly, we denote $\delta$ to be first
$n-m$ components of hyperbolicity direction vector $d$, and
$d_\mathcal{B}$ its last $m$ components.

Note that if $x$ is affine feasible, we may re-write $Ax=b$ as
$x_\mathcal{B}=\widetilde{A}\xi+x^*_\mathcal{B}$ for some
$\widetilde{A}\in\Re^{m\times n-m}$; similarly
$d_\mathcal{B}=\widetilde{A}\delta+x^*_\mathcal{B}$. Recall that
ODE~\ref{Dynamics:ODE} is affine invariant; thus, to understand
$d(t)$ we may equivalently analyze trajectories $\delta(t)$ of
\begin{eqnarray}\label{Dynamics:ODEnm}
\begin{array}{c}
\dot{\delta}=\xi(\delta)-\delta,\\
\delta|_{t=0}=\delta(0).
\end{array}
\end{eqnarray}

\begin{defn} An open linear segment $\mathcal{L}\subset\{x\in\Re^n_{++}:Ax=b\}$
whose closure contains $x^*$ is called \emph{the central line} if
for any $d\in\mathcal{L}$ we have $x(d)=x^*$, and $\mathcal{L}$ is
not a proper subset of any other linear segment with the above
properties.
\end{defn}
Since we will be mostly working with first $n-m$ coordinate
parametrization of $\{x\in\Re^n:Ax=b\}$ and, in particular,
ODE~\ref{Dynamics:ODEnm}, we allow for a minor abuse of notation by
using the same symbol $\mathcal{L}$ for the central line wether when
referring to a subset of $\Re^n$ as in the definition above, or its
projection onto first $n-m$ coordinates.

\begin{prop}\label{CL:exists}
The central line exists.
\end{prop}
\begin{proof}
Rewriting the first equation of conditions~\ref{HP:KKT_p} in the
first $n-m$ coordinates
\begin{eqnarray*}
[I; \widetilde{A}]^T \Diag(\mathbf{1}./d) \nabla_z
E_{m+1}(z)|_{z=x./d}=\tau [I; \widetilde{A}]^T c,
\end{eqnarray*}
and observing that $\left(\nabla_z E_{m+1}(z)\right)_i=E_m(z_{-i}),
i=1,\ldots,n$, at $x=x^*$ we get
\begin{eqnarray}\label{CL:grad}
E_m(d_\mathcal{B})\:\mathbf{1}./\delta=\tau [I; \widetilde{A}]^T c,
\end{eqnarray}
recalling that $x^*$ has only last $m$ non-zeros.

Since $x^*$ is a unique minimizer for $LP$, we must have that
\[
\xi^T [I; \widetilde{A}]^T c>0,\:\forall\xi\in\Re^{n-m}_+,
\]
and so $[I; \widetilde{A}]^T c\in\Re^{n-m}_{++}$ by elementary LP
conic duality.

Set $\widetilde{\delta}=\mathbf{1}./([I; \widetilde{A}]^T
c)\in\Re^{n-m}_{++}$ and observe that $\Delta_{\max}>0$ may be
chosen such that the linear segment
$\mathcal{L}=\{d\in\Re^n_{++}:d=x^*+[I;\widetilde{A}]\widetilde{\delta}\cdot\Delta,\:\Delta\in(0,\Delta_{\max})\}$
is the largest possible. In turn, any $\delta$
corresponding to $d\in\mathcal{L}$ will result in positive $\tau$
in~\ref{CL:grad}; moreover, clearly
$x^*\in\partial\mathcal{K}_{r,d}$ since, again, $x^*$ has precisely
$m$ zeros, and so~\ref{HP:KKT_p} is satisfied at $x^*$. Therefore,
for any $d\in\mathcal{L}$ we have $x(d)=x^*$.
\end{proof}
Due to affine invariance of ODE~\ref{Dynamics:ODE} for the purpose
of its analysis, without loss of generality, we assume that $[I;
\widetilde{A}]^Tc=\mathbf{1}$ and consequently
$\delta=\Delta\mathbf{1}, \Delta>0$, when $d\in\mathcal{L}$: if not,
simply re-scale the first $n-m$ coordinates accordingly.


Observe that the elementary symmetric polynomial in
$x=[\xi;x_{\mathcal{B}}]\in\Re^n$ satisfies
\begin{eqnarray}\label{ESP:decomp}
E_{m+1}(x)=E_1(\xi)E_m(x_{\mathcal{B}})+E_2(\xi)E_{m-1}(x_{\mathcal{B}})+
E_3(\xi)E_{m-2}(x_{\mathcal{B}})+\cdots.
\end{eqnarray}

\begin{prop}\label{xi:Jacobian}
The Jacobian of $\xi(\delta)$ for $\delta\in\mathcal{L}$ is of the
form
\begin{eqnarray*}
J_{\xi}(\delta)=-\frac{E_m(x^*_{\mathcal{B}}./d_{\mathcal{B}})}{E_{m-1}(x^*_\mathcal{B}./d_{\mathcal{B}})}
\left(I-\frac{\mathbf{1}\mathbf{1}^T}{n-m}\right).
\end{eqnarray*}
\end{prop}
\begin{proof}
In order to compute the derivative of $\xi(\delta)$ for
$\delta\in\mathcal{L}$ we implicitly differentiate~\ref{HP:KKT_p}.
Consider a vector $\dot{\xi}_{\delta_1}$ of partial derivatives
$\frac{\partial\xi(\delta)}{\partial\delta_1}$ -- the first column
of $J_{\xi}(\delta)^T$.

Differentiating second equation of~\ref{HP:KKT_p} and re-writing it
terms of $\delta,\xi$, recalling that~\ref{CL:grad} implies
$\nabla_{\xi}E_{m+1}(x./d)|_{\xi=\mathbf{0}}$ is a positive multiple
of $\mathbf{1}$ since $\delta$ is also a positive multiple of
$\mathbf{1}$, it follows that $\dot{\xi}_{\delta_1}$ is orthogonal
to $\mathbf{1}$, that is,
\[
\mathbf{1}^T\: \dot{\xi}_{\delta_1}=0.
\]

In order to differentiate the first equation in~\ref{HP:KKT_p}, we
first revisit the expression of the gradient of $E_{m+1}(x./d)$ in
coordinates $\delta,\xi$: note that~\ref{ESP:decomp} implies
\begin{eqnarray}\label{CL:grad}
\begin{array}{rl}
\nabla_{\xi}E_{m+1}(x./d) = &
\nabla_{\xi}E_{m+1}([\xi;x_{\mathcal{B}}]./[\delta;d_{\mathcal{B}}])\\
= & \nabla_{\xi}\left(
E_1(\xi./\delta)E_m(x_{\mathcal{B}}./d_{\mathcal{B}}) +
E_2(\xi./\delta)E_{m-1}(x_{\mathcal{B}}./d_{\mathcal{B}}) +
\cdots \right)\\
= &
\Diag(1./\delta)\mathbf{1}\: E_m(x_{\mathcal{B}}./d_{\mathcal{B}}) +\\
& E_1(\xi./\delta)\:
\widetilde{A}^T\Diag(1./d_{\mathcal{B}})\nabla_z
E_m(z)|_{z=x_{\mathcal{B}}./d_{\mathcal{B}}} +\\
& \Diag(1./\delta)\nabla_z E_2(z)|_{z=\xi./\delta}\:
E_{m-1}(x_{\mathcal{B}}./d_{\mathcal{B}}) +\\
& E_2(\xi./\delta)\:
\widetilde{A}^T\Diag(1./d_{\mathcal{B}})\nabla_z
E_{m-1}(z)|_{z=x_{\mathcal{B}}./d_{\mathcal{B}}} + \cdots.
\end{array}
\end{eqnarray}
Differentiating the above and evaluating at $\xi=\mathbf{0}$ we get
\[
\begin{array}{rl}
\frac{\partial}{\partial\delta_1} \left(\Diag(1./\delta)\mathbf{1}\:
E_m(x_{\mathcal{B}}./d_{\mathcal{B}})\right) = &
\Diag([-1./\delta_1^2;\mathbf{0}])\mathbf{1}\:E_m(x^*_{\mathcal{B}}./d_{\mathcal{B}}) +\\
& \Diag(1./\delta)\mathbf{1}\: \nabla_z
E_m(z)|_{z=x^*_{\mathcal{B}}./d_{\mathcal{B}}}^T
\left((\widetilde{A}\dot{\xi}_{\delta_1})./d_{\mathcal{B}}\right) +\\
& \Diag(1./\delta)\mathbf{1}\: \nabla_z
E_m(z)|_{z=x^*_{\mathcal{B}}./d_{\mathcal{B}}}^T
\left(-x^*_{\mathcal{B}}./d_{\mathcal{B}}.^2.\cdot\widetilde{A}_{:,1}\right),
\end{array}
\]
\[
\frac{\partial}{\partial\delta_1} \left(E_1(\xi./\delta)\:
\widetilde{A}^T\Diag(1./d_{\mathcal{B}})\nabla_z
E_m(z)|_{z=x_{\mathcal{B}}./d_{\mathcal{B}}}\right) =\\
\mathbf{1}^T\left(\dot{\xi}_{\delta_1}./\delta\right)\:
\widetilde{A}^T\Diag(1./d_{\mathcal{B}}) \nabla_z
E_m(z)|_{z=x^*_{\mathcal{B}}./d_{\mathcal{B}}},
\]
\[
\frac{\partial}{\partial\delta_1} \left(\Diag(1./\delta)\nabla_z
E_2(z)|_{z=\xi./\delta}\: E_{m-1}(x_{\mathcal{B}}./d_{\mathcal{B}})\right) =\\
\Diag(1./\delta)
(\mathbf{1}\mathbf{1}^T-I)\Diag(1./\delta)\dot{\xi}_{\delta_1}\:
E_{m-1}(x^*_{\mathcal{B}}./d_{\mathcal{B}}),
\]
and
\[
\frac{\partial}{\partial\delta_1} \left(E_2(\xi./\delta)\:
\widetilde{A}^T\Diag(1./d_{\mathcal{B}})\nabla_z
E_{m-1}(z)|_{z=x_{\mathcal{B}}./d_{\mathcal{B}}}\right)=\mathbf{0}
\]
with all the remaining ``higher-order'' in $\xi$ terms in the
expression for $\frac{\partial}{\partial\delta_1}
\nabla_{\xi}E_{m+1}(x./d)$ being zero. As a result, at
$\delta=\Delta\mathbf{1}\in\mathcal{L}$ where $\Delta>0$ and
corresponding $\xi=\mathbf{0}$, the first equation of
differentiated~\ref{HP:KKT_p} becomes
\begin{eqnarray}\label{CL:Hessian}
\begin{array}{rl}
\frac{\partial\tau}{\partial\delta_1} [I; \widetilde{A}]^Tc =
& [-1/\Delta^2;\mathbf{0}]\: E_m(x^*_\mathcal{B}./d_{\mathcal{B}}) +\\
& \frac{1}{\Delta}
\nabla_z E_m(z)|_{z=x^*_{\mathcal{B}}./d_{\mathcal{B}}}^T
\left((\widetilde{A}\dot{\xi}_{\delta_1})./d_{\mathcal{B}} -
x^*_{\mathcal{B}}./d_{\mathcal{B}}.^2.\cdot\widetilde{A}_{:,1}
\right)
\mathbf{1} -\\
& \frac{1}{\Delta^2}\dot{\xi}_{\delta_1}\:
E_{m-1}(x^*_\mathcal{B}./d_{\mathcal{B}}),
\end{array}
\end{eqnarray}
observing the cancelations due to the orthogonality condition
$\mathbf{1}^T\: \dot{\xi}_{\delta_1}=0$; note that the affine
feasibility requirement is satisfied by the choice of coordinates.

Finally, to compute $\dot{\xi}_{\delta_1}$ we need to
solve~\ref{CL:Hessian} together with $\mathbf{1}^T\:
\dot{\xi}_{\delta_1}=0$ -- a system of $n-m+1$ equations in $n-m+1$
variables $\dot{\xi}_{\delta_1},
\frac{\partial\tau}{\partial\delta_1}$.
Pre-multiplying both sides of the expression by $\mathbf{1}^T$,
recalling $[I;\widetilde{A}]^Tc=\mathbf{1}$, we obtain
\[
\frac{\partial\tau}{\partial\delta_1}  =
\frac{-E_m(x^*_\mathcal{B}./d_{\mathcal{B}})}{(n-m)\Delta^2} +
\frac{\omega}{\Delta}.
\]
where
\[
\omega=\nabla_z E_m(z)|_{z=x^*_{\mathcal{B}}./d_{\mathcal{B}}}^T
\left((\widetilde{A}\dot{\xi}_{\delta_1})./d_{\mathcal{B}} -
x^*_{\mathcal{B}}./d_{\mathcal{B}}.^2.\cdot\widetilde{A}_{:,1}
\right).
\]
Now, using the expression for $\frac{\partial\tau}{\partial\delta_1}$
we may re-write~\ref{CL:Hessian} as
\[
\frac{-E_m(x^*_\mathcal{B}./d_{\mathcal{B}})}{(n-m)\Delta^2}\:\mathbf{1}
= [-1/\Delta^2;\mathbf{0}]\: E_m(x^*_\mathcal{B}./d_{\mathcal{B}}) -
\frac{1}{\Delta^2}\dot{\xi}_{\delta_1}\:
E_{m-1}(x^*_\mathcal{B}./d_{\mathcal{B}}),
\]
resulting in
\[
\dot{\xi}_{\delta_1}=-\frac{E_m(x^*_\mathcal{B}./d_{\mathcal{B}})}
{E_{m-1}(x^*_\mathcal{B}./d_{\mathcal{B}})}\left(e_{(1)}-\frac{\mathbf{1}}{n-m}\right)
\]
where $e_{(1)}=[1;\mathbf{0}]\in\Re^{n-m}$ is the first unit vector.

Similarly, we derive the expressions for
$\frac{\partial\xi(\delta)}{\partial\delta_i},\: i=2,\ldots,n-m.$
\end{proof}

The Jacobian of $\xi(\delta)$ for $\delta\in\mathcal{L}$ may be
interpreted as a negative projection onto the null space of
$\mathbf{1}^T$ with a corresponding multiple
$\frac{E_m(x^*_\mathcal{B}./d_{\mathcal{B}})}
{E_{m-1}(x^*_\mathcal{B}./d_{\mathcal{B}})}>0$ -- recall that
$d=[\delta;d_{\mathcal{B}}]\in\Re^n_{++}$ for
$\delta\in\mathcal{L}$; note that
$\frac{E_m(x^*_\mathcal{B}./d_{\mathcal{B}})}
{E_{m-1}(x^*_\mathcal{B}./d_{\mathcal{B}})}$ is finite for all
$\delta\in\mathcal{L}$. In turn, this implies that, up to first
order, a small deviation of $\delta$ from $\mathcal{L}$ in the
direction orthogonal to $\mathbf{1}$, that is, orthogonal to the
central line, results in the displacement of the corresponding
$\xi(\delta)$ in precisely the opposite direction, see
Figure~\ref{Hyper:Fig5}.

The last observation suggests that $\mathcal{L}$ might be an
attractor set: when considering the dynamics
of~\ref{Dynamics:ODEnm}, note that small deviations of $\delta$ away
form $\mathcal{L}$ appear to be counter-acted by corresponding
changes in $\xi(\delta)$ away from $\mathbf{0}$, thus, forcing
$\delta(t)$ to cross-over the central line. In what follows we will
see that indeed this is the case.

\subsection{Trajectories near central line}





It is convenient to introduce the following orthogonal decomposition
of $\delta\in\Re^{n-m}_{++}$:
\[
\delta=\delta_\parallel+\delta_\perp, \mbox{ where }
\delta_\parallel=\Delta\mathbf{1}, \Delta>0, \mbox{ and
}\mathbf{1}^T\delta_\perp=0.
\]

Intuitively, if
\[
\left\{
\begin{array}{lll}
\dot{\delta}_\parallel & \approx & -1\cdot\delta_\parallel,\\
\dot{\delta}_\perp & \approx & -(1+\theta)\cdot\delta_\perp,
\end{array}\right.
\]
for some $\theta>0$ and the approximation above is ``accurate
enough'', we expect
\[
\left\{
\begin{array}{lll}
\delta_\parallel(t) & \approx & \delta_\parallel(0)\cdot e^{-t},\\
\delta_\perp(t) & \approx & \delta_\perp(0)\cdot e^{-(1+\theta)t},
\end{array}\right.
\]
and so
\[
\left\|\frac{\delta_\perp(t)}{\delta_\parallel(t)}\right\|\approx
\left\|\frac{\delta_\perp(0)}{\delta_\parallel(0)}\right\|\cdot
e^{-\theta t}\rightarrow0\mbox{ as }t\rightarrow\infty.
\]
Note that the Jacobian of $\xi(\delta)$ for $\delta\in\mathcal{L}$
suggests that the system~\ref{Dynamics:ODEnm} indeed assumes the
form of ODE as above, at least in some vicinity of the central line.
However, as we witness in this subsection, although our intuition
proves to be correct, we have to be quite careful
since $\xi(\delta)$ governing~\ref{Dynamics:ODEnm} might easily fail
to be differentiable at $\delta=\mathbf{0}$.

For the next lemma we allow for a slight abuse of notation using
$\delta_\parallel$ to denote the first coordinate $\delta_1$ of a
vector $\delta$, and $\delta_\perp$ to denote the vector of the
remaining coordinates $\delta_{-1}$ in some orthonormal basis; note
that this is consistent with, say, equipping $\Re^{n-m}$ with a
system of orthonormal coordinates where the first coordinate axis is
aligned with $\mathbf{1}$. The quality of the approximation in the
ODE above that suffices for our purposes may be characterized by the
following statement.
\begin{lem}
Let $\delta$ be governed by the following ODE with locally Lipschitz
continuous right-hand side
\begin{eqnarray*}
\left\{
\begin{array}{lll}
\dot{\delta}_\parallel & = & -1\cdot\delta_\parallel+
O\left(\frac{\|\delta_\perp\|^2}{\delta_\parallel}\mathbf{1}\right),\\
\dot{\delta}_\perp & = & -(1+\theta)\cdot\delta_\perp+
O\left(\frac{\|\delta_\perp\|^2}{\delta_\parallel}\mathbf{1}\right),
\end{array}\right.
\end{eqnarray*}
for some $\theta>0$. Then for any fixed $\Delta_1>0$ there exists
$\epsilon>0$ such that for any initial
$\delta(0)=\delta_\parallel(0)+\delta_\perp(0)$ in the \emph{central
wedge}
\[
\mathcal{W}=\{\delta:\|\delta_\perp\|<\delta_\parallel\cdot\epsilon
\mbox{ and }\Delta\in(0,\Delta_1) \}
\]
we have
\[
\delta_\parallel(t)\leq \delta_\parallel(0)\: e^{-\nu t}
\]
for some $\nu>0$, and
\[
\left\|\frac{\delta_\perp(t)}{\delta_\parallel(t)}\right\|\leq
\left\|\frac{\delta_\perp(0)}{\delta_\parallel(0)}\right\|\: e^{-\omega t},
\]
for some $\omega>0$, and so $\delta(t)\in\mathcal{W}$ for all
$t\geq0$; moreover, $\omega\rightarrow\theta, \nu\rightarrow 1$ as
$\left\|\frac{\delta_\perp(0)}{\delta_\parallel(0)}\right\|\rightarrow0$.
\end{lem}
\begin{proof}
Since the right-hand side of the ODE above is locally Lipschitz
continuous, the unique and continuously-differentiable solution
$\delta(t)$ exists for any choice of initial
$\delta_\parallel(0)\neq0$ and arbitrary $\delta_\perp(0)$, at least
on some open interval of $t$ containing 0. Consider
\[
\frac{1}{2}\:\frac{d}{dt}\:\left\|\frac{\delta_\perp}{\delta_\parallel}\right\|^2
= \left(\frac{\dot{\delta}_\perp
\delta_\parallel-\dot{\delta}_\parallel
\delta_\perp}{\delta_\parallel^2}\right)^T \frac{\delta_\perp}{\delta_\parallel}
= -\theta\:\frac{{\delta_\perp}^T\:\delta_\perp}{\delta_\parallel^2}
-\frac{\delta_\perp^T
O\left(\frac{\|\delta_\perp\|^2}{\delta_\parallel}\mathbf{1}\right)}{\delta_\parallel^2}
-O\left(\frac{\|\delta_\perp\|^2}{\delta_\parallel}\right)\:
\frac{{\delta_\perp}^T\:\delta_\perp}{\delta_\parallel^3}
\]
and note that by Cauchy-Schwarz inequality
\[
\omega\:\left.\left\|\frac{\delta_\perp}{\delta_\parallel}\right\|^2\right|_{t=0}
\leq-\left.\frac{1}{2}\:\frac{d}{dt}\:\left\|\frac{\delta_\perp}{\delta_\parallel}\right\|^2\right|_{t=0},
\]
where
\[
\omega=\left.\left(\theta
-K_1\left\|\frac{\delta_\perp}{\delta_\parallel}\right\|
-K_2\left\|\frac{\delta_\perp}{\delta_\parallel}\right\|^2\right)\right|_{t=0}
\]
and the constants $K_1,K_2>0$ satisfy
\[
\left\|O\left(\frac{\|\delta_\perp\|^2}{\delta_\parallel}\mathbf{1}\right)\right\|\leq K_1\:
\frac{\|\delta_\perp\|^2}{|\delta_\parallel|}
\]
and
\[
\left|O\left(\frac{\|\delta_\perp\|^2}{\delta_\parallel}\right)\right|\leq K_2\:
\frac{\|\delta_\perp\|^2}{|\delta_\parallel|}.
\]
Clearly, $\omega>0$ provided
$\left\|\frac{\delta_\perp(0)}{\delta_\parallel(0)}\right\|$ is
small enough, e.g., $\delta(0)\in\mathcal{W}$ for sufficiently small
$\epsilon$.

Continuity of $\delta(t)$ implies
$\frac{d}{dt}\:\left\|\frac{\delta_\perp}{\delta_\parallel}\right\|^2<0$
for $t\in[0,\tau)$ for some $\tau>0$, and so
$\left\|\frac{\delta_\perp}{\delta_\parallel}\right\|^2$ is
decreasing on $[0,\tau)$; therefore, if $\delta(0)\in\mathcal{W}$
and, in addition, $\delta_\parallel(t)$ is non-increasing, then
$\delta(t)\in\mathcal{W}$ for all $t\in[0,\tau)$. Moreover, for any
fixed $\kappa>0$ we may choose $\tau$  such that for all
$t\in[0,\tau]$ we have
\[
0<\frac{\omega}{1+\kappa}\:\left.\left\|\frac{\delta_\perp}{\delta_\parallel}\right\|^2\right|_{t=0}\leq -\frac{1}{2}\:\frac{d}{dt}\:\left\|\frac{\delta_\perp}{\delta_\parallel}\right\|^2,
\]
and so
\[
\left.\left\|\frac{\delta_\perp}{\delta_\parallel}\right\|^2\right|_{t=\tau}-
\left.\left\|\frac{\delta_\perp}{\delta_\parallel}\right\|^2\right|_{t=0}=
\int_{t=0}^\tau \frac{d}{dt}\:\left\|\frac{\delta_\perp}{\delta_\parallel}\right\|^2\:dt\leq -\frac{2\omega}{1+\kappa}\:\left.\left\|\frac{\delta_\perp}{\delta_\parallel}\right\|^2\right|_{t=0}\:\tau.
\]

Similarly, differentiating $\delta_\parallel^2$ and choosing $\epsilon$ small
enough in $\mathcal{W}\ni\delta(0)$, we may show that for sufficiently small $\tau>0$ we have
\[
-\frac{1}{2}\:\frac{d}{dt}\:\delta_\parallel^2=-\dot{\delta}_\parallel\:\delta_\parallel=
\delta_\parallel^2-
\delta_\parallel\:O\left(\frac{\|\delta_\perp\|^2}{\delta_\parallel}\mathbf{1}\right)\geq
\delta_\parallel^2\:\left(1-K_2\left\|\frac{\delta_\perp}{\delta_\parallel}\right\|^2\right)>0
\]
for $t\in[0,\tau]$, and so $\delta_\parallel$ is monotone-decreasing
on $[0,\tau]$ implying $\delta(t)\in\mathcal{W}$, and
\[
\left.\delta_\parallel^2\right|_{t=\tau}-
\left.\delta_\parallel^2\right|_{t=0}\leq
-\frac{2\nu}{1+\kappa}\:\left.\delta_\parallel^2\right|_{t=0}\:\tau
\]
where
\[
\nu=\left.\left(1-K_2\left\|\frac{\delta_\perp}{\delta_\parallel}\right\|^2\right)\right|_{t=0}.
\]

Noting that since $\delta(\tau)\in\mathcal{W}$, the argument above may be repeated
at $t=\tau$ treating it as $t=0$, we observe that the solution $\delta(t)\in\mathcal{W}$
with the above properties may be extended to any $t\geq0$. Finally, it is left to recognize
the exponents in the bounds for
$\left.\left\|\frac{\delta_\perp}{\delta_\parallel}\right\|^2\right|_{t=\tau}-
\left.\left\|\frac{\delta_\perp}{\delta_\parallel}\right\|^2\right|_{t=0}$
and
$\left.\delta_\parallel^2\right|_{t=\tau}-
\left.\delta_\parallel^2\right|_{t=0}$
while letting $\tau\rightarrow0$, followed by $\kappa\rightarrow0$.
\end{proof}

Observing that our proof relies on the big-$O$ form of the ODE only
in some central wedge $\mathcal{W}$, that $\theta$ only needs to be
bounded away from 0 on $\mathcal{W}$, and choosing the coordinate
system for $\delta\in\Re^{n-m}$ so that the first coordinate is
aligned with $\delta_\parallel$, we may state the following result;
recall $\delta=\delta_\parallel+\delta_\perp$ forms orthogonal
decomposition of $\delta$.
\begin{cor}\label{ODE:mainCor}
If in some \emph{central wedge} 
\[
\widetilde{\mathcal{W}}=\{\delta\in\Re^{n-m}:\delta=\delta_\parallel+\delta_\perp,
\delta_\parallel=\Delta\mathbf{1}, \mathbf{1}^T\delta_\perp=0,
\|\delta_\perp\|<\|\delta_\parallel\|\cdot\widetilde{\epsilon},
\mbox{ and }\Delta\in(0,\Delta_1) \}
\]
$\xi(\delta)$ is locally Lipschitz continuous, and the
ODE~\ref{Dynamics:ODEnm} may be re-written as
\begin{eqnarray}\label{Dynamics:ODEtarget}
\dot{\delta}=-\delta-\theta\cdot\delta_\perp+O\left(\frac{\|\delta_\perp\|^2}{\|\delta_\parallel\|}\:\mathbf{1}\right)
\end{eqnarray}
where $\theta=\theta(\delta)>0$ is bounded away from 0 on
$\widetilde{\mathcal{W}}$, then there is a possibly smaller central
wedge $\mathcal{W}\subseteq\widetilde{\mathcal{W}}$ corresponding to
$\epsilon\leq\widetilde{\epsilon}$,
\[
\mathcal{W}=\{\delta\in\Re^{n-m}:\delta=\delta_\parallel+\delta_\perp,
\delta_\parallel=\Delta\mathbf{1}, \mathbf{1}^T\delta_\perp=0,
\|\delta_\perp\|<\|\delta_\parallel\|\cdot\epsilon, \mbox{ and
}\Delta\in(0,\Delta_1) \},
\]
such that for some fixed $\nu,\omega>0$ we have
\[
\begin{array}{ccc}
\vspace{2mm}\left\|\delta_\parallel(t)\right\| & \leq & \left\|\delta_\parallel(0)\right\|\: e^{-\nu t},\\
\left\|\frac{\delta_\perp(t)}{\delta_\parallel(t)}\right\| & \leq &
\left\|\frac{\delta_\perp(0)}{\delta_\parallel(0)}\right\| \:
e^{-\omega t},
\end{array}
\]
for any $\delta(0)\in\mathcal{W}$.
\end{cor}

We say that $\delta$ \emph{converges exponentially} to $\mathcal{L}$
if $\left\|\frac{\delta_\perp(t)}{\delta_\parallel(t)}\right\|\leq
\left\|\frac{\delta_\perp(0)}{\delta_\parallel(0)}\right\|\:
e^{-\omega t}, \omega>0$, and
$\|\delta(t)\|\leq\|\delta(0)\|\:e^{-\eta t}, \eta>0$, see
Figure~\ref{Hyper:Fig5}; we say that $\xi=\xi(\delta)$
\emph{converges exponentially} to $\mathbf{0}$ if
$\|\xi(t)\|\leq\|\xi(0)\|\:e^{-\varpi t}, \varpi>0$. The
\textbf{main result} of this section is as follows.
\begin{thm}\label{Dynamics:MainTheorem}
For any $\Delta_1\in(0,\Delta_{\max})$, there is a corresponding
central wedge $\mathcal{W}$ such that if $\delta(0)\in\mathcal{W}$
then $\delta(t)$ converges exponentially to the central line
$\mathcal{L}$. Moreover, the corresponding $\xi(t)=\xi(\delta(t))$ converges
exponentially to $\mathbf{0}$.
\end{thm}
\noindent Observe that to prove the theorem, by
Corollary~\ref{ODE:mainCor} it is left to exhibit that indeed the
ODE~\ref{Dynamics:ODEnm} may be written in the
form~\ref{Dynamics:ODEtarget} in some central wedge
$\widetilde{\mathcal{W}}$.

We start by investigating the behavior of $\xi(\delta)$ for $\delta$
\emph{near} the central line. Recall that the necessary and
sufficient conditions~\ref{HP:KKT_p} for $x(d)$ may be re-written in
terms of $\delta,\xi$ to characterize $\xi(\delta)$ by
\[
\left\{\begin{array}{l}
\nabla_{\xi}E_{m+1}([\xi;x_{\mathcal{B}}]./[\delta;d_{\mathcal{B}}]) - \tau\mathbf{1}=\mathbf{0},\\
E_{m+1}([\xi;x_{\mathcal{B}}]./[\delta;d_{\mathcal{B}}])=0,
\end{array}\right.
\]
where the expression for
$\nabla_{\xi}E_{m+1}([\xi;x_{\mathcal{B}}]./[\delta;d_{\mathcal{B}}])$
may be found in~\ref{CL:grad}, $\tau>0$, and, additionally,
$[\xi;x_{\mathcal{B}}]\in\mathcal{K}_{n-m-1,d}$ captured, for
example, via Corollary~\ref{Hyper:ConeBdry}; note that $\tau>0$
guarantees that $\xi$ corresponds to the minimum and not the maximum
in $HP_{r,d}$.

First, we drop the positivity requirement on $\tau$ and consider the
conditions for the extremum of $HP_{r,d}$, treating $\delta$ as a
fixed parameter, which may be written as
\[
f(\xi) = \left(\begin{array}{c}
\proj_{\mathbf{1}^\perp} \left(\nabla_{\xi}E_{m+1}([\xi;x_{\mathcal{B}}]./[\delta;d_{\mathcal{B}}])\right)\\
E_{m+1}([\xi;x_{\mathcal{B}}]./[\delta;d_{\mathcal{B}}])
\end{array}\right) = \mathbf{0}\in\Re^{n-m}
\]
where $\proj_{\mathbf{1}^\perp}$ is a projection onto the subspace
orthogonal to $\mathbf{1}$ 
in some suitable basis, e.g., $\proj_{\mathbf{1}^\perp} (\omega) =
P^T\omega,\:P=[\mathbf{1}^T;-I]\in\Re^{(n-m)\times (n-m-1)}$. That
is, for a fixed $\delta$, $\xi(\delta)$ corresponds to a root of
$f(\xi)$ and the above produces $n-m$ polynomial equations in $n-m$
variables. Precisely for this reason we do not hope to obtain a
closed-form algebraic expression for $\xi(\delta)$, as it is well
known that even a single-variate polynomial of degree five and
higher in general is not solvable in radicals. Instead, we attempt
to approximate $\xi(\delta)$. The two main tools that we rely on are
the Implicit Function Theorem and Newton's method.

For fixed $\delta\in\Re$ corresponding to strictly $LP$-feasible
$d$, the Jacobian of $f(\xi)$ at $\mathbf{0}$,
\begin{eqnarray*}
J_f(\mathbf{0}) = \left.\left(\begin{array}{l}
P^T \nabla^2_\xi E_{m+1}([\xi;x_{\mathcal{B}}]./[\delta;d_{\mathcal{B}}])\\
\nabla_\xi E_{m+1}([\xi;x_{\mathcal{B}}]./[\delta;d_{\mathcal{B}}])^T
\end{array}\right)\right|_{\xi=\mathbf{0}}
\end{eqnarray*}
may be inverted by solving
\begin{equation}\label{NS:eqJinv}
J_f(\mathbf{0})\cdot\Delta_\xi=-f
\end{equation}
for an arbitrary vector
$f\in\Re^{n-m}$. Observe that just as the first equation in
$f(\xi)=\mathbf{0}$,
\[
\proj_{\mathbf{1}^\perp}
\left(\nabla_{\xi}E_{m+1}([\xi;x_{\mathcal{B}}]./[\delta;d_{\mathcal{B}}])\right)
=\mathbf{0},
\]
is satisfied if and only if there is $\tau$ such that
\[
\nabla_{\xi}E_{m+1}([\xi;x_{\mathcal{B}}]./[\delta;d_{\mathcal{B}}])-\tau\mathbf{1}
=\mathbf{0},
\]
same holds true for the linearization of this equation with respect to $\xi$.
That is, while solving $J_f(\mathbf{0})\cdot\Delta_\xi=-f$ for $\Delta_\xi$,
we may equivalently consider
\begin{eqnarray*}
\left\{\begin{array}{l}
\nabla^2_\xi E_{m+1}([\xi;x_{\mathcal{B}}]./[\delta;d_{\mathcal{B}}])
\cdot\Delta_\xi - \tau\mathbf{1}=-P (P^TP)^{-1}f_{-(n-m)},\\
\nabla_\xi E_{m+1}([\xi;x_{\mathcal{B}}]./[\delta;d_{\mathcal{B}}])^T
\cdot \Delta_\xi = -f_{n-m},
\end{array}\right.
\end{eqnarray*}
where the gradient and Hessian are evaluated at $\xi=\mathbf{0}$, which, in turn, becomes
\begin{eqnarray*}
\left\{\begin{array}{l}
E_{m-1}(x^*_{\mathcal{B}}./d_{\mathcal{B}})\cdot
\left(\mathbf{1}./\delta\:(\zeta./\delta)^T + \zeta./\delta\:
(\mathbf{1}./\delta)^T\right)\cdot\Delta_\xi +\\
E_{m-1}(x^*_{\mathcal{B}}./d_{\mathcal{B}})\cdot
\Diag(\mathbf{1}./\delta)(\mathbf{1}\mathbf{1}^T-I)\Diag(\mathbf{1}./\delta)
\cdot\Delta_\xi - \tau\mathbf{1}=-P (P^TP)^{-1}f_{-(n-m)},\\
E_m(x^*_{\mathcal{B}}./d_{\mathcal{B}})\: (\mathbf{1}./\delta)^T
\cdot \Delta_\xi = -f_{n-m},
\end{array}\right.
\end{eqnarray*}
with
\begin{eqnarray}\label{NS:zeta}
\zeta=\frac{1}{E_{m-1}(x^*_{\mathcal{B}}./d_{\mathcal{B}})}
\Diag(\delta)\left(\widetilde{A}^T\Diag(1./d_{\mathcal{B}})\nabla_z
E_m(z)|_{z=x^*_{\mathcal{B}}./d_{\mathcal{B}}}\right)
\end{eqnarray}
and $E_{m-1}(x^*_{\mathcal{B}}./d_{\mathcal{B}})\neq 0$. Noting that
the second rank-1 term, $\zeta./\delta\:(\mathbf{1}./\delta)^T$, of
the Hessian in the equation for $\Delta_\xi$ above may be replaced by
$\frac{-f_{n-m}}{E_m(x^*_{\mathcal{B}}./d_{\mathcal{B}})}\:\zeta./\delta$
due to the second equation in the above, pre-multiplying the first equation by
$\Diag(\delta)$, we get
\begin{eqnarray}\label{NS:eq}
\left\{\begin{array}{l}
E_{m-1}(x^*_{\mathcal{B}}./d_{\mathcal{B}})\:
\left(\mathbf{1}\zeta^T + \mathbf{1}\mathbf{1}^T-I\right)
\cdot\widetilde{\Delta_\xi} = \tau\delta+\widetilde{f}
,\\
E_m(x^*_{\mathcal{B}}./d_{\mathcal{B}})\: \mathbf{1}^T \cdot
\widetilde{\Delta_\xi} = -f_{n-m},
\end{array}\right.
\end{eqnarray}
where
\begin{eqnarray}\label{NS:DeltaScale}
\widetilde{\Delta_\xi}=\Diag(\mathbf{1}./\delta)\Delta_\xi,
\end{eqnarray}
and
\begin{eqnarray}\label{NS:fTilde}
\widetilde{f}=-\Diag(\delta)P(P^TP)^{-1}f_{-(n-m)}+
\frac{f_{n-m}E_{m-1}(x^*_{\mathcal{B}}./d_{\mathcal{B}})}
{E_m(x^*_{\mathcal{B}}./d_{\mathcal{B}})}\:\zeta.
\end{eqnarray}
Pre-multiplying the first equation by $\mathbf{1}^T$ and using the second equation, we get
\begin{eqnarray*}
E_{m-1}(x^*_{\mathcal{B}}./d_{\mathcal{B}})\:
\left((n-m)\zeta^T\widetilde{\Delta_\xi} - \frac{(n-m)f_{n-m}}
{E_{m}(x^*_{\mathcal{B}}./d_{\mathcal{B}})} +
 \frac{f_{n-m}}
{E_{m}(x^*_{\mathcal{B}}./d_{\mathcal{B}})}\right) =
\tau\:\mathbf{1}^T\delta+\mathbf{1}^T\widetilde{f}
\end{eqnarray*}
and so
\begin{eqnarray*}
\tau=\frac{1}{\mathbf{1}^T\delta}
\left(E_{m-1}(x^*_{\mathcal{B}}./d_{\mathcal{B}})
\left((n-m)\zeta^T\widetilde{\Delta_\xi} - \frac{(n-m)f_{n-m}}
{E_{m}(x^*_{\mathcal{B}}./d_{\mathcal{B}})} +
 \frac{f_{n-m}}
{E_{m}(x^*_{\mathcal{B}}./d_{\mathcal{B}})}\right) -
\mathbf{1}^T\widetilde{f}\right).
\end{eqnarray*}
Substituting the expression for $\tau$ back into the first equation
of~\ref{NS:eq} we have
\begin{eqnarray*}
D \cdot\widetilde{\Delta_\xi} = \frac{\delta}{\mathbf{1}^T\delta}\:
\left(-\frac{(n-m)f_{n-m}}
{E_{m}(x^*_{\mathcal{B}}./d_{\mathcal{B}})} + \frac{f_{n-m}}
{E_{m}(x^*_{\mathcal{B}}./d_{\mathcal{B}})} -
\frac{\mathbf{1}^T\widetilde{f}}{{E_{m-1}(x^*_{\mathcal{B}}./d_{\mathcal{B}})}}\right)
+\frac{\widetilde{f}}{E_{m-1}(x^*_{\mathcal{B}}./d_{\mathcal{B}})}
\end{eqnarray*}
where
\[
D=\left(\mathbf{1}-\frac{(n-m)}{\mathbf{1}^T\delta}\delta\right)\zeta^T
+ \left(\mathbf{1}\mathbf{1}^T-I\right).
\]
In turn, the above may be resolved relying on the Sherman–-Morrison
formula for the rank-1 update for the inverse of $D$, noting that
$\left(\mathbf{1}\mathbf{1}^T-I\right)^{-1}=-I-\frac{\mathbf{1}\mathbf{1}^T}{1-(n-m)}$:
\begin{eqnarray}\label{NS:Dinv}
D^{-1}=-I-\frac{\mathbf{1}\mathbf{1}^T}{1-(n-m)} +
\frac{\left(-\mathbf{1}^T\delta\:\mathbf{1}+(n-m)\:\delta\right)
\cdot
\left(\zeta+\frac{\mathbf{1}^T\zeta}{1-(n-m)}\:\mathbf{1}\right)^T}
{\rho}
\end{eqnarray}
provided
\begin{eqnarray}\label{NS:cond}
\rho=\mathbf{1}^T\delta-\mathbf{1}^T\delta\:\mathbf{1}^T\zeta+(n-m)\zeta^T\delta\neq0,
\end{eqnarray}
and so, assuming~\ref{NS:cond}, we may compute
$\widetilde{\Delta_\xi}$ and, consequently, $\Delta_\xi$; finally,
this allows us to recover the inverse of $J_f(\mathbf{0})$ from the
solution of~\ref{NS:eqJinv}.
Note that if we write $\delta=\delta_\parallel+\delta_\perp,\:
\delta_\parallel=\Delta\mathbf{1}, \mathbf{1}^T \delta_\perp=0$, the
expression for $\rho$ becomes
\begin{eqnarray}\label{NS:rho}
\rho=(n-m)\:\left(\Delta+\zeta^T\delta_\perp\right)=(n-m)\left(1+\zeta^T\frac{\delta_\perp}{\Delta}\right)\:\Delta,
\end{eqnarray}
and thus~\ref{NS:cond} may be easily satisfied by choosing $\delta$
with $\|\delta_\perp\|/\Delta$ small enough, that is, by choosing
$\delta$ close enough to $\mathcal{L}$.

Recall that according to our assumptions
$\mathcal{L}=\{\delta\in\Re^{n-m}:\delta=\Delta\mathbf{1},\Delta\in(0,\Delta_{\max})\}$.
One may formulate the following simple technical proposition.
\begin{prop}\label{NS:IFTprop}
For any fixed $(\Delta_0,\Delta_1)\subset(0,\Delta_{\max})$ there
exists $\epsilon>0$ such that for any
$\delta$ 
in the \emph{truncated central wedge}
\[
\mathcal{T}=\{\delta\in\Re^{n-m}:\delta=\delta_\parallel+\delta_\perp,
\delta_\parallel=\Delta\mathbf{1}, \mathbf{1}^T\delta_\perp=0,
\|\delta_\perp\|<\|\delta_\parallel\|\cdot\epsilon, \mbox{ and
}\Delta\in(\Delta_0,\Delta_1) \}
\]
$\xi(\delta)$ is smooth and we have
\[
\xi(\delta)=J_\xi(\delta_\parallel)\cdot\delta_\perp+O\left(\|\delta_\perp\|^2\:\mathbf{1}\right).
\]
\end{prop}
\begin{proof}
For a moment, consider $f(\xi)$ as a function $f(\delta;\xi)$ of two
vector variables $\delta$ and $\xi$; note that $f$ is $C^\infty$
with respect to both $\delta$ and $\xi$ on the interior of $LP$
feasible region. Since $\xi(\delta_\parallel)=\mathbf{0}$ and
$J_f(\mathbf{0})$ is non-singular on $\mathcal{L}$, by the Implicit
Function Theorem for any $\Delta\in[\Delta_0,\Delta_1]$ there is a
smooth function $\xi_\Delta(\delta)$ such that
$f(\delta;\xi_\Delta(\delta))=0$ for $\delta$ in some open ball
$B_\varepsilon(\delta_\parallel)$ of radius $\varepsilon>0$ centered
around $\delta_\parallel=\Delta\mathbf{1}$; clearly, the union of
such balls over all $\Delta\in[\Delta_0,\Delta_1]$ forms an open
cover of $[\Delta_0,\Delta_1]$. By compactness of
$[\Delta_0,\Delta_1]$, we may choose an open finite sub-cover
$\mathcal{U}=\bigcup_{i=2,k}B_{\varepsilon_i}(\delta_{\parallel,i})\supset[\Delta_0,\Delta_1]$
and since $B_{\varepsilon_i}(\delta_{\parallel,i})$ overlap with one
another, we may construct a smooth function
$\xi_{\:\mathcal{U}}(\delta)$ for $\delta$ in an open neighborhood
$\mathcal{U}$ of $[\Delta_0,\Delta_1]$. In particular,
$\xi_{\:\mathcal{U}}(\delta)$ is twice continuously differentiable
and $\epsilon>0$ may be chosen small enough so that the truncated
wedge $\mathcal{T}\subset\mathcal{U}$. The result follows from
Taylor's expansion of $\xi_{\:\mathcal{U}}(\delta)$.

Finally, observe that $\xi=\xi_{\:\mathcal{U}}(\delta)$ indeed
corresponds to the minimizer of $HP_{r,d}$ for some fixed $d=[\delta;d_\mathbf{B}]$,
and is not an arbitrary root of $f$. To show that
$x=[\xi;x_\mathcal{B}]\in\mathcal{K}_{n-m-1,d}$ note that
Proposition~\ref{Hyper:MultRoot} implies that the branches of $E_{m+1}(x./d)=0$
are distinct except for at the faces of $\Re^n_+$, and thus, if $\xi_{\:\mathcal{U}}(\delta)$
switches branches, then it cannot be smooth or even continuous; therefore, $\xi_{\:\mathcal{U}}(\delta)$ must result
in $x(d)$ that additionally satisfies the conditions of Corollary~\ref{Hyper:ConeBdry}.
Similarly, since the feasible region of $LP$ is assumed to have a non-empty
interior, and so, in particular contains an open ball of some radius $\varepsilon>0$,
$\xi_{\:\mathcal{U}}(\delta)$ cannot switch from being a minimizer
at $\xi(\mathbf{0})=\mathbf{0}$ to being a maximizer at some other $\delta\in\mathcal{T}$
and yet stay smooth, because $HP_{r,d}$ is a relaxation of $LP$ and consequently
$c^T x(d)$ must be smaller than the maximum of $HP_{r,d}$ by at least $\|c\|\varepsilon>0$.
\end{proof}
\begin{rem}
Strictly speaking, when discussing the Jacobian of $\xi(\delta)$ in
previous subsection, we should have justified the existence of
differentiable $\xi(\delta)$ first as in the above proposition. We
intentionally delayed this discussion till the present subsection in
an attempt to keep our motivation more transparent.
\end{rem}
\noindent The above and Proposition~\ref{xi:Jacobian} result in the
following straightforward consequence.
\begin{cor}
For any $(\Delta_0,\Delta_1)\subset(0,\Delta_{\max})$, the truncated
central wedge $\mathcal{T}$ may be chosen so that $\xi(\delta)$ is
continuously differentiable for $\delta\in\mathcal{T}$ and
\[
\xi(\delta)=-\theta(\delta)\cdot\delta_\perp+O\left(\frac{\|\delta_\perp\|^2}{\|\delta_\parallel\|}\:\mathbf{1}\right).
\]
for some $\theta(\delta)\geq\theta>0$.
\end{cor}
\noindent Note that the big-$O$ constant in the above might depend
on the choice of $\Delta_0,\Delta_1$.

If we could show that $\xi(\delta)$ is continuously differentiable
in some neighborhood $\mathcal{U}$ of $\delta=\mathbf{0}$, by
combining $\mathcal{U}$ with $\mathcal{T}$ the last corollary would
imply that for any $0<\Delta_1<\Delta_{\max}$ we may chose the
central wedge
$\widetilde{\mathcal{W}}\subset\:\mathcal{U}\cup\mathcal{T}$ as in
Corollary~\ref{ODE:mainCor}, where $\xi(\delta)$ is continuously
differentiable, and so is locally Lipschitz continuous. Moreover, if
$\xi(\delta)$ was well-behaved on $\mathcal{U}$ in a sense
of~\ref{Dynamics:ODEtarget}, this would imply our main result.

That is, currently, not only we cannot guarantee that the quality of
big-$O$ approximation in the above corollary for
$\delta\in\mathcal{T}$ does not deteriorate too fast as $\Delta_0$
gets closer and closer to $0$, we are not even guaranteed that
$\xi(\delta)$ is smooth enough in any central wedge
$\widetilde{\mathcal{W}}$ to guarantee the existence of the solution
to ODE~\ref{Dynamics:ODEnm} near $\delta=\mathbf{0}$.

In particular, recall that the existence of the inverse of
$J_f(\mathbf{0})$ depends on~\ref{NS:cond}, and so, potentially may
be compromised in the limit as $\delta\rightarrow\mathbf{0}$,
preventing us from being able to extend $\mathcal{T}$ in the above
proposition and corollary to enclose $\delta=\mathbf{0}$. Indeed, as
the following example illustrates, $\xi(\delta)$ may fail to be
differentiable at $\mathbf{0}$.

Example:
\begin{itemize}
\item considering the same problem as at the beginning of Section~\ref{SecSWIntro}
with its relaxation, $\min_x\{(1, 1, 0)^T x:\: \mathbf{1}^Tx=3, x\in\Re^3_+\}$ and
$\min_x\{(1, 1, 0)^T x:\: \mathbf{1}^Tx=3, x\in\mathcal{K}_{1,d}\}$, we derive the explicit
expression for $\xi(\delta)$.

$LP$ optimum is $x^*=(0, 0, 3)$, so $\mathcal{B}=\{3\}$
and $\delta=(d_1, d_2), \xi=(x_1, x_2)$. Recall that the boundary satisfies
$\frac{1}{2}\xi^T Q \xi+r^T\xi+s=0$ with $Q, r, s$ as before, namely,
\begin{eqnarray*}
\begin{array}{cc}
Q=\frac{6}{\delta_1\delta_2(3-\delta_1-\delta_2)}
\left(
\begin{array}{cc}
-2\delta_2 & 3-2(\delta_1+\delta_2)\\
3-2(\delta_1+\delta_2) & -2\delta_1
\end{array}
\right),&
r=\frac{6}{\delta_1\delta_2(3-\delta_1-\delta_2)}
\left(
\begin{array}{cc}
\delta_2\\
\delta_1
\end{array}
\right).
\end{array}
\end{eqnarray*}
The optimality
conditions for $\xi(\delta)$ correspond to
\[
\nabla_\xi \left(\frac{1}{2}\xi^T Q \xi+r^T\xi+s\right)=Q \xi+r=\tau\mathbf{1},\:\tau>0,
\]
and recalling $\det(\widetilde{Q})=-9+12(\delta_1+\delta_2)-4(\delta_1^2+\delta_2^2)-4\delta_1\delta_2$,
for small enough $\delta\in\Re^2_{++}$ may be
equivalently re-written with $\widetilde{\tau}>0$ as
\[
\xi =
\frac{1}{\det(\widetilde{Q})}
\left(
\begin{array}{cc}
-2\delta_1 & 2(\delta_1+\delta_2)-3\\
2(\delta_1+\delta_2)-3 & -2\delta_2
\end{array}
\right)\cdot
\left(\widetilde{\tau}\mathbf{1}-
\left(
\begin{array}{cc}
\delta_2\\
\delta_1
\end{array}
\right)
\right).
\]
Substituting $\xi$ back into the boundary condition to get $\widetilde\tau$ we get
\[
\widetilde{\tau}^2\:\mathbf{1}^T\widetilde{Q}^{-1}\mathbf{1}-
\widetilde{r}^T \widetilde{Q}^{-1}\widetilde{r} = 0,
\]
with $\widetilde{r}=\frac{\delta_1\delta_2(3-\delta_1-\delta_2)}{6}\:r$, and so
\[
\widetilde\tau=\sqrt{\delta_1\delta_2}
\]
as out of the two quadratic roots we are interested in positive $\widetilde\tau$.
Finally, observe that $\widetilde{\tau}$ results in $\xi(\delta)$ not being
differentiable at $\delta=\mathbf{0}$. The figure below illustrates a position of
$\xi(\delta)$ in relationship to its Jacobian-based approximation for one particular $\delta$.
\begin{figure}[h!]
\begin{center}
\epsfig{file=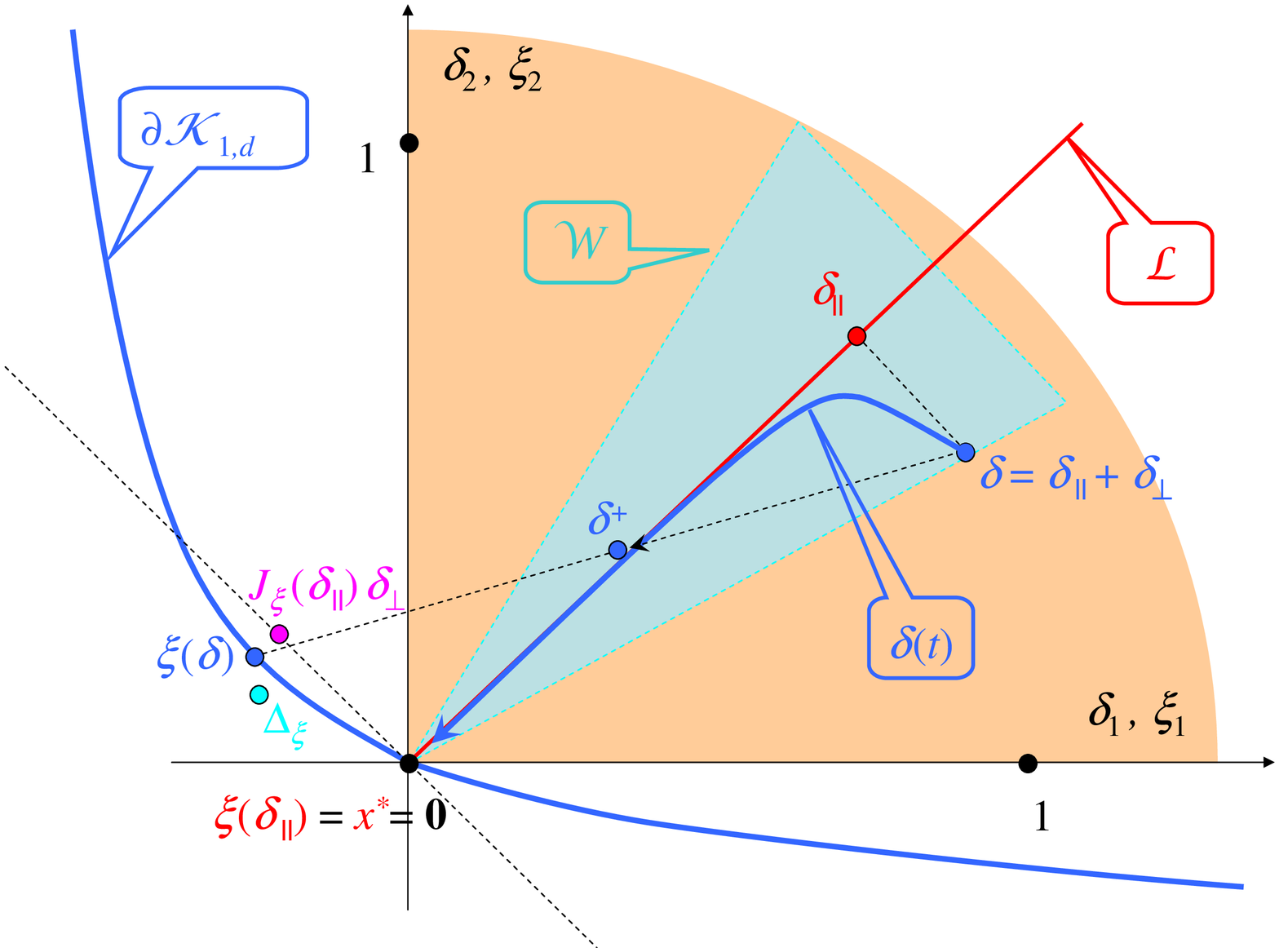,height=10.0cm} \caption{Shrink-Wrapping dynamics close-up}\label{Hyper:Fig5}
\end{center}
\end{figure}
\end{itemize}

\begin{rem} Non-differentiability of $\xi(\delta)$ at $\delta=\mathbf{0}$
also prevents us from relying on a standard ODE \emph{sink}-type
argument~\cite{Arnold:ODE}, as clearly, in the way it is defined,
$\xi(\delta)$ does not even exist beyond the nonnegative orthant. It
is conceivable that from purely algebraic point of view one may
extend $\xi(\delta)$ beyond $\Re^{n-m}_{++}$ as, say, a solution to
the polynomial system of equations. However, the basic problem of
non-differentiability at $\delta=\mathbf{0}$ is still likely to
persist if we continue using Euclidian coordinates for $\xi,\delta$.
Along the latter lines, in~\cite{Zinchenko:thesis} it has been
suggested that perhaps a non-liner change of coordinates, namely,
spherical coordinates, might be a more suitable choice to address
the problem; in particular, such a choice allows to overcome this
difficulty in the example above and subsequently permits the usage
of a sink. To justify the existence of a continuously differentiable
$\xi(\delta)$ in spherical coordinates beyond nonnegative orthant
one may attempt to use the Implicit Function Theorem. However, here
we choose to follow Newton's method-based analysis as, hopefully, it
may subsequently be used to lay down the ground work for the
path-following in the actual optimization algorithm.
\end{rem}

To remedy the situation, we rely on a different approximation to
$\xi(\delta)$ for $\delta$ near $\mathbf{0}$, namely, we use the
first iterate of Newton's method and its error analysis as
in~\cite{ShubSmale:Complexity}.

Fix $\widetilde{\delta}\in\mathcal{L}$ so that the corresponding
$\xi(\widetilde{\delta})=\mathbf{0}$ and consider
$\delta\in\Re^{n-m}_{++},\: \delta\neq\widetilde{\delta}$. If
$\delta$ is chosen close to $\widetilde{\delta}$, we may attempt to
approximate $\xi(\delta)$ by finding an approximate root
$\xi(\widetilde{\delta})+\Delta_\xi=\Delta_\xi$ of $f(\xi)$ based on
linearization at $\mathbf{0}$. That is, we solve the
following equation for the Newton step $\Delta_\xi$:
\[
f(\xi)\approx f(\mathbf{0}) + J_f(\mathbf{0})\cdot
\Delta_\xi=\mathbf{0},
\]
where $J_f(\mathbf{0})$ is the Jacobian of $f(\xi)$ at
$\xi=\mathbf{0}$.

Intuitively, in the limit as $\delta\rightarrow\widetilde{\delta}$
the Newton step $\Delta_{\xi}$ must resemble the first-order
approximation to $\xi(\delta)$ obtained with the Jacobian
$J_{\xi}(\widetilde{\delta})$: both rely on linearizations of
$\xi(\delta)$ but at ever so slightly different points $\delta$ and
$\widetilde{\delta}$. The latter provides a key motivation for
working with Newton iterates, as we hope to obtain a usable
approximation to $\xi(\delta)$ that takes on a nearly-projection
form similar to $J_{\xi}(\widetilde{\delta})$ acting on
$\delta-\widetilde{\delta}$.

To compute the Newton step $\Delta_\xi$ we specialize $f$ to
$f(\mathbf{0})$ in~\ref{NS:eqJinv}, that is,
\begin{equation*}
f = \left(\begin{array}{c}f_{-(n-m)}\\f_{n-m}\end{array}\right) =
\left.\left(\begin{array}{c}
P^T \nabla_{\xi}E_{m+1}([\xi;x_{\mathcal{B}}]./[\delta;d_{\mathcal{B}}])\\
0 \end{array}\right)\right|_{\xi=\mathbf{0}} =
\left(\begin{array}{c}
P^T E_{m}(x^*_{\mathcal{B}}./d_{\mathcal{B}})\:\mathbf{1}./\delta\\
0 \end{array}\right).
\end{equation*}
With the above in mind, we have
\begin{eqnarray*}
\begin{array}{rll}
\widetilde{f} = & -\Diag(\delta)P(P^TP)^{-1}P^T E_{m}(x^*_{\mathcal{B}}./d_{\mathcal{B}})\:\mathbf{1}./\delta & = \\
& -E_{m}(x^*_{\mathcal{B}}./d_{\mathcal{B}})\Diag(\delta)
\left(I-\frac{\mathbf{1}\mathbf{1}^T}{n-m}\right)\:\mathbf{1}./\delta & = \\
& -E_{m}(x^*_{\mathcal{B}}./d_{\mathcal{B}})
\left(\mathbf{1}-\frac{\mathbf{1}^T\mathbf{1}./\delta}{n-m}\:\delta\right),
\end{array}
\end{eqnarray*}
and
\[
\frac{\delta}{\mathbf{1}^T\delta}\: \left(-
\frac{\mathbf{1}^T\widetilde{f}}{E_{m-1}(x^*_{\mathcal{B}}./d_{\mathcal{B}})}\right)
+\frac{\widetilde{f}}{E_{m-1}(x^*_{\mathcal{B}}./d_{\mathcal{B}})} =
\frac{E_m(x^*_{\mathcal{B}}./d_{\mathcal{B}}}{E_{m-1}(x^*_{\mathcal{B}}./d_{\mathcal{B}}}
\left(-\mathbf{1}+\frac{n-m}{\mathbf{1}^T\delta}\:\delta\right).
\]
Thus, using the earlier expression for $D^{-1}$, distributing all
the terms and simplifying, we get the following expression for the
scaled Newton step $\widetilde{\Delta_\xi}$:
\begin{eqnarray*}
\begin{array}{rl}
\widetilde{\Delta_\xi} & =
-\frac{E_m(x^*_{\mathcal{B}}./d_{\mathcal{B}})}{E_{m-1}(x^*_{\mathcal{B}}./d_{\mathcal{B}})}
\left(-\mathbf{1}+\frac{n-m}{\mathbf{1}^T\delta}\:\delta
+\frac{(n-m)\zeta^T\delta-\mathbf{1}^T\delta\:\mathbf{1}^T\zeta}{\gamma}
\left(\mathbf{1}-\frac{n-m}{\mathbf{1}^T\delta}\:\delta\right)
\right)\\
& =
-\frac{E_m(x^*_{\mathcal{B}}./d_{\mathcal{B}})}{E_{m-1}(x^*_{\mathcal{B}}./d_{\mathcal{B}})}
\:\frac{1}{\gamma}\left((n-m)\:\delta-\mathbf{1}^T\delta\:\mathbf{1}\right).
\end{array}
\end{eqnarray*}

Using $\delta=\delta_\parallel+\delta_\perp$, $\delta_\parallel=\Delta\mathbf{1}, \mathbf{1}^T\delta_\perp=0$, and observing
\[
\mathbf{1}^T\delta-\mathbf{1}^T\delta\:\mathbf{1}^T\zeta+(n-m)\zeta^T\delta=(n-m)\:(\Delta+\zeta^T\delta_\perp),
\]
we get
\[
\widetilde{\Delta_\xi}=-\frac{E_m(x^*_{\mathcal{B}}./d_{\mathcal{B}})}{E_{m-1}(x^*_{\mathcal{B}}./d_{\mathcal{B}})}
\:\frac{\delta_\perp}{\Delta+\zeta^T\delta_\perp},
\]
and so, re-scaling by $\Diag(\delta)$ we finally have
\begin{eqnarray}\label{NS:delta}
\Delta_\xi=-\frac{E_m(x^*_{\mathcal{B}}./d_{\mathcal{B}})}{E_{m-1}(x^*_{\mathcal{B}}./d_{\mathcal{B}})}
\:\left(\frac{\Delta}{\Delta+\zeta^T\delta_\perp}\:\delta_\perp+\frac{1}{\Delta+\zeta^T\delta_\perp}\:(\delta_\perp).^2\right).
\end{eqnarray}

For real analytic $f$, it is well known that under mild
non-degeneracy assumptions, namely the invertibility of the Jacobian
of $f$ at the root $\xi$, Newton's method converges quadratically to
the associated root $\xi$ if the initial iterate $z$ is chosen close
enough to $\xi$. To this extent we formulate a slightly more
specialized and simple result following the analysis
in~\cite{ShubSmale:Complexity}, introducing two auxiliary quantities
\[
\beta_z=\|f(z)^{-1} f(z)\|,
\]
which corresponds to the length of the Newton step at $z$, and
\[
\gamma_z=\sup_{k\geq2}\left\|\frac{f'(z)^{-1}
f^{(k)}(z)}{k!}\right\|^{\frac{1}{k-1}}.
\]
\begin{lem}\label{Newton:Lemma}
There is a universal constant $\alpha_1>0$ such that if
\[
\alpha_z=\beta_z\cdot\gamma_z<\alpha_1
\]
then the distance from $z$ to the associated zero $\xi$ decreases
quadratically with each Newton iteration starting from $z$, that is, denoting $z_{(0)}=z$ and
\[
z_{(i+1)}=z_{(i)}-f'\left(z_{(i)}\right)^{-1} f\left(z_{(i)}\right),
i>0,
\]
for all $i\geq 0$ we have
\[
\|z_{(i+1)}-\xi\|\leq
\frac{4\gamma_z}{(5-\sqrt{17})\Psi(2\alpha_z)(1-2\alpha_z)}
\cdot\|z_{(i)}-\xi\|^2,
\]
where $\Psi(u)=2u^2-4u+1$, and, moreover, $\|z-\xi\|\leq2\beta_z$.
\end{lem}
\begin{proof}
Pick
$\alpha_1>0$ 
so that $2\alpha_1$ is less then the real root of
$2u^3-6u^2+\left(5+\frac{4}{5-\sqrt{17}}\right)u-1$; since $\Psi(u)$
is monotone decreasing for $u<1-\frac{\sqrt{2}}{2}$, for
$u\in[0,2\alpha_1)$ we have
$\frac{u}{\Psi(u)(1-u)}<\frac{5-\sqrt{17}}{4}$. By the cubic root
formula it may be verified that the decimal expansion of the real
root of the above polynomial, truncated to first five significant
digits, is $.11218$.

Note that $\alpha_1<\alpha_0=\frac{1}{4}(13-3\sqrt{17})\approx
.15767$, where $\alpha_0$ is the best known value for the constant
in Theorem 2 in Section 8 of~\cite{ShubSmale:Complexity}; therefore,
the theorem implies
\[
\|z-\xi\|\leq2\beta_z=\frac{2\alpha_z}{\gamma_z}<\frac{2\alpha_1}{\gamma_z}.
\]

Since $\|z-\xi\|\cdot\gamma_z<2\alpha_1<1-\frac{\sqrt{2}}{2}$,
Proposition 3 in Section 8 of~\cite{ShubSmale:Complexity} implies
\[
\gamma_\xi\leq\frac{\gamma_z}{\Psi(2\alpha_z)(1-2\alpha_z)}.
\]

Since $\|z-\xi\|\cdot\gamma_\xi\leq \frac{2\alpha_z}{\gamma_z}\cdot
\frac{\gamma_z}{\Psi(2\alpha_z)(1-2\alpha_z)} <
\frac{5-\sqrt{17}}{4}$, by Proposition 1 in Section 8
of~\cite{ShubSmale:Complexity}
\[
\|z_{(1)}-\xi\|\leq
\frac{\gamma_\xi}{\Psi(\|z-\xi\|\gamma_\xi)}\cdot\|z-\xi\|^2
<\frac{4\gamma_z}{(5-\sqrt{17})\Psi(2\alpha_z)(1-2\alpha_z)}\cdot\|z-\xi\|^2
\]
follows. Observing that the above inequality, in particular, implies
$\|z_{(1)}-\xi\|<\|z-\xi\|$, the last proposition may be re-applied to estimate
$\|z_{(2)}-\xi\|$ since $\|z_{(1)}-\xi\|\cdot\gamma_\xi<\frac{5-\sqrt{17}}{4}$
and so on, thus, completing the statement of our lemma.
\end{proof}
\noindent In particular, we may choose $\alpha_1=.11218$, in which
case
\begin{eqnarray}\label{Newton:quadBound}
\left\|\left(z-f'(z)^{-1}f(z)\right)-\xi\right\|< 20\gamma_z
\cdot\|z-\xi\|^2<80\gamma_z\cdot\beta_z^2.
\end{eqnarray}

\begin{prop}
For any fixed $0<\Delta_1<\Delta_{\max}$ there exists
$\widetilde\epsilon>0$ such that for any
$\delta$ 
in the central wedge
\[
\widetilde{\mathcal{W}}=\{\delta\in\Re^{n-m}:\delta=\delta_\parallel+\delta_\perp,
\delta_\parallel=\Delta\mathbf{1}, \mathbf{1}^T\delta_\perp=0,
\|\delta_\perp\|<\|\delta_\parallel\|\cdot\widetilde\epsilon, \mbox{
and }\Delta\in(0,\Delta_1) \}
\]
$\xi(\delta)$ is smooth and we have
\[
\xi(\delta)=\Delta_\xi+O\left(\frac{\|\Delta_\xi\|^2}{\|\delta_\parallel\|}\:\mathbf{1}\right).
\]
\end{prop}
\begin{proof}
We rely on the result of the previous lemma,
namely,~\ref{Newton:quadBound}. Note that $\zeta$, as defined
by~\ref{NS:zeta}, remains bounded from above on any $LP$ strictly
feasible closure of $\widetilde{\mathcal{W}}$. Consequently, the
existence of the Newton step $\Delta_\xi$ guaranteed
by~\ref{NS:cond}, recalling~\ref{NS:rho}, may be ensured for any
$\delta\in\widetilde{\mathcal{W}}$ by choosing
$\widetilde{\epsilon}>0$ sufficiently small. So, for a moment, fix
$\widetilde{\epsilon}$ such that $\|\zeta\|<K,\:K>0$, for all
$\delta\in\widetilde{\mathcal{W}}$, and $1-K\widetilde\epsilon>1/2$.
Then for
$\delta=\delta_\parallel+\delta_\perp\in\widetilde{\mathcal{W}}$
by~\ref{NS:delta} we have
$\beta_{\mathbf{0}}=\|\Delta_\xi\|=\widetilde{\epsilon}\:O(\|\delta_\parallel\|)$.
If necessary, we will refine our choice of $\widetilde{\mathcal{W}}$
at a later point by further reducing $\widetilde{\epsilon}$.

It is left to analyze $\gamma_{\mathbf{0}}$; recall that for an
operator $F$ its (induced) norm is defined as
$\sup_{\chi\neq\mathbf{0}} \frac{\|F(\chi)\|}{\|\chi\|}$. Since $f$
is polynomial, for $k>m+1$ the differential $f^{(k)}$ vanishes. For
$k\leq m+1$, the $k$-order differential of $f$, evaluated at a fixed
$k$-tuple $\chi$ is a vector whose first $n-m-1$ components are of
order $1/\|\delta\|^{k+1}$ with respect to $\delta$, and the last
component is of order $1/\|\delta\|^{k}$. Recalling that
$f'(\mathbf{0})^{-1}=J_f(\mathbf{0})^{-1}$ may be recovered from the
solution to~\ref{NS:eqJinv}, observing from~\ref{NS:Dinv} that
$\|D^{-1}\|=O(1)$ on $\widetilde{\mathcal{W}}$, and recalling the
definition for $\widetilde{f}$ as in~\ref{NS:fTilde}, particularly,
that the first $n-m$ components of $\widetilde{f}$ are scaled by
$\Diag(\delta)$, and the fact that from~\ref{NS:DeltaScale} we have
$\Delta_\xi=\Diag(\delta)\:\widetilde{\Delta_\xi}$, we conclude that
the composite $k$-linear operator $f'(z)^{-1} f^{(k)}(z)$ acting on
a fixed $k$-tuple $\chi$ results in a vector of order
$1/\|\delta\|^{k-1}$. Now, applying the induced norm and taking
$k-1$ root we conclude that
$\gamma_{\mathbf{0}}=O=(1/\|\delta\|)=O(1/\|\delta_\parallel\|)$ on
$\widetilde{\mathcal{W}}$.

Combining our estimates for $\beta_{\mathbf{0}},
\gamma_{\mathbf{0}}$, we get
$\alpha_{\mathbf{0}}=\beta_{\mathbf{0}}\cdot
\gamma_{\mathbf{0}}=\widetilde{\epsilon}\:O(1)$ on
$\widetilde{\mathcal{W}}$. So, if needed, $\widetilde{\epsilon}$ and
the corresponding central wedge $\widetilde{\mathcal{W}}$ may be
further reduced to result in $\alpha_{\mathbf{0}}<\alpha_1$ on
$\widetilde{\mathcal{W}}$, completing our estimate on $\xi(\delta)$.

Finally, observe that $\alpha_{\mathbf{0}}<\alpha_1$ on $\widetilde{\mathcal{W}}$ in
particular implies that $J_f(\xi(\delta))$ is invertible since $\gamma_{\xi(\delta)}$
is finite, see the above lemma. Invoking the argument similar to that of
Proposition~\ref{NS:IFTprop}, considering a finite open cover of the closure of
$\widetilde{\mathcal{W}}$, the Implicit Function Theorem implies the existence
of a smooth $\xi(\delta)$ defined on $\widetilde{\mathcal{W}}$ that corresponds to
the minimizers of $HP_{r,d}$.
\end{proof}
\noindent The above combined with~\ref{NS:delta}
result in the following straightforward consequence.
\begin{cor}
For any $0<\Delta_1<\Delta_{\max}$, the central wedge
$\widetilde{\mathcal{W}}$ may be chosen so that $\xi(\delta)$ is
continuously differentiable for $\delta\in\widetilde{\mathcal{W}}$
and
\[
\xi(\delta)=-\theta(\delta)\cdot\delta_\perp+O\left(\frac{\|\delta_\perp\|^2}{\|\delta_\parallel\|}\:\mathbf{1}\right).
\]
for some $\theta(\delta)\geq\theta>0$.
\end{cor}
\noindent The last corollary completes the proof of our main result
-- Theorem~\ref{Dynamics:MainTheorem}; the behavior of $\xi(t)$ is a
straightforward consequence of exponential convergence of
$\delta(t)$ to $\mathcal{L}$.

\begin{rem}
The actual basin of exponential convergence to $\mathcal{L}$, that
is, a subset of $LP$ feasible region starting from which the
trajectories $\delta(t)$ converge exponentially to $\mathcal{L}$,
and consequently, to the $LP$ optimum at $\delta=\mathbf{0}$, might
be far more complicated than simply a central wedge $\mathcal{W}$;
for once, such a set must necessarily contain the union of all the
central wedges as in Theorem~\ref{Dynamics:MainTheorem}, each
corresponding to different $0<\Delta_1<\Delta_{\max}$. Moreover,
instead of relying on Newton's method-based analysis of
$\xi(\delta)$, alternatively we could combine the central wedge
$\mathcal{W}$ for small $\|\delta\|$, addressing the potential
non-differentiability of $\xi(\delta)$ at $\mathbf{0}$, and the
truncated central wedge $\mathcal{T}$ for $\|\delta\|$ relatively
large. Note that intuitively, for
$\delta=\delta_\parallel+\delta_\perp$ close to $\mathcal{L}$ in a
sense of small $\|\delta_\perp\|/\|\delta_\parallel\|$, the
approximation to $\xi(\delta)$ based on the Jacobian
$J_\xi(\delta_\parallel)$ becomes
\[
\xi(\delta)\approx
J_\xi(\delta_\parallel)\:(\delta-\delta_\parallel)\approx
-\frac{E_m(x^*_{\mathcal{B}}./d_{\mathcal{B}})}{E_{m-1}(x^*_{\mathcal{B}}./d_{\mathcal{B}})}
\:\delta_\perp,
\]
while the Newton step approximation with small $\|\delta\|,
\|\delta_\perp\|/\|\delta_\parallel\|$ also results in
\[
\xi(\delta)\approx\Delta_\xi\approx
-\frac{E_m(x^*_{\mathcal{B}}./d_{\mathcal{B}})}{E_{m-1}(x^*_{\mathcal{B}}./d_{\mathcal{B}})}
\:\delta_\perp,
\]
so we expect the two approximations to act alike, see
Figure~\ref{Hyper:Fig5}. In addition to the above, if we were to
rely solely on Newton's method, it is well known that the basin of
convergence for the method alone may very well be extraordinarily
complicated -- see, for example, Mandelbrot
set~\cite{ShubSmale:Complexity}.
\end{rem}

We conjecture that $\delta(t)$ enters the central wedge $\mathcal{W}$
for some $t>0$ if $\xi(\delta(0))$ exists.

To illustrate the kind of implications continuous trajectories
$\delta(t), \xi(t),\:t\geq0$ might have for the resulting
optimization algorithm, which would most certainly operate on
discrete iterates $\delta^{(i)}, \xi^{(i)},\:i\geq0$, consider the
following proposition; as before, for simplicity we assume that
$\xi(\delta)$ is easily available given $\delta$. Once again, it is
convenient to adapt the following notation: $\delta_\parallel$
denotes the first coordinate of $\delta$ where the first coordinate
axis is aligned with $\mathbf{1}\in\Re^{(n-m)}$, $\delta_\perp$
denotes the remaining $(n-m)-1$ orthonormal coordinates of the
vector $\delta$ in this new coordinate system.
\begin{prop}
Given the initial iterate $\delta^{(0)}, \xi^{(0)}=\xi(\delta^{(0)})$, consider a simple
bisection-type scheme for determining the values of $\delta$:
\[
\delta^{(i+1)}=\delta^{(i)}+\frac{m}{m+1}\left(\xi^{(i)}-\delta^{(i)}\right),\:i>0,
\]
where $\xi^{(i)}=\xi(\delta^{(i)})$. If $\delta^{(0)}$ satisfies
$\|\delta^{(0)}_\perp\|/\delta^{(0)}_\parallel<\epsilon,
\delta^{(0)}_\parallel>0$ for sufficiently small $\epsilon>0$, i.e.,
is inside the properly chosen central wedge $\mathcal{W}$, and, in
addition, $\|\delta^{(0)}\|$ is sufficiently small, then the
iterates $\delta^{(i)}$ converge at least $R$-linearly and
$\xi^{(i)}$ converge $R$-super-quadratically to the $LP$ optimum
$\mathbf{0}$, in particular, for some $K>0$
\[
\|\delta^{(i)}\|\leq
\|\delta^{(0)}\|\left(\frac{1}{2}\right)^{i}\mbox{ and  }
\|\xi^{(i)}\|\leq K\left(\frac{1}{2}\right)^{2^{i}+i},\:i\geq0.
\]
\end{prop}
\begin{proof}
For brevity of notation we use $\delta$ to denote $\delta^{(0)}$,
and $\delta^+$ to denote $\delta^{(1)}$. According to the previous
corollary, choosing $\delta\in\widetilde{\mathcal{W}}$ we can write
\[
\xi(\delta)=\left(\begin{array}{c}
O\left(\frac{\|\delta_\perp\|^2}{\delta_\parallel}\right) \\
-\left(\frac{1}{m}+O(\|\delta\|)\right)\:\delta_\perp+
O\left(\frac{\|\delta_\perp\|^2}{\delta_\parallel}\:\mathbf{1}\right)
\end{array}\right)
\]
since the limit of differentiable function
$\frac{E_m(x^*_{\mathcal{B}}./d_{\mathcal{B}})}
{E_{m-1}(x^*_{\mathcal{B}}./d_{\mathcal{B}})}$ is $\frac{1}{m}$ as
$\delta\rightarrow\mathbf{0}$, 
and so
\[
\delta^+=\delta+\frac{m}{m+1}\left(\xi(\delta)-\delta\right)=
\frac{1}{m+1}\left(\begin{array}{c}
\delta_\parallel \\
-m O(\|\delta\|)\:\delta_\perp
\end{array}\right) +
\frac{m}{m+1}\left(\begin{array}{c}
O\left(\frac{\|\delta_\perp\|^2}{\delta_\parallel}\right) \\
O\left(\frac{\|\delta_\perp\|^2}{\delta_\parallel}\:\mathbf{1}\right)
\end{array}\right).
\]
Clearly, we may choose $\|\delta\|$ sufficiently small so that in
the expression above we have
\[
m O(\|\delta\|)<\frac{1}{3},
\]
and thus, considering the first and the last $(n-m)-1$ components of
$\delta^+$ we can write
\begin{eqnarray}\label{Dynamics:IneqDiscrete}
\begin{array}{c}
\delta^+_\parallel=\frac{\delta_\parallel}{m+1}\:\left(1+O(\epsilon^2)\right),\\
\|\delta^+_\perp\|=\frac{\delta_\parallel}{m+1}\:\left(\frac{\epsilon}{3}+O(\epsilon^2)\right).
\end{array}
\end{eqnarray}
So, if necessary, we may further reduce $\epsilon>0$ in
$\|\delta_\perp\|/\delta_\parallel<\epsilon<1$ to guarantee
\[
\delta^+_\parallel\leq\frac{1}{2}\delta_\parallel, \mbox{ and }
\frac{\|\delta^+_\perp\|}{\delta^+_\parallel}\leq\frac{\epsilon}{2}.
\]

Repeating the argument at $\delta=\delta^{(1)},
\delta^+=\delta^{(2)}$, and observing that now the quantity $m
O(\|\delta\|)<\frac{1}{3}$ also gets at least halved,
from~\ref{Dynamics:IneqDiscrete} we have
\[
\delta^{(2)}_\parallel\leq\left(\frac{1}{2}\right)^2\delta_\parallel^{(0)},
\mbox{ and }
\frac{\|\delta^{(2)}_\perp\|}{\delta^{(2)}_\parallel}\leq\frac{1}{2}\left(\frac{1}{2}\right)^2\epsilon,
\]
and so on for $i>2$, ultimately resulting in
\[
\delta^{(i)}_\parallel\leq\left(\frac{1}{2}\right)^i\delta_\parallel^{(0)},
\mbox{ and }
\frac{\|\delta^{(i)}_\perp\|}{\delta^{(2)}_\parallel}\leq\left(\frac{1}{2}\right)^{2^i-1}\epsilon.
\]
The last two bounds combined with the expression for $\xi(\delta)$
in the previous corollary give us $R$-super-quadratic bound on
$\xi^{(i)}$ as claimed; the bound on $\|\delta^{(i)}\|$ follows
trivially.
\end{proof}
\noindent Note that the above bisection scheme in fact does not
require us to know the optimal basis $\mathcal{B}$ a priori: replace
$\delta$ iterates with the corresponding
$d^{(i+1)}=d^{(i)}+\frac{m}{m+1}\left(x^{(i)}-d^{(i)}\right),\:i>0$,
where $x^{(i)}=x(d^{(i)})$. For an illustration, see
Figure~\ref{Hyper:Fig5}.

We fully anticipate the criticism of the last proposition as being
reliant on very strong assumptions from any practical point of view,
e.g., the availability of $\xi(\delta)$. However, it should be
understood that the purpose of the latter proposition is, at this
point, solely illustrative. We would like to add that in our limited
computational experiments we observed that indeed $\xi$ appears to
be a much more promising candidate to follow numerically, as the
iterates $\xi^{(i)}$ seem to converge to the optimum much sooner
than the corresponding $\delta^{(i)}$. This suggests that when
designing the actual optimization algorithm based on the
Shrink-Wrapping setting one might benefit from focusing on
$\xi^{(i)}$ rather than $\delta^{(i)}$; in the subsequent section we
will see that the iterates $\delta^{(i)}$ have an existing analogue
in the interior-point methods, while, in contrast, $\xi^{(i)}$
appear to be quite unique to Shrink-Wrapping.

\begin{rem}
A natural direction in refining the last proposition towards making
it more or less practically meaningful is to consider the Newton
based approximation to $\xi^{(i+1)}$ from the previous iterate
$\xi^{(i)}$; the latter is consistent with numerical path-following
approach commonly employed by the interior-point methods. Also note
that for large $m$, the ratio $m/(m+1)$ is very close to 1, e.g.,
when $m=99$, the Shrink-Wrapping iterates $\delta^{(i)}$ would
traverse at least 99\% of the distance to the boundary of the $LP$
feasible region with each step. This, again, is consistent with the
so-called predictor-corrector-type interior-point methods. Moreover,
in order to get fast convergence of $\xi^{(i)}$ iterates, most
probably we can get away with requiring the multiplier in front of
$\left(\xi^{(i)}-\delta^{(i)}\right)$ to approach the value
$m/(m+1)$ only asymptotically. Lastly, amongst many other immediate
potential research directions, developing an intrinsic proximity
measure of an iterate to the central line appears to be of great
importance. However, we believe that these questions go well beyond
the scope of this paper.
\end{rem}


\section{Pathological central paths vs. Shrink-Wrapping}
\label{SecSWvsIPM} In this section, we contrast the behavior of the
central path to the Shrink-Wrapping trajectories $d(t), t\geq0$, for
some known LP instances with large total curvature of the central
path. Namely, we consider the following three LP instances:
\emph{Megiddo-Shub simplex}~\cite{MS:simplex}, \emph{DTZ
snake}~\cite{DTZ:snake}, and \emph{redundant Klee-Minty
cube}~\cite{DNT:KMcube},~\cite{DTZ:KMcube}.

The total curvature of a smooth curve --here, the central path-- is
defined as a definite integral over the total length of the curve of
the norm of the curvature vector, where the latter corresponds to
the second derivative of the curve equation parameterized by its
arc-length, see, for example,~\cite{DTZ:snake}. In a sense, the
total curvature tells us how far is the curve from being a straight
line: for a straight line the total curvature is 0, for a planar
curve that coincides with a $\pi/2$-segment of the boundary of the
unit circle the total curvature is $\pi/2$, etc. Intuitively, if we
were to attempt to follow the curve numerically using, say, a
predictor-corrector type scheme, where one tries to make a predictor
step as close to the curve's tangent as possible, the total
curvature gives us some idea of how difficult it might be to
traverse such a curve. For example, one can traverse a linear
segment of 0 total curvature with just one predictor step knowing
the exact tangent, while on the opposite end of the spectrum, it
might take many steps to follow the curve of large total curvature
that makes many sharp turns.

In the context of path-following interior-point methods one
typically attempts to follow the central path that leads us to the
optimal solution, starting from the problem's analytic center. As
such, one may expect to witness many iterations of the optimization
algorithm when dealing with an LP instance where the central path is
known to have large total curvature, e.g., see~\cite{DTZ:KMcube}.

Therefore, given the above motivation, our goal is to investigate
the total curvature of the Shrink-Wrapping trajectories $d(t)$ as
compared to that of the central path, and get some feeling how the
two differ at least from the numerical perspective. We hope that the
latter would shed some light onto how efficient an algorithm based
on the Shrink-Wrapping setting might turn out. We focus on $d(t)$
rather than $x(t)=x(d(t))$ as the dynamics for $d$ defined
by~\ref{Dynamics:ODE} seems immediately suitable for defining the
corresponding discrete predictor-corrector scheme, see second
subsection of Section~\ref{SecSWIntro}, while it is not yet clear
what would be the natural setup for tracing $x(t)$ alone.

Due to the nature of the LP's considered, it is convenient to re-write
the problem in the so-called dual form
\begin{equation}\label{LP:dual}
\max_y\{f^T y:\: Gy\leq h\}.
\end{equation}
Note that this does not mean that we take the dual problem to the
$LP$ under consideration, but rather simply re-write its constraints
in the inequality form. The central path $\mathcal{P}$ corresponds
to the standard log-barrier and may be parameterized as
\begin{equation*}
\mathcal{P}=\{y\in\Re^\ell:\: y(\mu)=\arg\max_y\:\mu f^T
y+\sum_{i=1}^k \ln(h_i-G_{i,:}y)\mbox{ for some }\mu\in(0,\infty)\},
\end{equation*}
where $G_{i,:}$ is the $i^{th}$ row of $G\in\Re^{k\times\ell}$.

For the sake of comparison $\{d(t)\}_{t\geq0}$ trajectories,
developed for the $LP$ in standard equality form, are mapped to the
Shrink-Wrapping trajectories $\mathcal{D}$ in the same space of
$y$-variables as $\mathcal{P}$; since the transformation of $d(t)$
is affine, it does not change the qualitative nature of our
conclusions. If we report $x(t)$ or $x(d)$, we allow for a slight
abuse of notation and use same symbols for equivalent points in
$y$-basis. Both $\mathcal{P}$ and $\mathcal{D}$ are started at the
analytic center
\begin{equation*}
\chi=\arg\max_y\:\sum_{i=1}^k \ln(h_i-G_{i,:}y).
\end{equation*}
Basically, we aim to understand which of the two, the central path
$\mathcal{P}$ or the Shrink-Wrapping trajectory $\mathcal{D}$,
appear to be more straight. The presented findings are mostly
numerical and only suggest certain conclusions. Although, the
subsequent exposition is fairly lengthy, we believe that is is
important to provide enough details for the numerical experiments to
be repeated by the reader independently from us, if desired.

\subsection{Megiddo-Shub simplex}

For sufficiently small $\varepsilon>0$, $LP$ may be formulated as follows
\begin{equation*}
\min_x\{c^T x:\: \mathbf{1}^Tx=1, x\in\Re^n_+\}
\end{equation*}
where
\[
c_i=\left\{\begin{array}{lc}
-(1+\varepsilon)^{i-1} & , i<n,\\
0 & , i=n,
\end{array}\right.
\]
and re-written in the dual or inequality form with $f=-c_{-n}$ and
$y\in\Re^{n-1}$ as
\begin{equation*}
\max_y\{f^T y:\: [\mathbf{1}^T; -I] y\leq [1;\mathbf{0}]\}.
\end{equation*}

\begin{figure}[h!]
\begin{center}
\epsfig{file=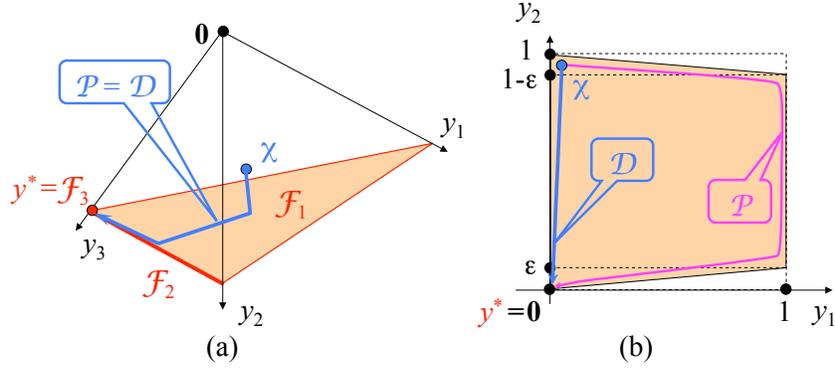,height=10.0cm} \caption{Shrink-Wrapping
dynamics close-up}\label{Hyper:Fig6}
\end{center}
\end{figure}

With $e^{(j)}$ denoting the $j^{th}$ unit vector, i.e.,
$e^{(j)}_j=1$ and $e^{(j)}_{-j}=\mathbf{0}$,
the optimal solution is $y^*=e^{(n-1)}\in\Re^{n-1}$. 
For $\varepsilon$ small enough, the central path $\mathcal{P}$ is
known to make $n-2$ sharp nearly-$\pi/2$ turns. Let
$\mathcal{F}_j=\conv\{e^{(j)}, e^{(j+1)},\ldots, e^{(n-1)}\}$ denote
$(n-j-1)$-dimensional face of $\{y\in\Re^{n-1}:\: [\mathbf{1}^T; -I]
y\leq [1;\mathbf{0}]\}$ spanned by $e^{(j)}, \ldots, e^{(n-1)}$. The
path starts at the analytic center $\chi=\frac{1}{n}\mathbf{1}$ and
first proceeds nearly orthogonal to the face $\mathcal{F}_1$. Next,
the path moves almost inside $\mathcal{F}_1$ and nearly orthogonal
to $\mathcal{F}_2$, until it nearly reaches $\mathcal{F}_3$, at
which point the path again makes a nearly-$\pi/2$ turn towards the
next face $\mathcal{F}_4$, and so on, until $\mathcal{P}$ reaches
$y^*=\mathcal{F}_{n-1}$; see Figure~\ref{Hyper:Fig6}(a).
Respectively, the total curvature of $\mathcal{P}$ is of order
$n\pi/2$; the lower bound may be established using the technique
of~\cite{DTZ:KMcube}, the upper $O(n)$ bound on the total curvature
follows from the bound on the so-called average total curvature of
$\mathcal{P}$ established in~\cite{DMS:aveCurv}.

The corresponding hyperbolic relaxation $HP_{r,d}$ is a convex
quadratic optimization problem, that is, since $m=1$ we have
$r=n-2$, and so the boundary of $HP_{r,d}$ is characterized by
$E_2(x./d)=0$, and thus $x(d)$ may be computed explicitly.

Renegar has observed that in case of $m=1, r=n-2$, the
Shrink-Wrapping trajectory for the $LP$ in equality form, started at
a point on the central path, coincides with the portion of
$\mathcal{P}$ from that point on, namely, it can be shown that with
$d\in\mathcal{P}$, $\dot{d}$ as in~\ref{Dynamics:ODE} produce a
direction tangential to $\mathcal{P}$. From the characterization of
the central path it follows that a tangent vector $d'$ to
$\mathcal{P}$ at a point $d$ is given by
\begin{eqnarray}\label{CP:tangent}
\Diag(\mathbf{1}./d.^2)\:d'=-\dot{\mu}\:c+A^Tu,\:\:Ad'=\mathbf{0},
\end{eqnarray}
for some $0\neq\dot{\mu}\in\Re, u\in\Re^m$; moreover, $\dot{\mu}>0$
corresponds to the direction of increasing $\mu$, that is, improving
the objective value along $\mathcal{P}$. At the same time, the first
of KKT conditions~\ref{HP:KKT_p} for $x=x(d)$ implies
\[
x=\Diag(d.^2)\:(\tau c+A^Tv),\:\tau>0, v\in\Re^m,
\]
while the second condition in~\ref{HP:KKT_p} used to determine the
precise value for $\tau$ for now may be ignored. One may verify that
the direction $d'=\dot{d}=x(d)-d$ indeed solves~\ref{CP:tangent}:
\begin{eqnarray*}
\begin{array}{cl}
\Diag(\mathbf{1}./d.^2)\:d' & =\Diag(\mathbf{1}./d.^2)\:
\left(\Diag(d.^2)\:(\tau c+A^Tv)-d\right) \\
& =\tau c+A^Tv-\mathbf{1}./d = (\tau-\mu)c+A^T(v+w),
\end{array}
\end{eqnarray*}
noting that $d\in\mathcal{P}$ implies $-\mathbf{1}./d=-\mu\:c+A^Tw,
\mu>0$. Note that $\tau\neq\mu$, as if it was, we could write $\tau
c=\mu c=\mathbf{1}./d-A^Tw$ and so $x./d.^2=\tau
c+A^Tv=\mathbf{1}./d+A^T(v-w)$ resulting in $x./d.^2-\mathbf{1}./d=
\Diag(\mathbf{1}./d.^2)(x-d)=A^T(v-w)$, that is, we could scale the
vector $x-d$, which belongs to the null space of $A$, by
pre-multiplying it with a positive-definite matrix
$\Diag(\mathbf{1}./d.^2)$, and obtain a vector in the range space of
$A^T$, $A^T(v-w)$, which is impossible as the null space of $A$ and
the range of $A^T$ are orthogonal subspaces of $\Re^n$. Furthermore,
since the $LP$ objective is monotone along both $\mathcal{P}$ and
$\mathcal{D}$, we must have $\tau<\mu$. Lastly,
$A\dot{d}=\mathbf{0}$ is trivial.

So, in this case, $\mathcal{P}=\mathcal{D}$, see
Figure~\ref{Hyper:Fig6}(a), and consequently the Shrink-Wrapping
trajectory is bound to have relatively large total curvature on the
order of $n$; note that both the dimension of the ambient space
$\ell=n-1$ containing the feasible region of the inequality-form
problem and the number of corresponding inequality constraints $k=n$
are almost the same.

\subsection{DTZ snake}

For this and the next subsection it is more natural to describe the
optimization problem in its dual form~\ref{LP:dual}. Since the
Shrink-Wrapping trajectories were developed for LP in standard
equality form, we start by describing the equivalent transformation
between the two formulations. Namely, given~\ref{LP:dual} we explain
how to formulate the \emph{equivalent} $LP$, equivalent in a sense
that any feasible point of~\ref{LP:dual} is uniquely mapped into
$LP$-feasible point and vice-versa, including the optimal solutions.
For simplicity we assume $G\in\Re^{k\times\ell}, k>\ell,$ to be
full-rank.

Observe that $Gy\leq h$ may be re-written as
\begin{eqnarray}\label{LP:y2x}
x=h-Gy,\: x\in\Re^k_+,
\end{eqnarray}
that is, $(h-x)$ belongs to the column-space of $G$ for some
nonnegative $x$. Let rows of $A\in\Re^{(k-\ell)\times k}$ form a
basis of the null space of columns of $G$, then
\[
A(h-x)=AGy=\mathbf{0}\in\Re^{k-\ell}
\]
and thus $Gy\leq h$ may be re-written as
\[
Ax=b, x\in\Re^k_+,
\]
where $b=Ah$. Likewise, from $h-x=Gy$, given $x$ we may easily
recover $y$ by pre-multiplying both sides with $G^T$:
\begin{eqnarray}\label{LP:x2y}
y=(G^TG)^{-1}G^T (h-x),
\end{eqnarray}
and so minimizing $-f^Ty$ corresponds to minimizing
$-f^T(G^TG)^{-1}G^T (h-x)=-c^Th+c^Tx$ with $c=G(G^TG)^{-T}f$; note
that $c^Th$ is a constant term that does not depend on $x$.
Therefore,~\ref{LP:dual} may be re-written as $LP$ with $A,b,c$ as
above and $m=k-\ell, n=k$.

The detailed DTZ snake construction in inequality form and the
subsequent analysis of the central path's geometry may be found
in~\cite{DTZ:snake}. The equivalent $LP$ may be constructed
according to the procedure above. For illustration purposes we
consider the case of $k=6, \ell=2$, in which case we have
$f=(0,-1)$, 
\[
\begin{array}{c}
G_{1,1}=0, G_{1,2}=1, G_{2,1}=1, G_{2,2}=-1/10, G_{3,1}=-1,
G_{3,2}=-1/3,\\
G_{i,1}=(-1)^i, G_{i,2}=-\frac{10^{i-2}}{11},\:i\geq4,
\end{array}
\]
\[
\begin{array}{c}
h_1=1, h_2=1/2, h_3=1/3,\\
h_i=\frac{5}{11}-\frac{10^{-4}i}{k^2},\:i\geq4.
\end{array}
\]
and $m=4, n=6$ for the equivalent $LP$.

\begin{figure}[h!]
\begin{center}
\epsfig{file=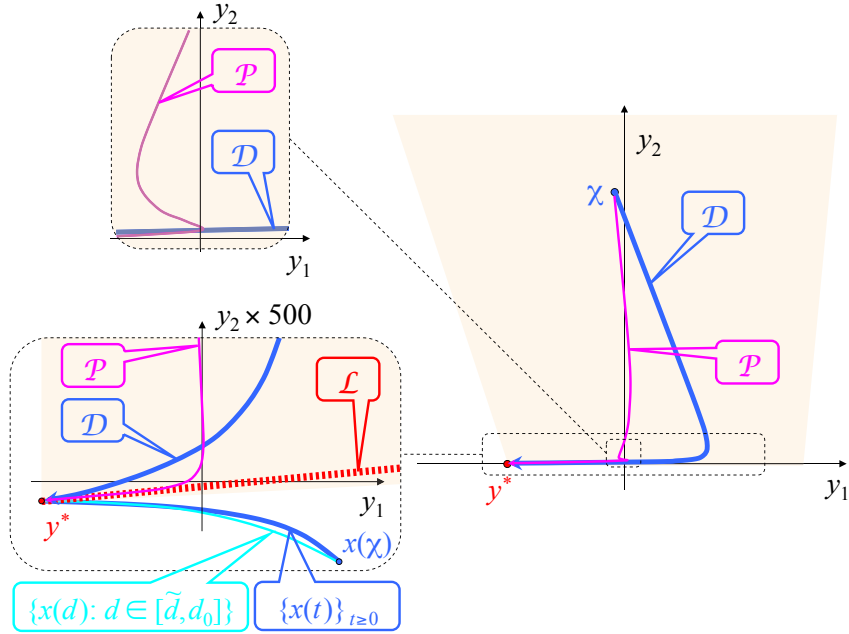,height=10.0cm} \caption{Dynamics close-up for
DTZ snake}\label{Hyper:Fig7}
\end{center}
\end{figure}

For $k$-even, $y^*$ is given by the intersection of the third and
last inequality producing
\[
-f^Ty^*=y^*_2=\frac{\frac{26}{33}-\frac{10^{-4}}{k}}{-\frac{10^{k-2}}{11}-\frac{1}{3}},\:
y^*_2=-\frac{1}{3}(y^*_2+1);
\]
for $k$-odd, the solution corresponds to the intersection of
$k^{th}$ and $(k-1)^{th}$ inequality,
\[
-f^Ty^*=y^*_2=\frac{1}{10^{k-3}}\left(-\frac{10}{11}+\frac{10^{-4}}{k^2}(2k-1)\right),\:
y^*_2=-\frac{10^{k-2}}{11}y^*_2-\frac{5}{11}+\frac{10^{-4}}{k}.
\]
The optimal $x^*$ may be computed according to~\ref{LP:y2x}; note in
case $k=6, \ell=2$, $x^*$ has four basic and two non-basic
variables, namely, $x^*_3=x^*_6=0$. As $k\rightarrow\infty$, the
central path $\mathcal{P}$ is known to make almost $k$ nearly-$\pi$
sharp turns, see Figure~\ref{Hyper:Fig7}. Respectively, the total
curvature of $\mathcal{P}$ is at least of order $k$. As mentioned
in~\cite{DTZ:snake}, the construction may be easily generalized to
arbitrary $\ell$; also, DTZ-snake may be modified to make all the
constraints non-redundant.

The corresponding hyperbolic relaxation $HP_{r,d}$ corresponds to
the first hyperbolic derivative cone of $\Re^n_{++}$, that is, $r=1$
and the boundary of $HP_{r,d}$ corresponds to $E_{n-1}(x./d)=0$, so
no explicit formula for $x(d)$ seems likely to exist.

Unfortunately, for DTZ-snake construction we could not establish an
analytic relationship between the total curvature of $\mathcal{P}$
and $\mathcal{D}$ unlike for the case of Megiddo-Shub simplex.
Instead, here we resort to numerics.

First, we describe our computational methodology for recovering
$\mathcal{P}$. A seemingly natural choice would be to use a
short-step path-following interior-point method, see, for
example~\cite{Renegar:IPM}. However, in our computational
experiments we observed that this approach suffers heavily from
numerical errors as the iterates approach the optimum,
in part, due to inherent ill-conditioning and the large bit-input
size of $F$. In turn, this causes significant problems while
attempting to recover $\mathcal{P}$ as for DTZ-snake the central
path starts to exhibit its pathological behavior only very close to
$y^*$, where the short-step method would typically fail due to
round-off errors.

The numerical stability problem is resolved by re-parameterizing
$\mathcal{P}$ with level sets of $f^Ty$:
$y(\nu)=\arg\max_{y:f^Ty=\nu}\:\sum_{i=1}^k \ln(h_i-G_{i,:}y)$,
which results in univariate maximization problem on a fixed
interval. The latter one-dimensional optimization problem for
finding the point $y(\nu)$ is handled with a simple bi-section
method, thus, avoiding the ill-conditioning problems associated with
the second derivative-based Newton's method; as a stopping criteria
for the bi-section scheme we use the length of the interval
containing $y(\nu)$ falling below a prescribed threshold. Since we
are interested in recovering the geometry of the path, the length of
the interval measured with respect to the Euclidian norm appears to
give us more accurate answer when approximating $y(\nu)$, as opposed
to working with the norm induced by the self-concordant barrier
typically used in the interior-point methods. The reason 
is the lack of scaling along any particular direction for the
Euclidean norm unlike for the barrier-induced norm.

In order to traverse $\mathcal{P}$, we gradually increase the
corresponding $\nu$ parameter starting from the value $f^T\chi$, and
generate a sequence of iterates $\{y_i\}_{i=0,K}$ in close Euclidian
proximity to the central path, until we reach the LP optimum.
Furthermore, to speed up computations of each subsequent $y_i$, we
warm-start the bi-section from $y_{i-1}$. Near the optimum the
discrete stepping of $\nu$ gets more and more refined to allow us to
capture sharp turns of $\mathcal{P}$. The first iterate $y_0$
corresponds to the approximate analytic center $\chi$, which is
computed using MATLAB `fsolve' routine: we attempt to zero out the
gradient of the log-barrier, starting from the initial approximation
that corresponds to the analytical center
$\left(0,\frac{k-3}{k-4}\right)$ of the perturbed LP with a feasible
region corresponding to
$\left\{y\in\Re^2:\:\left[I;-I;-{e^{(2)}}^T;\ldots;-{e^{(2)}}^T\right]y\leq(1,1,1,0,0,\ldots,0)\right\}$,
i.e., a planar unit cube centered at $(1/2, 1/2)$ with the bottom
face repeated $k-3$ times. The last iterate in the sequence is
$y_K=y^*$; $K$ is chosen so that $y_{K-1}$ is close enough to $y^*$.
The resulting \emph{approximate central path}
$\widetilde{\mathcal{P}}$ is a piece-wise linear interpolation of
$\mathcal{P}$ from $\{y_i\}_{i=0,K}$.

Next, we describe our computational methodology for recovering
$\mathcal{D}$. We compute the approximate Shrink-Wrapping trajectory
for $LP$ and map both $d(t)$ and $x(t)$ onto the feasible region
of~\ref{LP:dual} according to~\ref{LP:x2y}. To recover $d(t), x(t)$,
we employ standard discrete predictor-corrector scheme for tracing
the trajectory of the ODE given by~\ref{Dynamics:ODE}: given some
initial pair $(d_i, x_i), x_i\approx x(d_i)$, we set the next
iterate $d_{i+1}=d_i+\alpha (x_i-d_i)$ and $x_{i+1}\approx
x(d_{i+1})$, where $\alpha>0$ is some small constant. The
predictor-corrector scheme is known to converge to the true ODE
trajectory when $\alpha\rightarrow 0$ under some mild assumptions.
We experimented with several choices of $\alpha$. We found that the
most numerically stable approach is to normalize the predictor step
length along $(x_i-d_i)$ to have a prescribed length $\tilde\alpha$,
where $\tilde\alpha$ is either fixed on the order of
$10^{-2}-10^{-3}$, or is gradually decreasing as $d_i$ approach
$x^*$ to enforce $LP$ feasibility of $d_i$. Both step normalization
choices appear to attain virtually indistinguishable numerical
results. We generate the sequence $\{d_i\}_{i=0,K}$ with $d_0$
approximating the analytic center of $LP$ and $d_K=x^*$, similar to
the case of $\mathcal{P}$. The \emph{approximate Shrink-Wrapping
trajectory} $\widetilde{\mathcal{D}}$ is a piece-wise linear
interpolation from $\{d_i\}_{i=1,K}$.

In order to compute $x_{i+1}\approx x(d_{i+1})$ we use Newton's
method to find the root of $f(\xi)$ as defined in the previous
section, warm-started at $x_i$, with termination criteria being the
Euclidian norm of the gradient of $f(\xi)$ falling below a certain
threshold. The corresponding function evaluations and derivative
information may be computed using the FFT approach outlined
in~\cite{Renegar:Poly}. Given $d_0$ --equivalently,
$y_0\approx\chi$-- we recover the initial point $x_0\approx x(d_0)$
by performing a linear homotopy from another point on the central
line. That is, we numerically follow $x(d)$ using Newton's method as
$d$ traverses $[\widetilde{d},d_0]$, starting at
$\widetilde{d}\in\mathcal{L}$. Recall that at least in the vicinity
of $\mathcal{L}$ the Newton iterates are well defined. As $d$ gets
gradually changed from $\widetilde{d}$ to $d_0$, MATLAB does not
encounter any problems with ill-conditioning or non-invertibility of
derivative matrices. The latter and the homotopy path $x(d),
d\in[\widetilde{d},d_0]$ appearing rather smooth, see
Figure~\ref{Hyper:Fig7}, indicates that we indeed did not switch
branches of $f(\xi)$ and recovered the correct approximate to
$x(d_0)$; if desired, we may further confirm the validity of our
approximation by checking that $x(d)$ is in or close enough to
$\mathcal{K}_{r,d}$.

Note that unlike the central path iterates $y_i$, we do not use any
low-order method to recover $x_i$ because there appears to be no
suitable re-parametrization of $x(d)$ readily available. Thus,
hypothetically, our computations for $\mathcal{D}$ are more
susceptible to round-off errors. However, we are still fairly
confident in the results of our numerical findings due to the
following two reasons.

\begin{itemize}
\item Computational safeguard procedure: to make sure our numerical
approach produces no obvious nonsense results, we re-compute an
approximate central path relying on the equivalence of $\mathcal{P}$
and $\mathcal{D}$ for $r=n-2$, using the outlined numerical approach
for computing $\widetilde{\mathcal{D}}$ with $r=n-2$ as above. We
compare our results with the first approximation
$\widetilde{\mathcal{P}}$ to make sure both paths are consistent
with one another. Indeed, both methods seem to recover visually
indistinguishable approximate central paths. Moreover, we are not
overly concerned with approximating $x(d_i)$ with $x_i$ due to the
fact that numerical ill-conditioning of Newton's method, resulting
from the ill-conditioning of the derivative matrix, as reported by
MATLAB, manifests itself for the iterates $d_i$ only well past the
last sharp turn of the central path with respect to the LP objective
value. That is, by the time MATLAB begins to report the numerical
ill-conditioning for locating $x_i\approx x(d_i)$, the respective
points on the central path that correspond to the LP objective level
sets with values $c^T d_i$ are located well past the last sharp turn
of $\widetilde{\mathcal{P}}$.
\item Central line: since the transformation~\ref{LP:x2y} is linear,
the existence and the attractor-like properties of the invariant
central line persist through the equivalent transformation between
the LP formulations. In particular, for our example, in the basis of
$y$-variables the central line is extends from
$y^*\approx(-0.333,-0.000866)$ and passes through a point with
approximate coordinates $(.0528,-.0000171)$ -- coincidentally, the
point which we start the linear homotopy from to recover $x_0$, see
Figure~\ref{Hyper:Fig7}; the procedure to recover $\mathcal{L}$
given $x^*$ is outlined in Proposition~\ref{CL:exists}. By
Theorem~\ref{Dynamics:MainTheorem}, if the trajectory $\mathcal{D}$
at some point gets sufficiently close to the central line, from that
point on $\mathcal{D}$ gets pulled into the line very quickly, and,
most certainly, the central line may not be crossed over. Examining
our numerical results we see that indeed $\widetilde{\mathcal{D}}$
appears to get very close to the central line and straightens out
from that point on, see Figure~\ref{Hyper:Fig7}, concurring to our
intuition.
\end{itemize}
As the purpose of this section is mostly to gain some qualitative
insight into the behavior of Shrink-Wrapping trajectories
comparative to the central path, we do not attempt to further refine
or justify our numerical approach for tracing $\mathcal{D}$.

In summary, for DTZ construction, the Shrink-Wrapping trajectory
does not seem to exhibit any pathological behavior as compared to
the central path; both $d(t)$ and $x(t)$ trajectories appear to be
fairly straight and thus are likely to have small total curvature,
with $d(t)$ making only one turn. Unlike the case of $\mathcal{P}$,
where the source of the large total curvature is the constant
zigzagging of the central path, similar behavior for $\mathcal{D}$
is less likely due to the existence of the central line.

\subsection{Redundant Klee-Minty cube}

The detailed problem formulation in inequality form and the
subsequent analysis of the central path's geometry may be found
in~\cite{DTZ:KMcube}. With properly chosen parameters the central
path is known to make at least $2^\ell-2$ sharp nearly-$\pi/2$ turns
closely following the standard simplex method pivot sequence,
resulting in large total curvature of $\mathcal{P}$; $k$ is
exponential as a function of $\ell$. In particular, we use the
\emph{geometrically-decaying distance} model with $k=O\left(\ell^3
2^{2\ell}\right)$ and rely on Corollary 7.2 of~\cite{DTZ:KMcube}
that guarantees exponential order of the total curvature of
$\mathcal{P}$. For our illustration we consider planar redundant
Klee-Minty cube with $\ell=2$ and $k\approx 14,000$, given by
$f=(0,1)$,
\[
\begin{array}{c}
G_{1,1}=-1, G_{1,2}=0, G_{2,1}=1, G_{2,2}=0, G_{3,1}=\varepsilon,
G_{3,2}=-1, G_{4,1}=-\varepsilon, G_{4,2}=1,\\
G_{i,1}=-1, G_{i,2}=0, \:i\in[5,4+h_1],\\
G_{i,1}=\varepsilon, G_{i,2}=-1, \:i\in[5+h_1,4+h_1+h_2],
\end{array}
\]
\[
\begin{array}{c}
h_1=0, h_2=1, h_3=0, h_4=1,\\
h_i=r_1, \:i\in[5,4+h_1],\\
h_i=r_2, \:i\in[5+h_1,4+h_1+h_2],
\end{array}
\]
where $\varepsilon=.1$, $\delta=.05$ and $r_1=8, r_2=4, h_1=2963,
h_2=10766$ are computed according to~\cite{DTZ:KMcube}. Clearly, the
optimal solution is $y^*=\mathbf{0}$. We intentionally reduce
$\varepsilon, \delta$ to get sharper turns of $\mathcal{P}$: the
central path makes two sharp nearly-$\pi/2$ turns near vertices
$(1,1-2\varepsilon)$ and $(1,\varepsilon)$, see
Figure~\ref{Hyper:Fig6}(b). The corresponding $LP$ has $m=13731,
n=13733$.

The corresponding hyperbolic relaxation $HP_{r,d}$ corresponds to
the first hyperbolic derivative cone of $\Re^n_{++}$ with $r=1$.
Similar to DTZ construction, we cannot analyze the setting
analytically and resort to numerics: most of the numerical
considerations above may be carried over to the case of redundant
Klee-Minty construction. We implement several changes to better
address the nature of the problem.
\begin{itemize}
\item For the construction, the short-step path-following
interior-point method produces stable numerical results which are
consistent with theoretical findings in~\cite{DTZ:KMcube},
therefore, the method may be used to recover
$\widetilde{\mathcal{P}}$ for comparative purposes.
\item When recovering $\widetilde{\mathcal{D}}$, it is much more
efficient to re-cast Newton's method for tracing $x(d)$ into the
basis of $y$-variables, which gives us much smaller, and thus, less
prone to numerical errors, linear system that we need to work with.
\item Lastly, due to high degree of the $HP_{r,d}$-underlying hyperbolic
polynomial, the FFT approach seems not as effective as for DTZ
snake, mostly due to round-off errors. Instead, we use the ratio
$\frac{E_{n-1}(x./d)}{E_n(x./d)}$ to characterize the boundary of
$\mathcal{K}_{r,d}$ outside of $\partial\Re^n_+$. The latter allows
for an explicit and simple form of the derivatives needed to
implement Newton's method; also, see the subsequent discussion.
\end{itemize}

Instead of directly computing the central line, which in this case
seems to exhibit its attractor properties only towards the very end
of the Shrink-Wrapping trajectory, we claim that $HP_{r,d}$ produces
a very tight relaxation to $LP$ itself due to the high degree of the
underlying hyperbolic polynomial and the presence of many
remotely-positioned redundant constraints. In other words, for any
$LP$-strictly feasible $d$, $x(d)$ may not be far from $x^*$;
consequently, $d(t)$ gets driven to $x^*$ almost along the straight
line, see Figure~\ref{Hyper:Fig6}(b). For brevity we only sketch the
argument.

It is convenient to switch back and forth between the domains
of~\ref{LP:dual} and $LP$; to this end we introduce three pairs of
vector variables in $\Re^\ell$ and $\Re^n$ spaces respectively,
$y\equiv d, \gamma\equiv\Delta, y(t)=y+\gamma t\equiv x(t)$ with
$t\in\Re$, where $\equiv$ is the equivalence relationship given
by~\ref{LP:y2x},~\ref{LP:x2y}. We allow for a slight abuse of
notation referring with $HP_{r,d}$ to both the primal and $y$-space
re-formulation of the hyperbolic relaxation problem.

Fix $y$ in the interior of Klee-Minty cube defined by the first four
constraints of $Gy\leq h$; denote the corresponding 1-dimensional
faces --hyperplanes-- by $\mathcal{F}_1, \mathcal{F}_2,
\mathcal{F}_3, \mathcal{F}_4$. In order to understand how close is
$HP_{r,d}$ to $LP$, consider when
\[
q(t)=\frac{E_{n-1}(x(t)./d)}{E_n(x(t)./d)}=\frac{d_1}{x_1(t)}+\frac{d_2}{x_2(t)}+\cdots+\frac{d_n}{x_n(t)}
\]
crosses 0, that is, when $E_{n-1}(x(t)./d)=0$; recall that the
feasible region of $HP_{r,d}$ touches $LP$-feasible region precisely
at the vertices of Klee-Minty cube. In other words, we ask how far
along the ray $y(t)=y+\gamma t$ one needs to travel outside of the
Klee-Minty cube before encountering the boundary of
$HP_{r,d}$-feasible region.

A simple root $t$ of $E_{n-1}(x(t)./d)=0$ in the vicinity of the
boundary of Klee-Minty cube may occur only past the point when
$y(t)$ crosses either $\mathcal{F}_1, \mathcal{F}_2, \mathcal{F}_3,
\mathcal{F}_4$, and only one of these faces at a time, i.e., past
$t^*$ when only one of the corresponding $x_1(t^*), x_2(t^*),
x_3(t^*), x_4(t^*)$ becomes zero. By the root interlacing property
of polynomials with all real roots applied to $E_n, E_{n-1}$ we know
that $E_n, E_{n-1}>0$ inside $LP$-feasible region, and $E_n<0,
E_{n-1}>0$ just on the outside. So
\[
q(t)=\frac{d_1}{x_1(t)}+\frac{d_2}{x_2(t)}+\cdots+\frac{d_n}{x_n(t)}<0
\]
just outside of $LP$-feasible region; in fact, $q(t)\rightarrow
-\infty$ as $t\rightarrow -\frac{d_1}{\Delta_1}$,
$-\frac{d_2}{\Delta_2}$, $-\frac{d_3}{\Delta_3}$ or
$-\frac{d_4}{\Delta_4}$ from the left, while $y(t)$ remains feasible
with respect to the remaining three faces. Re-writing
\[
q(t)=\sum_{i=1}^4 \frac{d_i}{d_i+\Delta_i t}+\sum_{i=5}^n
\frac{d_i}{d_i+\Delta_i t}
\]
we note that for the second summand, recalling the redundant
constraints, we have
\[
M=h_1\frac{8}{8+(1+2\cdot.1)}+h_2\frac{4}{4+(1+2\cdot.1)}
<\sum_{i=5}^n
\frac{d_i}{d_i+\Delta_i t}
\]
for all $y(t)$ within $.1$ or lesser Euclidian distance from
Klee-Minty cube.

Since only one $\mathcal{F}_i,\: i=1,4$, is being crossed-over, say,
$\mathcal{F}_1$, we can write
\[
0\leq q(t)\leq \frac{d_1}{d_1+\Delta_1 t}+M
\]
as long as
\[
t\geq -\frac{d_1}{\Delta_1}\left(1+\frac{1}{M}\right),
\]
noting that $t$ is such that $d_1+\Delta_1 t<0$ and assuming the
remaining $\mathcal{F}_2, \mathcal{F}_3, \mathcal{F}_4$ feasible --
note that two faces of the cube may not be consecutively traversed
by $y(t)$ without encountering a root of $q(t)$ -- a root of
$E_{n-1}$. So, $q(t)$ changes sign from minus to plus as $t$ goes
through $\left(-\frac{d_1}{\Delta_1}\left(1+\frac{1}{M}\right),
-\frac{d_1}{\Delta_1}\right)$. So, from any $LP$-strictly feasible
$y$ along any $\gamma$ one needs to traverse precisely
$t=-\frac{d_1}{\Delta_1}$ to get to the boundary of $LP$-feasible
region, and at most
$t=-\frac{d_1}{\Delta_1}\left(1+\frac{1}{M}\right)$ to reach the
boundary of $HP_{r,d}$-feasible region.

In other words, $HP_{r,d}$-feasible region is a ``slightly inflated
version'' of Klee-Minty cube, ``inflated'' by a factor of at most
$1+\frac{1}{M}\approx1.0001$, so, indeed $y(t)$ remains in close
proximity to the cube, particularly, is not farther than $.1$ in
Euclidian distance. Moreover, recall that $HP_{r,d}$ is convex and
touches $LP$-feasible region only at the vertices of Klee-Minty
cube. Simple geometric considerations may be used to complete the
claim, resulting in $x(d)\approx x^*$. Note that the ratio $q$ is
the reciprocal of the concave ratio functional briefly mentioned in
Section~\ref{SecSWIntro}.

Going back to our numerical results, again, unlike $\mathcal{P}$,
$\mathcal{D}$ appears to be much more straight which suggests it
having much lower total curvature than the central path, see
Figure~\ref{Hyper:Fig6}(b). Note that the redundant constraints
--the source of large total curvature for the central path-- are
handled exceptionally well in the Shrink-Wrapping setting, as
effectively they do not play any negative role in determining the
dynamics of $d(t)$. In fact, in this example, the presence of the
remotely-positioned redundant constraints does the opposite and
helps to straighten out $\mathcal{D}$. Also, recall that the central
line, which in turn seems to drive the limiting behavior of $d(t)$,
is defined only by the active constraints, so redundancy does not
negatively affect us here either.

\section{Conclusion}
Following the idea of Renegar, we introduce the Shrink-Wrapping
setting for solving linear programming problems based on hyperbolic
relaxations of the nonnegative orthant. We analyze the local
behavior of the Shrink-Wrapping trajectories that lead to the LP
optimum, provided a suitable choice of the initial point. A striking
difference between the standard path-following interior-point
methods and the Shrink-Wrapping setting is the existence of the
invariant with respect to dynamics of the trajectory set --the
central line-- in the latter case. The central line acts, at least
locally, as an attractor set for the Shrink-Wrapping trajectories,
which in turn guarantees extremely quick, i.e., $R$-super-quadratic
local convergence of a simple bi-section type discretization scheme
based on Shrink-Wrapping.

We attempt to analyze the behavior of the Shrink-Wrapping
trajectories comparative to the central path for three known
pathological linear programming instances, where the central path
has large total curvature. Partial theoretical analysis is
substantiated with numerics. Although, we encounter a negative
example when the Shrink-Wrapping trajectory and the central path
look identical, in most cases (2 out of 3) Shrink-Wrapping
trajectories appear to be much more straight than the central paths.
This suggests that the Shrink-Wrapping approach may result in more
efficient predictor-corrector type algorithms for solving the
underlying optimization problems.

A possible explanation to this distinctive difference between the
behavior of the central path and the Shrink-Wrapping trajectories
lies in a seemingly more appropriate choice of the degree of the
hyperbolic relaxation problem. While one may think of the central
path as the Shrink-Wrapping with degree fixed permanently, the
setting suggests that the degree should be chosen adaptively.
Namely, the proper choice of the relaxation degree results in the
dynamics of the trajectory being driven by a simple root of a
polynomial system of equations, giving rise to a number of favorable
properties, such as the earlier mentioned central line, while for
the central path such a root is almost always multiple. As the
choice of the degree of Shrink-Wrapping suggests much tighter fit of
the relaxation to the original problem, we expect the trajectories
to converge to the optimum sooner -- an intuition confirmed by the
numerics.

When the linear programming problem is re-written in inequality
form, optimistically, we hope that in case of Shrink-Wrapping the
total curvature of the trajectory is driven by the dimensionality of
the ambient space, rather than the number of constraints unlike for
the central path. Again, the investigated numerics support our bold
conjecture. We present the initial analysis of the newly proposed
setting. Much work remains to be done to convert these ideas into an
actual optimization algorithm.

As a side result, we provide the first, to our knowledge, proof of
convexity of hyperbolicity cones which does not rely on complex
variables.

\section*{Acknowledgement}
I would like to thank James
Renegar for his invaluable suggestions while supervising my Ph.D.
thesis, which served as a foundation for Sections 3--4 and was
supported by NSF Grant \#CCF-0430672.


\begin{thebibliography}{99}
\bibitem{Arnold:ODE}V. Arnold, \emph{Ordinary differential
equations}, The MIT Press, 1978.
%
\bibitem{Lewis:Poly}H. Bauschke, O. Guler, A. Lewis, and H. Sendov,
\emph{Hyperbolic polynomials and convex analysis}, Can. J. Math.
\textbf{53} (2001), 470--488.
%
\bibitem{BenTalN:ConvOpt}A. Ben-Tal, A. Nemirovski, \emph{Lectures
on modern convex optimization: analysis, algorithms, and engineering
applications}, MPS-SIAM, 2001.
%
\bibitem{BP:Finance}D. Bertsimas, I. Popescu, \emph{On the relation
between option and stock prices: an optimization approach}, Oper.
Res. \textbf{50} (2002), 358--374.
%
\bibitem{ShubSmale:Complexity}L. Blum, F. Cucker, M.
Shub, and S. Smale, \emph{Complexity and real computation},
Springer-Verlag, 1998.
%
\bibitem{SMY:IMRT}M. Chu, Y. Zinchenko, S. Henderson, and M.
Sharpe, \emph{Robust optimization for intensity modulated radiation
therapy treatment planning under uncertainty}, Phys. Med. Biol.
\textbf{50} (2005), 5463--5478.
%
\bibitem{Chua:Poly}C. Chua, \emph{Relating homogeneous cones and
positive definite cones via T-algebras}, SIAM J. Opt. \textbf{14}
(2003), 500--506.
%
\bibitem{Chvatal:LP}V. Chvatal, \emph{Linear programming}, W. H. Freeman
\& Company, 1983.
%
\bibitem{DMS:aveCurv} J-P. Dedieu, G. Malajovich, M. Shub, \emph{On the curvature of the central path of
linear programming theory}, Found. Comp. Math. \textbf{5} (2005),
145-—171.
%
\bibitem{DNT:KMcube} A. Deza, E. Nematollahi, T. Terlaky, \emph{How good are
interior point methods? Klee-Minty cubes tighten
iteration-complexity bounds}, Math. Pro. \textbf{113} (2008), 1--14.
%
\bibitem{DTZ:KMcube}A. Deza, T. Terlaky, Y. Zinchenko, \emph{Central path curvature and iteration-complexity for
redundant Klee-Minty cubes}, Advances in Mechanics and Mathematics
17 (D. Gao and H. Sherali eds.), Springer, 2009, pp. 223--256.
%
\bibitem{DTZ:snake}A. Deza, T. Terlaky, Y. Zinchenko, \emph{Polytopes and arrangements: diameter and curvature},
Oper. Res. Let. \textbf{36} (2008), 215--222.
%
\bibitem{Garding:Poly}L. Garding, \emph{An inequality for
hyperbolic polynomials}, J. Math. Mech. \textbf{8} (1959), 957--965.
%
\bibitem{GW:Maxcut} M. Goemans, D. Williamson, \emph{Improved approximation
algorithms for maximum cut and satisfiability problems using
semidefinite programming}, J. ACM \textbf{42} (1995), 1115--1145.
%
\bibitem{Guler:Poly}O. Guler, \emph{Hyperbolic polynomials and
interior point methods for convex programming}, Math. Oper. Res.
\textbf{22} (1997), 350--377.
%
\bibitem{Lasserre:Poly}J. Lasserre, \emph{Global optimization with
polynomials and the problem of moments}, SIAM J. Opt. \textbf{11}
(2001), 796--817.
%
\bibitem{MS:simplex}N. Megiddo, M. Shub, \emph{Boundary behavior
of interior point algorithms in linear programming}, Math. Oper.
Res. \textbf{14} (1989), 97–-146.
%
\bibitem{Renegar:IPM}J. Renegar, \emph{A mathematical view of
interior-point methods in convex optimization}, MPS-SIAM, 2001.
%
\bibitem{Renegar:Poly}J. Renegar, \emph{Hyperbolic programs, and their
derivative relaxations}, Found. Comp. Math. \textbf{6} (2006),
59--79.
%
\bibitem{Renegar:SW} J. Renegar, \emph{Central swaths (a generalization of the central
path)}, Opt. Online (June 2010).
%
\bibitem{Boyd:MC}J. Sun, S. Boyd, L. Xiao, and P.
Diaconis, \emph{The fastest mixing Markov process on a graph and a
connection to a maximum variance unfolding problem}, SIAM Rev.
\textbf{48} (2006), 681--699.
%
\bibitem{Zinchenko:ESP}Y. Zinchenko, \emph{On hyperbolicity cones
associated with elementary symmetric polynomials}, Opt. Let.
\textbf{2} (2007), 389--402.
%
\bibitem{Zinchenko:thesis} Y. Zinchenko, \emph{The local behavior of
the Shrink-Wrapping algorithm for linear programming}, Ph.D. thesis,
Cornell University, 2005.
\end{thebibliography}
\end{document}